\documentclass[a4paper]{amsart}
\usepackage[utf8]{inputenc}
\usepackage[english]{babel}
\usepackage{amsmath}
\usepackage{amsfonts}
\usepackage{amssymb}
\usepackage{amsthm}
\usepackage{wasysym}
\usepackage{amscdx}
\usepackage{graphicx}
\usepackage{comment}
\usepackage{chngcntr}
\usepackage{tikz-cd}
\usepackage{xcolor}
\usepackage{marvosym}
\usepackage{enumitem}
\usepackage{stix}
\usepackage{todonotes}

\usepackage{hyperref}
\hypersetup{colorlinks=true}
\makeatletter
\renewcommand\@mkboth[2]{\markboth{#1}{}}
\makeatother

\numberwithin{equation}{section}

\newtheorem*{defnintro}{Definition}
\newtheorem*{thmintro}{Theorem}

\newtheorem{thm}{Theorem}[section]
\newtheorem{lem}[thm]{Lemma}
\newtheorem{prop}[thm]{Proposition}
\newtheorem{cor}[thm]{Corollary}

\theoremstyle{definition}
\newtheorem{defn}[thm]{Definition}

\newtheorem{claim}{Claim}[thm]

\theoremstyle{remark}
\newtheorem{remark}[thm]{Remark}
\newtheorem{question}{Open Problem}

\DeclareMathOperator{\Ar}{Ar}
\DeclareMathOperator{\WAr}{WAr}
\DeclareMathOperator{\Kn}{Kn}

\DeclareMathOperator{\CKn}{CKn}
\DeclareMathOperator{\WCKn}{WCKn}
\DeclareMathOperator{\skr}{\precsim_{\Kn}^{<\omega}}
\DeclareMathOperator{\cskr}{\precsim_{\CKn}^{<\omega}}
\newcommand{\arc}[2]{(\bar #1, #2)}
\newcommand{\barc}[2]{(#1, #2)}

\DeclareMathOperator{\Lin}{\mathsf{Lin}}
\DeclareMathOperator{\Scat}{\mathsf{Scat}}
\DeclareMathOperator{\WO}{\mathsf{WO}}
\DeclareMathOperator{\Fin}{\mathsf{Fin}}
\DeclareMathOperator{\LO}{\mathsf{LO}}
\DeclareMathOperator{\CO}{\mathsf{CO}}

\DeclareMathOperator{\Ima}{Im}
\DeclareMathOperator{\dom}{dom}
\DeclareMathOperator{\Int}{Int}

\DeclareMathOperator{\Supp}{Supp}
\DeclareMathOperator{\Cvx}{Int}

\newcommand{\cvx}{\trianglelefteq_{\LO}}
\newcommand{\ncvx}{\ntrianglelefteq_{\LO}}
\newcommand{\cvxeq}{\mathrel{\underline{\bowtie}_{\LO}}}

\newcommand{\iso}{\cong_{\LO}}

\newcommand{\ccvx}{\trianglelefteq_{c}}
\newcommand{\nccvx}{\ntrianglelefteq_{c}}
\newcommand{\pccvx}{\trianglelefteq_{\CO}^{<\omega}}
\newcommand{\pcvx}{\mathrel{\trianglelefteq_{c}^{<\omega}}}
\newcommand{\pccvxeq}{\mathrel{\underline{\bowtie}^{<\omega}_{\CO}}}
\newcommand{\npccvx}{\mathrel{{\ntrianglelefteq}^{<\omega}_{\CO}}}

\DeclareMathOperator{\ot}{ot}

\newcommand{\z}{\zeta}
\renewcommand{\o}{\omega}
\newcommand{\g}{\gamma}
\renewcommand{\a}{\alpha}

\newcommand{\Z}{\mathbb{Z}} 
\newcommand{\Q}{\mathbb{Q}} 
\newcommand{\R}{\mathbb{R}} 
\newcommand{\N}{\mathbb{N}} 

\newcommand{\rest}{\!\restriction\!}

\newenvironment{enumerate-(a)}{\begin{enumerate}[label={\upshape (\alph*)}, leftmargin=2pc]}{\end{enumerate}}
\newenvironment{enumerate-(a)-r}{\begin{enumerate}[label={\upshape (\alph*)}, leftmargin=2pc,resume]}{\end{enumerate}}
\newenvironment{enumerate-(a)-5}{\begin{enumerate}[label={\upshape (\alph*)}, leftmargin=2pc,start=5]}{\end{enumerate}}
\newenvironment{enumerate-(A)}{\begin{enumerate}[label={\upshape (\Alph*)}, leftmargin=2pc]}{\end{enumerate}}
\newenvironment{enumerate-(A)-r}{\begin{enumerate}[label={\upshape (\Alph*)}, leftmargin=2pc,resume]}{\end{enumerate}}
\newenvironment{enumerate-(i)}{\begin{enumerate}[label={\upshape (\roman*)}, leftmargin=2pc]}{\end{enumerate}}
\newenvironment{enumerate-(i)-r}{\begin{enumerate}[label={\upshape (\roman*)}, leftmargin=2pc,resume]}{\end{enumerate}}
\newenvironment{enumerate-(I)}{\begin{enumerate}[label={\upshape (\Roman*)}, leftmargin=2pc]}{\end{enumerate}}
\newenvironment{enumerate-(I)-r}{\begin{enumerate}[label={\upshape (\Roman*)}, leftmargin=2pc,resume]}{\end{enumerate}}
\newenvironment{enumerate-(1)}{\begin{enumerate}[label={\upshape (\arabic*)}, leftmargin=2pc]}{\end{enumerate}}
\newenvironment{enumerate-(1)-r}{\begin{enumerate}[label={\upshape (\arabic*)}, leftmargin=2pc,resume]}{\end{enumerate}}
\newenvironment{itemizenew}{\begin{itemize}[leftmargin=2pc]}{\end{itemize}}
\newenvironment{enumerate-(star)}{\begin{enumerate}[label={\upshape{(\( \star_{ \arabic*} \))}}, leftmargin=2pc]}{\end{enumerate}}

\DeclareFontEncoding{LS1}{}{}
\DeclareFontSubstitution{LS1}{stix}{b}{it}
\DeclareSymbolFont{symbols4}{LS1}{stixbb}{b}{it}
\DeclareMathSymbol{\csube}{\mathrel}{symbols4}{"9F}
\DeclareFontEncoding{LS1}{}{}
\DeclareFontSubstitution{LS1}{stix}{b}{it}
\DeclareSymbolFont{symbols4}{LS1}{stixbb}{b}{it}
\DeclareMathSymbol{\csub}{\mathrel}{symbols4}{"9D}
\DeclareFontEncoding{LS1}{}{}
\DeclareFontSubstitution{LS1}{stix}{b}{it}
\DeclareSymbolFont{symbols4}{LS1}{stixbb}{b}{it}
\DeclareMathSymbol{\csupe}{\mathrel}{symbols4}{"A0}
\DeclareFontEncoding{LS1}{}{}
\DeclareFontSubstitution{LS1}{stix}{b}{it}
\DeclareSymbolFont{symbols4}{LS1}{stixbb}{b}{it}
\DeclareMathSymbol{\csup}{\mathrel}{symbols4}{"9E}

\begin{document}

\title[Convex Embeddability and Knots]{Convex Embeddability 
and Knot Theory}
\subjclass[2020]{03E15, 06A05, 57K10, 57M30}
\author[M.~Iannella]{Martina Iannella} \address{Dipartimento di scienze matematiche, informatiche e fisiche, Universit\`a di Udine, Via delle Scienze 208, 33100 Udine --- Italy} \email{martinaiannella2@gmail.com}
\author[A.~Marcone]{Alberto Marcone} \address{Dipartimento di scienze matematiche, informatiche e fisiche, Universit\`a di Udine, Via delle Scienze 208, 33100 Udine --- Italy} \email{alberto.marcone@uniud.it} 
\author[L.~Motto Ros]{Luca Motto Ros} \address{Dipartimento di matematica \guillemotleft{Giuseppe Peano}\guillemotright, Universit\`a di Torino, Via Carlo Alberto   10, 10121 Torino --- Italy} \email{luca.mottoros@unito.it}
\author[V.~Weinstein]{Vadim Weinstein} 
\address{Center for Ubiquitous Computing\\
Erkki Koiso-Kanttilan katu 3\\ 
door E P.O Box 4500\\
FI-90014 University of Oulu} 
\email{vadim.weinstein@iki.fi}

\thanks{The first three authors were partially supported by the Italian PRIN 2017 Grant ``Mathematical Logic: models, sets, computability".}

\begin{abstract}
  We consider countable linear orders and study the quasi-order of convex embeddability and its induced equivalence relation. 
  We obtain both combinatorial and descriptive set-theoretic results, and further extend our research to the case of circular orders. 
  These results are then applied to the study of arcs and knots, establishing combinatorial properties and lower bounds (in terms of Borel reducibility) 
  for the complexity of some natural relations between these geometrical objects.
\end{abstract}

\maketitle

\tableofcontents

\section{Introduction}

Knots are very familiar and tangible objects in everyday life, and they
also play an important role in modern mathematics. A mathematical knot
is a homeomorphic copy of $S^1$ embedded in $S^3$. The study of knots
and their properties is known as knot theory (see
e.g.~\cite{BZ03}). This paper uses discrete objects, such as linear
and circular orders, to gain insight into knots. This approach was
already exploited in~\cite{Kul17}, where it is shown
that isomorphism $\iso$ on the Polish space $\LO$ of
linear orders defined on \(\N\) strictly Borel reduces to equivalence
on knots \( \equiv_{\Kn} \). (Recall that Borel reducibility $\leq_B$ provides a hierarchy
of complexities for equivalence relations defined on Polish or
standard Borel spaces.)

The proof in \cite{Kul17} uses proper arcs (which intuitively are
obtained by cutting a knot) and their subarcs, called ``components''
in \cite{Kul17}, which are the analogues of convex
subsets of linear orders. Thus, to expand the previous results it is
natural to study the following relation between linear orders.

\begin{defnintro}
Given linear orders \(L \) and \( L'\), we set \(L \trianglelefteq L'\) if
and only if \(L\) is isomorphic to a convex subset of \(L'\).
\end{defnintro}

We call convex embeddability the relation \( \trianglelefteq \), which was already introduced and briefly studied in \cite{BCP73}. Even if convex embeddability is a very natural relation, as far as we know it has not received much attention in the last 50 years.

We first focus on the restriction of \(\trianglelefteq\) to \(\LO\), denoted by \(\cvx\).  We begin establishing that \(\cvx\) induces a structure on \(\LO\) very different from the one given by the usual embeddability relation: as conjectured by Fra\"{i}ss\'{e} in 1948 (\cite{Fra00}) and proved by Laver in 1971 (\cite{La71}), the latter is a well quasi-order (briefly: a wqo), i.e.\ there are no infinite descending chains and no infinite antichains; in contrast, we show that \(\cvx\) is not well-founded and has chains and antichains of size continuum (Proposition \ref{prop1}). We prove also other combinatorial properties of $\cvx$, computing in particular its unbounding number \( \mathfrak{b}(\cvxeq) \) and its dominating number \( \mathfrak{d}(\cvxeq) \) (Propositions~\ref{prop_WO} and~\ref{prop:dom_fam}). 

We then explore the problem of classifying \(\LO\) under the equivalence relation induced by \(\cvx\), which we call convex biembeddability and denote by \(\cvxeq\). We obtain the following results (Corollaries~\ref{red_iso_cvxeq} and~\ref{cvxeq_baire_iso}):

\begin{thmintro}
\begin{enumerate-(a)}
\item \( \iso \) is Borel reducible to \( \cvxeq \), in symbols \({\iso} \leq_B  {\cvxeq}\);    
\item\label{thmintro-b} $\cvxeq$ is Baire reducible to $\iso$, in symbols \( \cvxeq \leq_{\text{\scriptsize\textit{Baire}}} \iso \).
\end{enumerate-(a)}
\end{thmintro}

Actually the reduction of part \ref{thmintro-b} of the previous theorem is such that the preimage of any Borel set is a Boolean combination of analytic sets, and hence is also universally measurable. 
Thus, although we are not able to show that \({\cvxeq} \leq_B {\iso}\), the two equivalence relations are similar in some respect, e.g.\ no turbulent equivalence relation Borel reduces to $\cvxeq$, and \( E_1 \not \leq_B {\cvxeq} \) (Corollaries~\ref{turb} and~\ref{E1_cvxeq}).
In particular, \( \cvxeq \) is not complete for analytic equivalence relations and thus \( \cvx \) is not complete for analytic quasi-orders.

In Theorem \ref{thm:reduc_lo_arcs} we establish a connection between linear orders and the theory of proper arcs, showing that:

\begin{thmintro}
${\cvx} \leq_B {\precsim_{\Ar}}$, where
$\precsim_{\Ar}$ is the relation of subarc on the standard
Borel space $\Ar$ of proper arcs.
\end{thmintro}

This allows us to transfer some properties of $\cvx$ and $\cvxeq$ to the
corresponding relations on proper arcs (Corollary~\ref{cor:sumuparcs}).
Moreover, further elaborating on the techniques used to study \( \cvx \), we show that \( \precsim_{\Ar} \) is very complicated from the combinatorial point of view (Theorem~\ref{thm:umboundedanddominatingforarcs}, Corollary~\ref{cor:incomparablearcsandunboundedchains}, Theorems~\ref{thm:basisforarcs} and~\ref{thm:antichainsforarcs}) and it shares most of the features of \( \cvx \) mentioned above, such as the  existence of large antichains, the absence of finite bases, and so on.

If one wants to study knots (instead of proper arcs), it is more natural to consider orders which are ``circular'' (instead of linear). 
The notion of circular order, although not as widespread as that of linear order, is very natural and in fact has been rediscovered several times in
different contexts. 
The oldest mention we found is in \v{C}ech's 1936 monograph (see the English version \cite{Ce69}) and a sample of more recent work is \cite{Meg76, KM05, LM06, BR16, CMR18, PBGGR18, Mat21, GM21, CMMRS23}. 
There is a natural notion of convex subset of a circular order, but the obvious translation of convex embeddability to circular orders fails to be transitive. 
We thus consider its ``transitivization'', called piecewise convex embeddability $\pcvx$, which matches well with the topological notion of piecewise subknot of a knot. 
Then we study the restriction $\pccvx$ of $\pcvx$ to the Polish space \(\CO\) of circular orders with domain \(\N\), together with its
induced equivalence relation $\pccvxeq$. We show that $\pccvxeq$ is
strictly more complicated than $\cvxeq$ in terms of Baire reducibility (Corollary~\ref{cor:piecewiseiscomplex}). Indeed, while $E_1$ is not
Borel reducible to $\cvxeq$, we prove in Theorem \ref{no_action_group}
that

\begin{thmintro}
$E_1 \leq_B {\pccvxeq}$.
\end{thmintro}

We then obtain information about the relation of being a piecewise subknot of a knot, denoted by $\skr$,  
proving in Theorem \ref{thm:lowerforknots} that

\begin{thmintro}
${\pccvx} \leq_B {\skr}$.
\end{thmintro}

An interesting consequence of the above results is that the equivalence
relation associated to \( \skr \) is not induced by a Borel action of a Polish group.
This is in stark contrast with the relation $\equiv_{\Kn}$ of equivalence on knots, which is induced by a Borel action of the Polish group of homeomorphisms of \(S^3\) onto itself (see e.g.\ \cite[Proposition 1.10]{BZ03}).
We also prove a number of combinatorial results concerning \( \skr \) (see Proposition~\ref{prop:chains_skr}, Theorem~\ref{thm:unbounded_and_dominating_for_knots}, Corollary~\ref{cor:incomparable_knots_and_unbounded_chains}, Theorems~\ref{thm:basis_for_knots} and~\ref{thm:antichains_for_knots}),
mimicking those obtained for the quasi-order \( \precsim_{\Ar} \). 
\smallskip

Inspired by piecewise convex embeddability on circular orders, in \cite{IMMRW23} we define the notion of piecewise convex embeddability of linear orders, where the collection of pieces can range in various sets of countable linear orders. 
In the same paper, we also extend the analysis of $\trianglelefteq$ and of piecewise convex embeddability to uncountable linear orders.

The paper is organized as follows. In Section 2 we recall some
background about descriptive set theory and Borel
reducibility. Moreover, we introduce and prove basic properties of
(countable) linear and circular orders.  In Section 3 we define and
analyse convex embeddability on countable linear orders and countable
circular orders, proving some combinatorial properties of these
quasi-orders and several results about the complexity of the
corresponding equivalence relations with respect to\ Borel
reducibility. 
In the last section we look at the connections between
convex embeddability on \(\LO\) and piecewise convex embeddability on
\(\CO\) on one side, and the relations of subarc on proper arcs and
piecewise subknot on knots on the other one. We also establish some combinatorial properties of the subarc and piecewise subknot relations. 

\section{Preliminaries}

\subsection{Borel reducibility}\label{subsection:preliminaries-BR}

In this section we introduce some basic definitions and results from descriptive set theory that will be used in the sequel; the standard references are \cite{Kec95,Gao09}.

A \textbf{Polish space} is a separable and completely metrizable
topological space. 
A subset of a Polish space $X$ is \textbf{Borel} if it belongs to the smallest $\sigma$-algebra on $X$ containing all open subsets of $X$. Recall that \(A \subseteq X\) has the \textbf{Baire property} if there exists some open set \(U\) such that \(A \mathrel{\triangle} U\) is meager, i.e.\ a countable union of sets whose closure has empty interior. All Borel sets have the Baire property.

A \textbf{standard Borel space} is a pair $(X, \mathcal{B})$ where
$X$ is a set, $\mathcal{B}$ is a $\sigma$-algebra on $X$, and there
is a Polish topology on $X$ for which $\mathcal{B}$ is precisely the
collection of Borel sets. The elements of \( \mathcal{B} \) are
called Borel sets of~\(X\).
In particular, every Polish space is standard Borel when equipped with its \(\sigma\)-algebra of Borel sets. 

Let $X$ and $Y$ be Polish or standard Borel spaces. A function $\varphi \colon X \to Y$ is \textbf{Borel} if the preimage of any Borel subset of $Y$ is Borel in $X$.
If $X$ is Polish we say that $\varphi$ is \textbf{Baire measurable} if the preimage of any Borel subset of $Y$ has the Baire property.

Let \( X \) be a standard Borel space. A subset   $A\subseteq X$ is \textbf{analytic} (or $\mathbf{\Sigma}_1^1$) if it is the Borel image of a standard Borel space, and it is \textbf{coanalytic} (or $\mathbf{\Pi}_1^1$) if   $X\backslash A$ is analytic. 
By \(D_2(\mathbf{\Pi}_1^1)\) we denote the class of sets which are the intersection of an analytic set and a coanalytic set.

Let \(X\) and \(Y\) be topological spaces and \(A \subseteq X\), \(B \subseteq Y\). We say that \(A\) is \textbf{Wadge reducible} to \(B\), in symbols \(A \leq_W B\), if there is a continuous map \(\varphi\colon X \to Y\) such that $x \in A \iff \varphi(x) \in B$, for all $x \in X$.

Let \(\Gamma\) be a class of sets in Polish spaces. If \(Y\) is a
Polish space, we say that the subset \(A\) of \(Y\) is
\textbf{\(\Gamma\)-hard} if \(B \leq_W A\) for any
\(B \in \Gamma(X)\) with \(X\) a zero-dimensional Polish
space. If moreover \(A \in \Gamma(Y)\), we say that \(A\) is
\textbf{\(\Gamma\)-complete}.\medskip

An important line of research within descriptive set theory is the
study of definable equivalence relations, which are typically compared
using the next definition.

Let $X$ and $Y$ be sets and consider $E$ and $F$ equivalence
relations on $X$ and $Y$, respectively. A function
$\varphi \colon X \to Y$ is called a \textbf{reduction} from
$E$ to $F$ if
${x_1}\mathrel{E}{x_2}\iff {\varphi(x_1)}\mathrel{F}{\varphi(x_2)},$
for all $x_1, x_2 \in X$.
We say that $E$ is \textbf{Borel reducible} to $F$, and write
$E \leq_B F$, if \(X\) and \(Y\) are standard Borel spaces and there
exists a Borel map \(\varphi\) reducing \(E\) to \(F\). The
equivalence relations $E$ and $F$ are \textbf{Borel bireducible},
${E} \sim_B {F}$ in symbols, if both $E \leq_B F$ and $F \leq_B E$.
Finally, we say that \(E\) is \textbf{Baire reducible} to \(F\), and we write
\({E} \leq_{\text{\scriptsize\textit{Baire}}} {F}\), if
\(X\) and \(Y\) are topological spaces and there exists a Baire
measurable map \(\varphi\colon  X \to Y\) reducing \(E\) to \(F\).

Let \( \Gamma \) be a collection of equivalence relations on
standard Borel spaces. We say that an equivalence relation \(E\) is
\textbf{complete for \(\Gamma\)} (or \textbf{\(\Gamma\)-complete})
if it belongs to \(\Gamma\) and any other equivalence relation in
\(\Gamma\) Borel reduces to $E$. 

A topological group is \textbf{Polish} if its underlying topology is
Polish. Examples of Polish groups include the group of permutations of
natural numbers $S_\infty$ with the topology inherited as a subspace
of the Baire space \( \N^\N \), and the group of homeomorphisms of
\(S^3\) into itself with the topology induced by the uniform metric.

If the Polish group $G$ acts on the standard Borel space $X$ we denote by $E_G^X$ the orbit equivalence relation induced by the action. 
%
%
%
An important class of analytic equivalence relations are those induced by a Borel action of \(S^\infty\).
Among these are all the isomorphism relations on the countable models of a first-order theory.
An analytic equivalence relation is \textbf{\( S_\infty \)-complete} if it is complete for the class of equivalence relations arising from a Borel action of $S_\infty$ on a standard Borel space. 
We say that an equivalence relation \(E\) is \textbf{classifiable by countable structures} if \(E\) is Borel reducible to some \(E_{S_\infty}^Y\).

\begin{thm}[{H. Friedman-Stanley, see~\cite{FS89,Gao09}}]
\label{thm_0.9}
Let $\cong_{\mbox{\tiny \emph{C-GRAPH}}}$ and
$\cong_{\mbox{{\tiny \emph{GRAPH}}}}$ be the isomorphism
relations on, respectively, the Polish space \emph{C-GRAPH} of countable connected graphs and the Polish space \emph{GRAPH} of  countable
graphs. Then
${\cong_{\mbox{\tiny \emph{C-GRAPH}}}} \sim_B {\cong_{\mbox{\tiny \emph{GRAPH}}}}$,
and both equivalence relations are \( S_\infty \)-complete.
\end{thm}


Let $E$ and $F$ be equivalence relations on standard Borel spaces $X$ and $Y$ respectively. Let $A\subseteq X$. We say that   \textbf{$E \rest A$ is Borel reducible to $F$} if there is a Borel map $\varphi \colon X\to Y$, still called a Borel reduction of \( E \rest A \) to \( F \), such that for every \(x,y \in A\), \(x\mathrel{E} y \iff \varphi(x)\mathrel{F}\varphi(y)\). This definition is equivalent to the one given in \cite{CMMR18,CMMR20} (where $\varphi$ is required to be defined only on $A$) by a theorem of Kuratowski (see \cite[Theorem 12.2]{Kec95}).


\begin{defn}\label{def:QuasiClass}
We say that an equivalence relation $E$ on a Polish space $X$ is
\(\mathbb{\sigma}\)\textbf{-classifiable by countable structures} if there exists a countable partition 
$(X_i)_{i \in I}$ of $X$ such that for
all $i \in I$:
\begin{enumerate-(i)}
\item $X_i$ is closed under $E$ (i.e.\ if $x \in X_i$ and $y \mathrel{E} x$
then $y \in X_i$);
\item $X_i$ has the Baire property;
\item $E \rest X_i$ is Borel reducible to \({\cong_{\mbox{{\tiny \emph{GRAPH}}}}}\).
\end{enumerate-(i)}
\end{defn}

Clearly, if an equivalence relation \(E\) is classifiable by countable
structures then it is \(\sigma\)-classifiable by countable structures.

\begin{prop}\label{Baire_red_for_qclas_count_str}
Let \(E\) be an equivalence relation defined on a Polish space
\(X\). If \(E\) is \(\sigma\)-classifiable by countable structures,
then
\({E} \leq_{\text{\scriptsize \textit{Baire}}}
{\cong_{\mbox{\tiny \emph{GRAPH}}}}\).
\end{prop}

\begin{proof}
Assume that \(E\) is \(\sigma\)-classifiable by countable structures and fix sets \(X_i\) witnessing this. Then by Theorem \ref{thm_0.9} for each \(i \in I\) there exists a Borel reduction \(\varphi_i\) from \(E\rest X_i\) to ${\cong_{\mbox{\tiny C-GRAPH}}}$, so that $\varphi_i(x)$ is an infinite connected graph for every $x \in X_i$ (in particular, it is not isomorphic to the graph consisting of a single isolated vertex). Let \( \tilde{\varphi}_i \colon X \to \text{GRAPH} \) be defined by
\[
\tilde{\varphi}_i(x) = \varphi_i(x) \sqcup A_i,
\]
where \( A_i \) is the graph consisting of \( i \)-many isolated vertices. It is easy to check that \( \tilde{\varphi}_i \) is still a Borel function and it reduces \( E \rest X_i \) to \( \cong_{\mbox{\tiny GRAPH}} \). Finally, define \(\varphi \colon X \to \mbox{GRAPH}\) by setting \( \varphi(x)=\tilde{\varphi}_i(x) \), where \(i\) is the unique index of the subset \(X_i\) of \(X\) to which \(x\) belongs.

We first show that \(\varphi\) is a reduction. Let \(x, y\) be two elements of \(X\) such that \(x \mathrel{E} y\). Since \(X_i\) is closed under \(E\) for every \(i \in I\), there exists \(i_0 \in I\) such that \(x,y \in X_{i_0}\). Then \(\varphi_{i_0}(x) \cong_{\mbox{\tiny C-GRAPH}} \varphi_{i_0}(y)\), and so \(\varphi(x) \cong_{\mbox{\tiny GRAPH}} \varphi(y)\). Conversely, suppose that \(\varphi(x)= \varphi_i(x) \sqcup A_i \cong_{\mbox{\tiny GRAPH}} \varphi_j(y) \sqcup A_j =\varphi(y)\), for some \(i,j \in I\), \( x \in X_i \), and \( y \in X_j \). Since isomorphism between graphs preserves connected components, we must have \(i=j\) because \( \varphi(x) \) contains \( i \)-many isolated vertices and \( \varphi(y) \) contains \( j \)-many isolated vertices, and moreover \(\varphi_i(x) \cong_{\mbox{{\tiny C-GRAPH}}} \varphi_i(y)\) because those are the only infinite connected components in \( \varphi(x) \) and \( \varphi(y) \), respectively. Since \( \varphi_i \) was a reduction we get \(x \mathrel{E} y\), as desired.

Now take a Borel subset \(A\) of \(\mbox{GRAPH}\). Then
\[
\varphi^{-1}(A)= \bigcup_{i \in I} \big(\varphi^{-1}(A) \cap X_i\big) = \bigcup_{i \in I} \big(\tilde{\varphi}_i^{-1}(A) \cap X_i\big).
\] 
Since \(X_i\) has the BP and \(\tilde{\varphi}_i\) is Borel for every \(i \in I\), we have that \(\tilde{\varphi}_i^{-1}(A) \cap X_i\) has the BP for each \(i\). Hence also \(\varphi^{-1}(A)\) has the BP and \(\varphi\) is a Baire measurable reduction.
\end{proof}

Not all orbit equivalence relations are Borel reducible, or even Baire 
reducible, to an $S_\infty$-complete equivalence relation: Hjorth
isolated a sufficient condition for this failure, called \textbf{turbulence}.


\begin{thm}[\cite{Hjo00}, Corollary 3.19]\label{turb_act}
There is no Baire measurable reduction of a turbulent orbit
equivalence relation to any \(E_{S_\infty}^Y\).
\end{thm}

Let \(E_1\) be the equivalence relation defined on \(\R^\N\) by \((x_n)_{n \in \N} \mathrel{E_1} (y_n)_{n \in \N}\) if and only if there exists $m$ such that $x_n = y_n$ for all $n \geq m$.
We also use the tail version \(E_1^t\), defined by setting \((x_n)_{n \in \N} \mathrel{E_1^t} (y_n)_{n \in \N}\) if and only if there exist $n,m$ such that $x_{n+k} = y_{m+k}$ for all $k$. 
Notice that $E_1$ and $E_1^t$ are Borel bireducible with the analogous relations defined on $(2^\N)^\N$, called $E_0(2^\N)$ and $E_t(2^\N)$ in \cite{DJK94}. In the proof of \cite[Theorem 8.1]{DJK94} it is shown that \(E_t(2^\N) \leq_B E_0(2^\N)\), while the opposite reduction is mentioned in the observation immediately following that proof. This yields:

\begin{prop}\label{E1}
\({E_1} \sim_B {E_1^t}\).
\end{prop}

The following result of Shani about $E_1$ generalizes a classical
theorem by Kechris and Louveau \cite{KL97}. (The additional part follows from the fact that by~\cite[Theorem 8.38]{Kec95} every Baire measurable map between Polish spaces is actually continuous on a comeager \( G_\delta \) set.)

\begin{thm}[{\cite[Theorem 4.8]{Sha21}}]\label{E1_orbit}
The restriction of \(E_1\) to any comeager subset of \(\R^\N\) is not
Borel reducible to an orbit equivalence relation. Thus in particular $E_1 \not\leq_{\text{\scriptsize \textit{Baire}}} {\iso}$.
\end{thm}


Let $E$ be an equivalence relation on a standard Borel space $X$. The \textbf{Friedman-Stanley jump} of $E$, introduced by Friedman and Stanley in \cite{FS89} and denoted by $E^{+}$, is the equivalence relation on the space $X^\N=\{(x_n)_{n \in \N} \mid x_n \in X \}$ defined by
$$(x_n)_{n \in \N} \mathrel{E^+} (y_n)_{n \in \N} \iff \{[x_n]_E \mid n \in \N\}=\{[y_n]_E \mid n \in \N\}.$$  


\begin{prop}[see \cite{Gao09}]\label{prop 0.10}
Let \( E \) and \( F \) be equivalence relations on standard Borel spaces. Then $E\leq_B E^+$, and if $E\leq_B F$ then $E^+\leq_B F^+$.
\end{prop}

One can transfer many of the above definitions concerning equivalence relations to the wider context of binary relations and, in particular, analytic quasi-orders. We just recall a few results in this direction.

\begin{thm}[\cite{LR05}]\label{thm:compl_graph}
Every analytic quasi-order Borel reduces to the embeddability
relation between countable (connected) graphs, i.e.\ the latter
relation is complete for analytic quasi-orders.
\end{thm}

Every analytic quasi-order \( R \) on a standard Borel space \( X \)
canonically induces the analytic equivalence relation \( E_R \) on the
same space defined by \( x \mathrel{E_R} y \iff {x \mathrel{R} y} \wedge {y \mathrel{R} x}\). The complexities of \( R \) and \( E_R \) are linked by the following
result.

\begin{prop}[\cite{LR05}]\label{compl_qo}
If a quasi-order $R$ on a standard Borel space $X$ is complete for
analytic quasi-orders, then \( E_R \) is complete for analytic
equivalence relations.

\end{prop}

\subsection{Countable linear orders} \label{subsec:ctbllinord}

Any $L \in 2^{\N\times\N}$ can be seen as a code of a binary relation
on \(\N\), namely, the one relating \( n \) and \( m \) if and only if
$L(n,m) = 1$. Denote by $\LO$ the set of codes for linear orders on
$\N$, i.e.
\[
\LO = \{L \in 2^{\N\times\N} \mid L \mbox{ codes a reflexive linear order on } \N\}.
\]
When $L \in \LO$ we denote by $ \leq_L $ the order on \(\N\)
coded by \( L \), and by \( <_L \) its strict part.

It is easy to see that $\LO$ is a closed subset of the Polish space
$2^{\N\times\N}$, thus it is a Polish space as well. Given
\( L \in \LO \), a neighbourhood base of \( L \) in \( \LO \) is
determined by the sets
\[
\{ L' \in \LO \mid L'\rest n = L\rest n \}
\]
where \( n \) varies over \( \N \) and \( L\rest  n= L'\rest n \) means that \( m \leq_L m' \iff m \leq_{L'} m' \) for every  \( m, m' < n \). 
We also
denote by \(\WO\) the set of all well-orders on \(\N\), and recall that it is
a proper coanalytic subset of \( \LO \).

We denote by \( \preccurlyeq \) the quasi-order of embeddability
on linear orders, that is: $L \preccurlyeq L'$ if there exists an injection $f$ from \( L \) to \( L' \), called embedding, such that $n \leq_L m \Rightarrow f(n) \leq_{L'} f(m)$ (equivalently $n \leq_L m \iff f(n) \leq_{L'} f(m)$) for every $n,m \in L$.
The restriction \( \preccurlyeq_{\LO} \) of \( \preccurlyeq \) to \( \LO \) is clearly an analytic quasi-order. In contrast with 
Theorem~\ref{thm:compl_graph}, the relation \(\preccurlyeq_{\LO}\) is far from being complete because it is combinatorially simple and it is indeed a wqo. 
Moreover $\LO$ has a maximal element under \(\preccurlyeq_{\LO}\), the equivalence class of non-scattered linear orders (recall that a linear order is scattered if the rationals do not embed into it).

The isomorphism relation on $\LO$ is denoted by \( \iso \), and
it is an analytic equivalence relation. 

\begin{thm}[\cite{FS89}]\label{thm 2}
$\iso$ is $S_{\infty}$-complete.
\end{thm}

Recall that ${E} \leq_B {E^+}$ for any analytic equivalence relation
$E$ (Proposition~\ref{prop 0.10}). In the case of $\iso$, we
also have the converse.

\begin{prop}[Folklore]{\label{cor 2}}
$(\iso)^+ \sim_B {\iso}$.
\end{prop}

\begin{proof}
By Theorem \ref{thm_0.9} 
and Theorem \ref{thm 2},
we have that
${\cong_{\mbox{\tiny C-GRAPH}}} \sim_B {\cong_{\mbox{\tiny GRAPH}}}
\sim_B {\iso}$, so it is enough to prove that
$(\cong_{\mbox{\tiny C-GRAPH}})^+ \leq_B {\cong_{\mbox{\tiny
GRAPH}}}$ because \( (\cong_{\mbox{\tiny C-GRAPH}})^+ \sim_B (\iso)^+ \) 
by Proposition~\ref{prop 0.10}. Given a sequence of countable connected graphs
$(A_n)_{n \in \N}$, let
$G_A = \bigsqcup_{n, i \in \N} A_{n,i}$ be the disjoint union of
the graphs $A_{n,i}$, where $A_{n,i}\cong A_n$ for every $n,i \in
\N$. Then the Borel map from  \( ({\mbox{C-GRAPH}})^\N\) to \({{\mbox{GRAPH}}}\) which sends  \((A_n)_{n \in \N}\) to \(G_A\) is a reduction of $(\cong_{\mbox{\tiny C-GRAPH}})^+$ to~${\cong_{\mbox{\tiny GRAPH}}}$.
\end{proof}

We need to deal also with finite linear orders, which are missing in
\(\LO\). For this reason, we let \(\Lin\) be the subset of
$2^{\N\times\N}$ consisting of all (codes for) linear orders defined either on a finite subset
of \(\N\) or on the whole $\N$. Thus $\Lin$ is the union of $\LO$ and $\Fin$, where \(\Fin \subset \Lin\) is the set of (codes for)
finite linear orders. It is easy to see that $\Lin$ is a $F_\sigma$ subset of $2^{\N\times\N}$, and hence it is a standard Borel space, and that isomorphism on $\Lin$ is induced by a Borel action of $S_{\infty}$. 

If \(L \in \Lin\) we denote by \(L\) also its
domain. For convenience, sometimes we use the notation
\(n_L\) to emphasize that \(n\) is an element of the domain of \(L\).

We recall some isomorphism invariant operations on the class of linear orders that are
useful to build Borel reductions. They can all be construed as Borel maps from \( \Lin \), \( (\Lin)^n\), or \( \Lin^\N \)
to \( \Lin \), and their restriction to \( \LO \) has range contained in \( \LO \).
\begin{itemizenew}
\item The \textbf{reverse} \( L^* \) of a linear order \( L \) is the linear order on
the domain of \( L \) defined by setting
$x \leq_{L^*} y \iff y \leq_L x$.
\item If \( L \) and \( K \) are linear orders, their \textbf{sum} \( L+K \) is the linear order defined on
the disjoint union of \( L \) and \( K \) by setting
${x} \leq_{L+K} {y}$ if and only if either $x \in L$ and
$y\in K$, or $x, y \in L$ and $x \leq_{L} y$, or \( x,y \in K \) and \( x \leq_K y \). 
\item
In a similar way, given a linear order \(K\) and a sequence of linear orders
\((L_k)_{k \in K}\) we can define the \textbf{\(K\)-sum} $\sum_{k \in K} L_k$ on the disjoint union of the \( L_k \)'s by setting \( x \leq_{\sum_{k \in K} L_k} y \) if and only if there are \( k <_K k' \) such that \( x \in L_k \) and \( y \in L_{k'} \), or \( x,y \in L_k \) for the same \( k \in K \) and \( x \leq_{L_k} y \). Formally, \( \sum_{k \in K} L_k \) is thus defined on the set \( \{ (x,k) \mid k \in K , x \in L_k \} \) 	by stipulating that \( (x,k) \leq_{\sum_{k \in K} L_k} (x',k') \) if and only if \( k <_K k' \) or else \( k = k' \) and \( x \leq_{L_k} x' \).
\item The \textbf{product} \(LK\) of two linear orders \( L \) and \( K \) is the cartesian product \( L \times K\)
ordered antilexicographically. Equivalently, \( LK  = \sum_{k \in K} L \). 

\end{itemizenew}
%

For every \(n \in \N\), we denote by \(\boldsymbol{n}\) the element of
\(\Fin\) with domain \(\{0,...,n-1\}\) ordered as usual. Similarly, for every infinite  ordinal \( \alpha < \omega_1 \) we fix a well-order \( \boldsymbol{\alpha} \in \LO \) with order type \( \alpha \). 
We also fix computable copies of
\( (\N,{\leq}) \), \( (\Z,{\leq}) \) and \( (\Q,{\leq}) \) in \( \LO \), and denote them by \(\o\), \(\z\)
and \(\eta\), respectively.
We denote by \( \min L \) and \( \max L \) the minimum and maximum of \( L \), if they exist.
Finally, we denote by \(\Scat \subseteq \Lin\) the set of scattered linear orders.

\begin{defn}
A subset \( I \) of the domain of a linear order \( L \) is
\textbf{(\( L \)-)convex} if \( x \leq_L y \leq_L z \) with
\( x,z \in I \) implies \( y \in I \). An \( L \)-convex set is
\textbf{proper} if it is neither empty nor the entire \( L \).
\end{defn}

An \textbf{initial segment} of a linear order \( L \) is a subset \( I \) of its domain which is \( \leq_L \)-downward closed, i.e.\ \( x \in I \) whenever \( x \leq_L y \) for some \( y \in I \). Dually, \( I \subseteq L \) is a \textbf{final segment} of \( L \) if it is \( \leq_L \)-upward closed, i.e.\ if \( y \in I \) and \( y \leq_L x \) imply \( x \in I \). Clearly, initial and final segments are convex of \( L \).

If \(m,n \in L\), we adopt the notations $[m,n]_L$, $(m,n)_L$,
$(-\infty,n]_L$, $(-\infty,n)_L$, $[n,+\infty)_L$, and $(n,+\infty)_L$
to indicate the obvious \( L \)-convex sets. Notice however that not all $L$-convex sets are of one of these forms.

Given \(L \in \LO\), we write \(L_0 \subseteq L\) (resp.
\(L_0 \subset L\)) if \(L_0\) is a (resp.\ proper) sub-order of \(L\),
and \(L_0 \csube L\) (resp. \(L_0 \csub L\)) if \(L_0\) is a
(resp.\ proper) convex subset of \(L\). If
\(L_0,L_1 \subseteq L\), we write \(L_0 \leq_L L_1\) (resp. \(L_0 <_L L_1\)) iff
\(n \leq_L m\) (resp.\ \(n <_L m\)) for every \(n \in L_0\) and \(m \in L_1\). Notice that
if \(L_0 \leq_L L_1\) then either \(L_0\) and \(L_1\) are disjoint, in which case \(L_0 <_L L_1\), or
the only element in their intersection is $\max L_0 = \min L_1$.

We need to recall some other basic notions about linear orders
(see~\cite{Ros82}).
Let \(L\) be a linear order. The \textbf{(finite) condensation} of \(L\)
is determined by the map \(c_F^L \colon L \to \mathscr{P}(L)\)
defined by \(c_F^L(n)=\{m \mid [n,m]_L \cup [m,n]_L \text{ is finite}\} \) for
every \(n \in L\). It is immediate that if \(m \in c_F^L(n)\) then
\(c_F^L(m)=c_F^L(n)\), while if \( m \notin c_F^L(n) \) then \( c_F^L(n) \cap c_F^L(m) = \emptyset \). 
We call a set \(c_F^L(n)\) a
\textbf{condensation class}. 
A condensation class may be finite or infinite, and in the latter case its order type is one of \(\o\), \(\o^*\) and \( \z \). 
We denote by \(L_F\) the set of condensation
classes of~\(L\).
In the sequel we use the basic properties of condensation classes which are collected in the following proposition.

\begin{prop}\label{cond_classes_isom}
Let \( L \) be any linear order.
\begin{enumerate-(a)}
\item For every \(\ell \in L\), \(c_F^L(\ell)\) is convex.   
\item \(\bigcup_{\ell \in L} c_F^L(\ell) = L\), and 
\( c^L_F(\ell) \cap c^L_F(\ell') = \emptyset \) if \( c^L_F(\ell) \neq c^L_F(\ell') \); hence \( L_F \) is a partition of \(L \).
\item If \(c_F^L(\ell)\) and \(c_F^L(\ell')\) are two different
condensation classes, then \(\ell <_L \ell'\) if and only if
\(c_F^L(\ell) <_L c_F^L(\ell')\); hence \(L_F\) is linearly ordered.
\item Let \(L, L'\) be linear orders. If \(f\) is an isomorphism
from \(L\) to \(L'\) then the restriction of \(f\) to each
\(c_F^L(\ell)\) is an isomorphism between \(c_F^L(\ell)\) and
\(c_F^{L'}(f(\ell))\) and hence
\(|c_F^L(\ell)|=|c_F^{L'}(f(\ell))|\). Moreover, \(L_F \cong L'_F\) via the well-defined map \( c^L_F(\ell) \mapsto c^{L'}_F(f(\ell)) \).
\end{enumerate-(a)}
\end{prop}	

This condensation is useful to prove results as the next one.

\begin{lem}\label{lem:isom_zetaL}
Given two linear orders \(L,L'\), \(\z L \cong \z L'\) if and only if \(L \cong L'\).
\end{lem}
\begin{proof}
For the nontrivial direction, notice that \( c^{\z L}_F(i,n) = \z \times \{ n \} \), and similarly for the condensation classes of \( \z L' \).
It follows that \( (\z L)_F \cong L \) and \( (\z L')_F \cong L' \). By Proposition~\ref{cond_classes_isom}, if \( \z L \cong \z L' \) 
then \( (\z L)_F \cong (\z L')_F \), hence \( L \cong L' \).
\end{proof}

We conclude this section recalling the definition of the powers of \(\Z\) and some of their properties. When $\a$ is an ordinal we can define $\Z^\a$ in two equivalent ways: by induction on $\a$ (\cite[Definition 5.34]{Ros82}) and by explicitly defining a linear order on a set (\cite[Definition 5.35]{Ros82}); the latter can actually be used to define \(\Z^L\) for any linear order \(L\).

\begin{defn} \label{def:Z^alpha}
\begin{enumerate-(1)}
\item \(\Z^0 = \boldsymbol{1}\),
\item \(\Z^{\a+1}=(\Z^\a \o)^* + \Z^\a + \Z^\a \o\),
\item
\(\Z^\alpha = \big(\sum_{\beta < \alpha} \Z^\beta \o\big)^* +
\boldsymbol{1} + \sum_{\beta < \alpha} \Z^\beta \o\) if $\a$ is limit.
\end{enumerate-(1)}
\end{defn}

\begin{defn}\label{pow_z_def_theor}
Let \(L\) be a linear order. For any map \(f \colon L \to \Z\), we
define the support of \(f\) as the set
\(\Supp(f)=\{n \in L \mid f(n) \neq 0\}\) . The \(L\)-power of
\(\Z\), denoted by \(\Z^L\), is the linear order on
\(\{f\colon L \to \Z \mid \Supp(f) \text{ is finite}\}\) defined by
the following: if \(f,g\colon L \to \Z\) are maps with finite
support let \(f \leq_{\Z^L}g\) if and only if \(f=g\) or
\(f(n_0) <_{\Z} g(n_0)\) where
\(n_0 = \max\{n \in \Supp(f) \cup \Supp(g) \mid f(n) \neq g(n)\}\).
\end{defn}

Sometimes we need the following properties (see~\cite[Section~3.2]{CCM19}).

\begin{prop}\label{power_Z}
For all ordinals \(\beta<\alpha\), we have
\[
\Z^\alpha \cong \bigg(\sum_{\beta \leq \gamma < \alpha} \Z^\gamma
\o\bigg)^* + \bigg(\sum_{\beta \leq \gamma < \alpha} \Z^\gamma
\o\bigg).
\]
\end{prop}

\begin{prop}\label{power_Z^L}
For any linear orders \(L\) and \(L'\) we have
\begin{enumerate-(a)}
\item \label{power_Z^L-a}
\((\Z^L)^*=\Z^L\),
\item \(\Z^{L+L'} \cong \Z^L \Z^{L'}\),
\item if \(L\) is countable and not a well-order then there is a
countable ordinal \(\a\) such that \(\Z^L\cong \Z^\a \eta\).
\end{enumerate-(a)}
\end{prop}

\subsection{Circular orders}\label{sec:circular orders}

We now describe the basic notation and notions regarding circular orders. The prototype of a circular order is the unit circle $S^1$ traversed counterclockwise, which we denote by $C_{S^1}$. 


\begin{defn}\label{def1}(\cite[Definition 2.1]{KM05})
A ternary relation $C \subset X^3$ on a set $X$ is said to be a
\textbf{circular order} if the following conditions are satisfied  for every $x,y,z,w \in X$:
\begin{enumerate-(i)}
\item Cyclicity: $(x,y,z) \in C \Rightarrow (y,z,x) \in C$;
\item Antisymmetry and reflexivity: $(x,y,z) \in C \wedge (y,x,z) \in C \iff x=y \lor y=z \lor z=x$;
\item\label{def:co_trans} Transitivity: $(x,y,z) \in C \Rightarrow \forall t ((x,y,t) \in C \vee (t,y,z) \in C)$;
\item Totality: $(x,y,z) \in C \lor (y,x,z) \in C$.
\end{enumerate-(i)}
\end{defn}
Notice that, assuming the other conditions, \ref{def:co_trans} is equivalent to asserting that $(x,y,z) \in C$ and $(x,z,w) \in C$ imply $(x,y,w) \in C$ whenever $x \neq z$. In the sequel we often make use of this reformulation. Definition \ref{def1} is different from \cite[5.1]{Ce69}: indeed, the latter characterizes the strict relation associated to $C$, i.e.\ the set of all triples \( (x,y,z)\) such that $(x,y,z) \in C$ and \( x,y,z \) are all distinct.

By abuse of notation, when \(C\) is a circular order on \( X \) we write
$C(x,y,z)$ instead of 
\( (x,y,z) \in C \),
for  \(x, y, z \in X\).
The \textbf{reverse} \( C^* \) of a circular order \( C \) on \( X \) is the circular order on \( X \) defined by \( C^*(x,y,z) \iff C(z,y,x) \) for all \( x,y,z \in X \).

Let $C$ and $C'$ be circular orders on sets \(X\) and \(X'\),
respectively. We say that $C$ is \textbf{embeddable} into $C'$, and
write $C \preccurlyeq_{c} C'$, if there exists an injective function
$f\colon X \to X'$, called embedding, such that for every $x,y,z \in X$,
\( C(x,y,z) \Rightarrow C'(f(x),f(y),f(z)) \) (notice that by totality and antisymmetry, one also has \(C'(f(x),f(y),f(z)) \Rightarrow C(x,y,z)\)).
We say that $C$ and $C'$ are \textbf{isomorphic}, and write $C \cong_{c} C'$, if there exists $f$ as above which is a bijection (in which case \( f \) is called isomorphism).

For a circular order, the notions of successor and predecessor of an element are meaningless. However, we can still define a notion of immediate successor or immediate predecessor.

\begin{defn}
Given a circular order \(C\) on the set \(X\) and \(x, y \in X\), we say that \(x\) is the \textbf{immediate predecessor} (resp.\ \textbf{immediate successor}) of \(y\) in \(C\) if $x \neq y$ and \(C(x,y,z)\) (resp.\ \(C(y,x,z)\)) for every \(z \in X\).
\end{defn}

\begin{defn}\label{rel_LO-CO}
Given a linear order $L$, we define a circular order $C[L]$ by setting   $C[L](x,y,z)$ if and only if one of the following conditions is satisfied:
\[
x \leq_L y \leq_L z,\quad y \leq_L z \leq_L x,\quad z \leq_L x \leq_L y. 
\]
\end{defn}

Notice that every circular order \(C\) is of the form \(C[L]\) for some (in general non unique) linear order \(L\).
Clearly, for two linear orders \(L\) and \(L'\) such that \(L \preccurlyeq L'\) we have \(C[L] \preccurlyeq_{c} C[L']\). 

Denote by $\CO$ the set of codes for circular orders on $\N$, i.e. 
\[\CO = \{C \in 2^{\N \times \N \times \N} \mid C \text{ codes a circular order on }\N\}.\]
Since $\CO$ is a closed subset of the Polish space
$2^{\N \times \N \times \N}$, we have that it is a Polish space as
well. Denote by \(\preccurlyeq_{\CO}\) and \(\cong_{\CO}\) the
restriction of the relations of embeddability \( \preccurlyeq_c \) and isomorphism \( \cong_c \) to \(\CO\), respectively. It is immediate that both \( \preccurlyeq_{\CO} \) and \( \cong_{\CO} \) are analytic.

\begin{prop}
\(\preccurlyeq_{\CO}\) is a wqo.
\end{prop}

\begin{proof}
Recall that a quasi-order $(X, {\leq_X})$ is a wqo if for every sequence \((x_n)_{n \in \N}\) of elements of \(X\), there exist \(n<m\) such that \(x_n \leq_X x_m\).
Suppose that \((C_n)_{n \in \N}\) is a sequence of elements of
\(\CO\). For every \(n \in \N\) consider the linear order \(L_n\) defined by 
\[
x\leq_{L_n} y \iff C_n(0,x,y) \land (y=0 \Rightarrow x=0).
\] 
Notice that $C[L_n] = C_n$.

Since the embeddability relation
\(\preccurlyeq_{\LO}\) on \(\LO\) is wqo, there
are \(n<m\) such that
\(L_n\preccurlyeq_{\LO} L_m\) and hence 
\(C_n \preccurlyeq_{\CO} C_m\).
\end{proof}

The isomorphism $\cong_{\CO}$ is an equivalence relation on $\CO$. Clearly, for $L, L' \in \LO$, we have
that \({{L}\iso {L'}}\) implies
\({C[L]} \cong_{\CO} {C[L']}\).  The converse implication
is not true, as showed by \(C[\o +\boldsymbol{1}]\) and
\(C[\o]\), for which we have
$C[\o +\boldsymbol{1}] \cong_{\CO} C[\o]$, but
$\o + \boldsymbol{1} \ncong_{\LO} \o$.

\begin{thm}\label{thm:iso_co_bired_iso_lo}
${{\cong}_{\CO}} \sim_B {\iso}.$
\end{thm} 

\begin{proof}
For the Borel reduction from \(\cong_{\CO}\) to \(\iso\), it
is enough to note that \(\cong_{\CO}\) is an equivalence relation
arising from a Borel action of the group $S_\infty$. Then
\({{\cong}_{\CO}} \leq_B {\iso}\) by Theorem \ref{thm 2}.

For the converse, consider the Borel map
$\varphi\colon\LO \to \CO$ defined by
$$\varphi(L) = C[\boldsymbol{1}+\z L].$$
If $L\iso L'$ we have immediately that
\(\varphi(L) \cong_{\CO} \varphi(L')\). Suppose now that
\(\varphi(L) \cong_{\CO} \varphi(L')\) via the map $f$. Since
$\boldsymbol{1}$ is the only element which has no immediate
successor in both $\varphi(L)$ and $\varphi(L')$, we have that
\(f(\boldsymbol{1}) = \boldsymbol{1} \). Thus \(\z L \iso \z L'\)
and by Lemma \ref{lem:isom_zetaL} we obtain \(L \iso L'\).
\end{proof}

\section{Convex embeddability}\label{sec:cvx emb}

This is the main definition of the paper.

\begin{defn}[\cite{BCP73}]
Let $L$ and $L'$ be linear orders. We say that an embedding \(f\)
from $L$ to $L'$ is a \textbf{convex embedding} if $f(L)$ is an
\( L' \)-convex set. We write $L\trianglelefteq L'$ when such \(f\)
exists, and call \textbf{convex embeddability} the resulting binary relation.
\end{defn}

\begin{remark}\label{def cvx}
Notice that	\(L \trianglelefteq L'\) if and only if 
\[
L' \cong L_{l} + L + L_{r},
\]
for some (possibly empty) $L_l$ and $L_r$, if and only if \(L\) is
isomorphic to an \(L'\)-convex set.
\end{remark}

While \(L \preccurlyeq \eta\), for every countable linear order \(L\),
we have \(L \trianglelefteq \eta\) if and only if \(L\) has order type
\(\boldsymbol{1} \), \( \eta \), \( \boldsymbol{1}+\eta \), \( \eta +\boldsymbol{1}\) or
\(\boldsymbol{1}+\eta + \boldsymbol{1} \).

One easily sees that the restriction of convex embeddability to the
Polish space \( \LO \) is an analytic quasi-order, which we denote by
\(\cvx\).  
The strict part of \( \trianglelefteq_{\LO} \) is denoted by \( \triangleleft_{\LO} \), that is,
\( L \triangleleft_{\LO} L' \) if and only if \( L \trianglelefteq_{\LO} L' \) but \( L' \not{\trianglelefteq}_{\LO} L \).
We call \textbf{convex biembeddability}, and denote it by
$\cvxeq$, the equivalence relation on $\LO$ induced
by $\trianglelefteq_{\LO}$, that is
\[
{L} \cvxeq {L'} \iff {L}
\trianglelefteq_{\LO} {L'} \mbox{ and } {L'} \trianglelefteq_{\LO}
{L} .
\]

Clearly, if \({{L}\iso {L'}}\) then
\({{L} \cvxeq {L'}}\).  The converse
implication does not hold, as witnessed by $\z \o$ and
$\o + \z \o$.

Finally, notice that if \(L \cvxeq L'\) then \(L \equiv_{\LO} L'\),
where \(\equiv_{\LO}\) is the equivalence relation of biembeddability
on \(\LO\) induced by \( \preccurlyeq_{\LO} \). 
The converse is not true: the linear orders of the form \(\boldsymbol{k}\eta\), for \(k>0\), belong to the same \(\equiv_{\LO}\)-equivalence class, but they are pairwise $\trianglelefteq_{\LO}$-incomparable.

\subsection{Combinatorial properties}\label{sec:comb_prop}

In this section we explore the combinatorial properties of convex
embeddability, pointing out several differences between \(\cvx\) and the
embeddability relation \(\preccurlyeq_{\LO}\) on \(\LO\). For example, we show 
that $\cvx$ has antichains of size the continuum and chains of order type \( (\R, {\leq} ) \) (hence descending and ascending chains of arbitrary countable length),
that well-orders are unbounded with respect to \( \cvx \) (hence the unbounding number of \( \cvx \) is \( \aleph_1 \)), 
that \(\cvx\) has dominating number \( 2^{\aleph_0} \) (thus in particular there is no \( \cvx \)-maximal element), and
that all bases for \( \cvx \) have maximal size \( 2^{\aleph_0} \). 
This is in stark contrast with the fact that \(\preccurlyeq_{\LO}\) is a wqo (and hence has neither infinite antichains nor infinite descending chains), that \(\eta\) is the maximum with respect to \(\preccurlyeq_{\LO}\) (hence there are no \( \preccurlyeq_{\LO} \)-unbounded sets and the dominating number of \( \preccurlyeq_{\LO} \) is \( 1 \)), and that \( \{ \o, \o^* \} \) is a two-elements basis for  \(\preccurlyeq_{\LO}\).

Applying Proposition \ref{cond_classes_isom} and recalling that a convex embedding
\( f \colon L \to L' \) is just an isomorphism between \( L \) and a convex subset of \( L' \), 
we easily obtain the
following useful fact.

\begin{prop}\label{cond_classes}
Let \(L, L'\) be arbitrary linear orders. If 
\(L \trianglelefteq L'\) via some convex embedding \( f \colon L \to L' \), then the restriction of \(f\) witnesses
\(c_F^L(\ell) \cong c_F^{L'}(f(\ell)) \cap f(L)\) and hence
\(|c_F^L(\ell)|=|c_F^{L'}(f(\ell)) \cap f(L)| \leq |c_F^{L'}(f(\ell))| \) for every \(\ell \in
L\). Moreover, \(f(c_F^L(\ell)) = c_F^{L'}(f(\ell))\) for every
\(\ell \in L\), except for the first and last condensation classes of \( L \) (if
they exist). Finally, \(L_F \trianglelefteq L'_F\) via the well-defined map \( c^L_F(\ell) \mapsto c^{L'}_F(f(\ell)) \).
\end{prop}

Using the previous proposition and arguing as in the proof of Lemma \ref{lem:isom_zetaL}, it is straightforward to prove that \(\z L \trianglelefteq \z L'\) if and only if \(L \trianglelefteq L'\).

Given a map \(f \colon \Q \to \Scat\), let \(\eta_f \in \LO \) be (an isomorphic copy on \( \N \) of) the linear order \( \eta_f = \sum_{q \in \Q} f(q) \).

\begin{lem}\label{open_int}
There is an embedding from the partial order
\((\Cvx(\R),\subseteq)\) into \((\LO,\cvx)\), where \(\Cvx(\R)\) is
the set of the open intervals of \(\R\).
\end{lem}

\begin{proof}
Consider an injective map
\(f \colon \Q \to \{\boldsymbol{n} \mid n \in \N \setminus
\{0\}\}\) and consider the resulting linear order \(\eta_f \). Notice that for each \((\ell,q) \in \eta_f\),
\(|c_F^{\eta_f}(\ell,q)|= |f(q)|\) is finite. Moreover if \(q\) and
\(q'\) are distinct rational numbers then
\(|c_F^{\eta_f}(\ell,q)| \neq |c_F^{\eta_f}(\ell',q')|\) for every
\(\ell \in f(q)\) and \(\ell' \in f(q')\) by injectivity of \( f \).

An element of \(\Cvx(\R)\) is of the form \((x,y)\) where
$x \in \{-\infty\} \cup \R$ and $y \in \R \cup \{+\infty\}$ with
$x<y$. For such $(x,y)$ we define the linear order \( L_{(x,y)} \cong \sum_{q \in \Q \cap (x,y)} f(q) \) as the restriction of \( \eta_f \) to \( \{(\ell,q) \in \eta_f \mid q \in \Q \cap (x,y)\} \),
which
is a convex subset of \(\eta_f\) with no first and last condensation
class.

We show that, after canonically coding each \( L_{(x,y)} \) as an element of \( \LO \), 
the map \((x,y) \mapsto L_{(x,y)}\) is an embedding of
the partial order \((\Cvx(\R),\subseteq)\) into \((\LO,\cvx)\).
Clearly, if \((x,y) \subseteq(x',y')\), then
\(L_{(x,y)} \csube L_{(x',y')}\) and in particular we have
\(L_{(x,y)} \cvx L_{(x',y')}\).
Vice versa, take \((x,y),(x',y') \in \Cvx(\R)\), with
\((x,y) \not\subseteq (x',y')\) and fix
\(q \in (x,y) \setminus (x',y')\). The condensation class of
\((0,q)\) in \(L_{(x,y)}\) has cardinality \(f(q)\) and, by
injectivity of \(f\), no condensation class in \(L_{(x',y')}\) has
the same cardinality. Since there is no first and last condensation
class in \(L_{(x,y)}\), we get
\(L_{(x,y)} \ncvx L_{(x',y')}\)  by Proposition \ref{cond_classes}.
\end{proof}

\begin{prop}\label{prop1}
$\cvx$ has chains of order type \( (\R, {<} ) \), as well as antichains of size $2^{\aleph_0}$.
\end{prop}

\begin{proof}
This is immediate from Lemma \ref{open_int} and the fact that
\((\Cvx(\R), {\subseteq})\) has the same properties: consider e.g.\ the families \( \{ (x,+\infty)  \mid x \in \R \}\)
and \( \{ (x,x+1) \mid x \in \R \}\), respectively.
\end{proof}

Let \( \mathfrak{b}(\cvx) \) be the \textbf{unbounding number} of \( \cvx \), i.e.\ the smallest size of a family \( \mathcal{F} \subseteq \LO \) which is unbounded with respect to \( \cvx \).
Using
infinite (countable) sums of linear orders, one can easily prove that \( \mathfrak{b}(\cvx) > \aleph_0 \). The next result thus shows that \( \mathfrak{b}(\cvx) \) attains the smallest possible value.

\begin{prop}\label{prop_WO}
\(\WO\) is a maximal \(\omega_1\)-chain without an upper bound in \(\LO\)
with respect to \(\cvx\). Hence \( \mathfrak{b}(\cvx) = \aleph_1 \).
\end{prop}

\begin{proof}
Fix \(L \in \LO\) and for every \(n \in L\) define
\[
\alpha_{n,L} = \sup \{\ot(L')\mid L' \text{ is a well-order},\
L'\csube L, \text{ and } n = \min L'\}.
\]
Notice that \( \alpha_{n,L} \) is actually attained by definition of \( \csube \). Therefore, \( \alpha_{n,L} < \omega_1 \) because \( L \) is countable. 
Let
\(\alpha_L = \sup_{n \in L} \alpha_{n,L}< \o_1\).  
By construction, if \( L' \csube L\) and \( L \) is well-ordered, then \( \ot(L') \leq \alpha_L \),
thus
\(\boldsymbol{\alpha}_L + \boldsymbol{1} \ntrianglelefteq L\).  
Since \( L \) was arbitrary, we showed that for every \(L \in \LO\) there
exists \(L' \in \WO\) such that \(L'\ntrianglelefteq L\), i.e.\ that 
\( \WO \) is \( \cvx \)-unbounded in \( \LO \).

Clearly, \( \WO \) is a \( \cvx \)-chain: maximality then follows from unboundedness of \( \WO \), together with the observation that for \( \omega \leq \alpha < \omega_1 \) and \( L \in \LO \), if \(  L \cvx \boldsymbol{\alpha} \) and \( \boldsymbol{\beta} \triangleleft_{\LO} L \) for every \( \beta < \alpha \), then \( L \cong \boldsymbol{\alpha} \).
\end{proof}

An easy consequence of Proposition~\ref{prop_WO} is that no \( L \in \LO \) is a node with respect to \( \cvx \). This will be subsumed by Proposition~\ref{prop:everylobelongstoantichain}.

\begin{cor} \label{cor:nonoteincvx}
For every \( L \in \LO \) there is \( M \in \LO \) which is \( \cvx \)-incomparable with \( L \), i.e.\ \( L \not\trianglelefteq_{\LO} M \) and \( M \not\trianglelefteq_{LO} L \).
\end{cor}

\begin{proof}
If \( L \) is not a well-order, then it is enough to let \( M \in \WO \) be such that \( M \not\trianglelefteq_{\LO} L \) (the existence of such an \( M \) is granted by Proposition~\ref{prop_WO}). If instead \( L \) is a well-order, then it is enough to set \( M = \omega^* \).
\end{proof}

Another easy consequence of Proposition~\ref{prop_WO} is that \( \cvx \) has no maximal element. In fact, much more is true.

\begin{cor}\label{cor_no_max}
Every \( L \in \LO \) is the bottom element of a \( \cvx \)-unbounded chain of length \( \omega_1 \).
\end{cor}

\begin{proof}
For every \( \beta < \omega_1 \) set \( L_\beta = L+\boldsymbol{\beta} \) (in particular, \( L_0 = L \)), and consider the (not necessarily strictly) \( \cvx \)-increasing sequence \( \langle L_\beta \mid \beta < \omega_1 \rangle \). Since \( \boldsymbol{\beta} \cvx L_\beta \) for every \( \beta < \omega_1 \), the above sequence is \( \cvx \)-unbounded by Proposition~\ref{prop_WO}. Moreover, for every \( \beta < \omega_1 \) there is \( \beta' < \omega_1 \) such that \( L_\beta \triangleleft_{\LO} L_{\beta'} \). Indeed, it is enough to set \( \beta' = \alpha_{L_\beta} + 1 \), where \( \alpha_{L_\beta} \) is as in the proof of Proposition~\ref{prop_WO}: then \( \boldsymbol{\beta}' \not\trianglelefteq_{\LO} L_\beta \), and thus also \( L_{\beta'} \not\trianglelefteq_{\LO} L_\beta \).
This easily implies that \( \langle L_\beta \mid \beta < \omega_1 \rangle \) contains a strictly \( \cvx \)-increasing cofinal (hence \( \cvx \)-unbounded in \( \LO \)) chain of length \( \omega_1 \) beginning with \( L_0 \), as desired.
\end{proof}

We say that a collection \( \mathcal{B} \) of (infinite) linear orders on \( \N \) is a \textbf{basis} for \( \cvx \) if for every \( L \in \LO \) there is \( L' \in \mathcal{B} \) such that \( L' \cvx L \).
The next result shows that each basis with respect to \( \cvx \) is as large as possible.

\begin{prop} \label{prop:basisforcvx}
\begin{enumerate-(a)}
\item \label{prop:basisforcvx-1}
There are \(2^{\aleph_0}\)-many \( \cvx \)-incomparable \(\cvx\)-minimal elements in \( \LO \). In particular, if \( \mathcal{B} \) is 
a basis for \( \cvx \) then \( |\mathcal{B}| = 2^{\aleph_0} \).
\item \label{prop:basisforcvx-2}
There is a \( \cvx \)-decreasing \( \omega \)-sequence in \( \LO \) which is not \( \cvx \)-bounded from below.
\end{enumerate-(a)}
\end{prop}

\begin{proof}
\ref{prop:basisforcvx-1}
Consider an infinite subset \(S\subseteq \N\). Let
\(f_S\colon \Q \to \{ \boldsymbol{n} \mid n \in S \} \) be a  map such that
\[
\forall q,q'(q<q'\rightarrow \forall n \in S\ \exists q''(q<q''<q'
\wedge f_S(q'')= \boldsymbol{n})),
\]
so that in particular \( f_S \) is surjective, and consider the linear order \(\eta_{f_S} \). Let \( q  < q' \) be arbitrary rational numbers. By a back-and-forth argument
on the condensation classes, it is easy to see that by choice of \( f_S \) the linear order \( \eta_{f_S} \) is isomorphic to the restriction \( \eta_{f_S} \rest (q,q') \iso \sum_{q'' \in \Q \cap (q,q')} f_S(q'') \)
of \( \eta_{f_S} \) to 
its convex subset \( \{ ( \ell, q'') \in \eta_{f_S} \mid q < q'' < q' \} \). This implies that each \( \eta_{f_S} \) is \( \cvx \)-minimal, because by density of \( \eta \) and finiteness of the condensation classes of \( \eta_{f_S} \),  any infinite convex subset of \( \eta_{f_S} \) contains some \( \eta_{f_S} \rest (q,q') \). 
Finally, by the choice of \( f_S \) for every \( n \in S \) there are densely many condensation classes in \( (\eta_{f_S})_F \)  of size exactly \( n \). Thus if \( S \neq S' \) we have \( \eta_{f_S} \not\trianglelefteq_{\LO} \eta_{f_{S'}} \) and \( \eta_{f_{S'}} \not\trianglelefteq_{\LO} \eta_{f_{S}} \) by Proposition
\ref{cond_classes}, as desired.
%
%

\ref{prop:basisforcvx-2}
Consider the family \( \{ L_{(n,+\infty)} \mid n \in \N \}\), where \( L_{(n,+\infty)} \) is as in the proof of Lemma~\ref{open_int}. It is a strictly \( \cvx\)-decreasing chain, and we claim that it is \( \cvx \)-unbounded from below. To this aim, it is enough to consider any \( L \in \LO \) with \( L \cvx L_{(0,+\infty)}\), and show that \( L \ncvx L_{(m,+\infty)} \) for some \( m \in \N \). 
Since \( L \cvx L_{(0,+\infty)}\), all the condensation classes of \( L \) are finite by Proposition~\ref{cond_classes}. Let \( \ell \in L \) be such that \(c^L_F(\ell) \) is not the minimum or the maximum of \( L_F \), and let \( q \in \Q \) be such that \( f(q) = |c^L_F(\ell) | \), where \(f \colon \Q \to \{\boldsymbol{n} \mid n \in \N \setminus
\{0\}\}\) is the function used to define the linear orders \( L_{(n,+\infty)} \). Let \( m \in \N \) be such that \( q < m \). Then \( L \ncvx L_{(m,+\infty)} \) because otherwise by choice of \( \ell \) the latter would have a condensation class of size \( f(q) \) by Proposition~\ref{cond_classes}, which is impossible by choice of \( m \) and the fact that \( f \) is an injection.
\end{proof}

Proposition~\ref{prop:basisforcvx} allows us to considerably improve
Corollary~\ref{cor:nonoteincvx} as follows.

\begin{prop} \label{prop:everylobelongstoantichain} 
Every \( \cvx \)-antichain is contained in a \( \cvx \)-antichain of size \( 2^{\aleph_0} \).
In particular, there are no maximal \( \cvx\)-antichains of size smaller than \(  2^{\aleph_0} \), and every \( L \in \LO \) belongs to a \( \cvx \)-antichain of size \( 2^{\aleph_0} \).
\end{prop}

\begin{proof}
Let \( \mathcal{B} \) be a \( \cvx \)-antichain and assume that $|\mathcal{B}|<2^{\aleph_0}$ (otherwise the statement is trivial). Consider the antichain \( \mathcal{A} = \{ \eta_{f_S} \mid S \subseteq \N \wedge S \text{ is infinite} \} \) of size \( 2^{\aleph_0} \) from Proposition~\ref{prop:basisforcvx}. From \( \cvx \)-minimality of \( \eta_{f_S} \) it follows that \(\mathcal{B} \cup (\mathcal{A} \setminus \bigcup_{L \in \mathcal{B}} \{ K \in \mathcal{A} \mid K \cvx L \}) \) is a \( \cvx \)-antichain. To show that this antichain has size \( 2^{\aleph_0} \) it suffices to show that 
\begin{claim} \label{claim:everylobelongstoantichain}
  \( \{ K \in \mathcal{A} \mid K \cvx L \} \) is countable for every $L \in \LO$,
\end{claim}
\noindent
so that \( |\bigcup_{L \in \mathcal{B}} \{ K \in \mathcal{A} \mid K \cvx L \}| \leq \aleph_0 \cdot |\mathcal{B}| = \max \{ \aleph_0, |\mathcal{B}| \} < 2^{\aleph_0} \).

To prove the claim, suppose that \( S \subseteq \N \) is such that \( \eta_{f_S} \cvx L \), so that without loss of generality we can write \( L = L_l + \eta_{f_S} + L_r \). If \( f \) were a convex embedding of \( \eta_{f_{S'}} \) into \( L \) with \( f(\eta_{f_{S'}}) \cap \eta_{f_S} \neq \emptyset \), then by density of \( \eta \) and finiteness of the condensation classes of \( \eta_{f_{S'}} \) there would be rationals \( q < q' \) such that \( f(\eta_{f_{S'}} \rest (q,q')) \subseteq \eta_{f_S} \), and since \( \eta_{f_{S'}} \cong \eta_{f_{S'}} \rest (q,q') \) we would get \( \eta_{f_{S'}} \cvx \eta_{f_S} \). Thus if \( S \neq S' \), then \( f(\eta_{f_{S'}}) \cap \eta_{f_S} = \emptyset \). Since \( L \) is countable, this means that there are only countably many distinct \( S \subseteq \N \) for which \( \eta_{f_S} \cvx L \) can hold.

Finally, the additional part of the statement follows by viewing $L \in \LO$ as the element of an antichain of size $1$.
\end{proof}

We say that a collection \(\mathcal{F} \subseteq \LO\) is
a dominating family with respect to \(\cvx\) if and only if for every
\(L \in \LO\) there exists \(L' \in \mathcal{F}\) such that
\(L \cvx L'\). Let \( \mathfrak{d}(\cvx) \) be the \textbf{dominating number} of \( \cvx \), i.e.\ the least size of a dominating family with respect to \( \cvx \).
The next proposition shows that \( \mathfrak{d}(\cvx) \) is as large as it can be.

\begin{prop}\label{prop:dom_fam}
\( \mathfrak{d}(\cvx) = 2^{\aleph_0}\).
\end{prop}

\begin{proof}
Consider again the antichain \( \mathcal{A} = \{ \eta_{f_S} \mid S \subseteq \N \} \) from the proof of Proposition~\ref{prop:basisforcvx}. 
If \( \mathcal{F} \) were a dominating family with respect to \( \cvx \) of size \( \kappa < 2^{\aleph_0} \), then by \( |\mathcal{A}| = 2^{\aleph_0} \) there would be \( M \in \mathcal{F} \) such that \( \{ K \in \mathcal{A} \mid K \cvx M \} \) is uncountable, contradicting Claim~\ref{claim:everylobelongstoantichain}.
\end{proof}

\subsection{Complexity with respect to Borel reducibility}\label{sec:complexity of cvx}

At the beginning of Section \ref{sec:cvx emb} we introduced the equivalence
relation \(\cvxeq\) of convex biembeddability on \(\LO\), observing
that it is different from both isomorphism and biembeddability. We
now focus on determining the complexity of \(\cvxeq\) with respect to
Borel reducibility.

\begin{thm} \label{thm:red_iso_cvxeq}
The map \( \varphi \) sending a linear order \( L \) to \( \varphi(L) = \boldsymbol{1}+\z L + \boldsymbol{1} \) is such that
\begin{enumerate-(a)}
\item \label{thm:red_iso_cvxeq-a}
\( {L \cong L'} \iff { \varphi(L) \cong \varphi(L')} \iff {\varphi(L) \mathrel{\underline{\bowtie}} \varphi(L')}\iff {\varphi(L) \trianglelefteq \varphi(L')} \);
\item \label{thm:red_iso_cvxeq-b}
\( |\varphi(L)| = \max \{ \aleph_0, |L| \} \).
\end{enumerate-(a)}
\end{thm}

\begin{proof} 
We claim that \(\varphi\) reduces \(\iso\) to \(\cvxeq\). 
The second part is obvious, so let us concentrate on the first one.
It
is immediate that if \(L \cong L'\) then
\(\varphi(L) \cong \varphi(L')\) and hence
$\varphi(L) \mathrel{\underline{\bowtie}} \varphi(L')$, while $\varphi(L) \mathrel{\underline{\bowtie}} \varphi(L')$ clearly implies $\varphi(L) \trianglelefteq \varphi(L')$.

Let now $f$ witness
$\varphi(L)\trianglelefteq \varphi(L')$. The only elements of $\varphi(L)$
and $\varphi(L')$ without immediate successor and immediate
predecessor are their minimum and maximum, respectively. Therefore, we must have
$f(\min \varphi(L))=\min \varphi(L')$ and
$f(\max \varphi(L))=\max \varphi(L')$. Hence $f$ is also surjective (hence an isomorphism),
and $f\restriction (\z L)$ witnesses $\z L \cong \z L'$. Thus $L \cong L'$ by
Lemma \ref{lem:isom_zetaL}.
\end{proof}	

Noticing that when restricted to \( \LO \) the map from Theorem~\ref{thm:red_iso_cvxeq} is Borel, we immediately get

\begin{cor}\label{red_iso_cvxeq}
${\iso} \leq_B {\cvxeq}$.
\end{cor}

The main question now becomes whether \(\cvxeq \leq_B
{\iso}\). This is still open and the answer is not obvious because e.g.\ it is not even clear if \(\cvxeq\) is induced by a Borel action of \(S_\infty\). We now embark in a deeper analysis of \(\cvxeq\), leading at least to
\(\cvxeq \leq_{\text{\scriptsize \textit{Baire}}} {\iso}\).

In the spirit of the definition of convex embeddability and recalling
Remark \ref{def cvx}, we introduce the following notions.

\begin{defn}
Let $L$ be a linear order. We say that 
\begin{enumerate-(1)}
\item $L$ is \textbf{right compressible} if $L = L' + L_{r}$, with
$L' \cong L$ and $L_{r} \neq\emptyset$;
\item $L$ is \textbf{left compressible} if $L= L_{l} + L'$, with
$L' \cong L$ and $L_{l} \neq\emptyset$;
\item $L$ is \textbf{bicompressible} if it is both left compressible
and right compressible,
\item $L$ is \textbf{incompressible} if it is neither left nor right
compressible.
\end{enumerate-(1)}
\end{defn}

Notice that the set of right compressible linear orders is
invariant with respect to isomorphism. The same holds for the
set of left compressible linear orders, the set of bicompressible
linear orders, and that of incompressible linear orders.

It is clear that \( \omega^* \) is right compressible but not left
compressible, \( \omega \) is left compressible but not right
compressible, \(\omega + \omega^*\) and \( \eta \) are bicompressible,
and \( \z \) is incompressible.

The following characterizations of the above notions turn out to be useful.

\begin{lem}\label{strong_right}
Let $L$ be a linear order. Then
\begin{enumerate-(a)}
\item \label{strong_right-a}
$L$ is right compressible if and only if
$L = L_{l} + \widetilde{L} + L_{r}$, with $\widetilde{L} \cong L$
and $L_{r} \neq\emptyset$.
\item \label{strong_right-b}
$L$ is left compressible if and only if
$L = L_{l} + \widetilde{L} + L_{r}$, with $\widetilde{L} \cong L$
and $L_{l} \neq\emptyset$.
\item \label{strong_right-c}
$L$ is bicompressible if and only if
$L=L_{l} + \widetilde{L} + L_{r}$, with $\widetilde{L} \cong L$
and $L_l,L_{r} \neq\emptyset$.
\end{enumerate-(a)}
\end{lem}

\begin{proof}
\ref{strong_right-a} For the non trivial direction, suppose that
$L = L_{l} + \widetilde{L} + L_{r}$, with $L \cong \widetilde{L}$
via some $f \colon L \to \widetilde{L}$ and $L_{r} \neq\emptyset$. Let
$M_0 = f(L_r)\subseteq \widetilde{L}$ and for every $n \in \N$
define $M_{n+1} = f(M_n)\subseteq \widetilde{L}$. Let
$M = \bigcup_{n \in \N} M_n \subseteq \widetilde{L}$ and note
that $f\upharpoonright (M + L_r)\colon  M + L_r \to M$ is an
isomorphism. Then the map $g\colon L\to L_l + \widetilde{L}$
defined by
\[g(x) = 
\begin{cases}
f(x) & \mbox{if } x \in M + L_r,\\
x, & \mbox{otherwise }
\end{cases}\] 
is an isomorphism witnessing
\( L \cong L_l + \widetilde{L}\).  Thus, we can write
$L = L' + L_{r}$, with $L' = L_l + \widetilde{L}$.

\ref{strong_right-b} is similar to \ref{strong_right-a}.

\ref{strong_right-c} If $L= L_{l} + L' + L_{r}$ with $L' \cong L$ and $L_l,L_{r} \neq\emptyset$, then by \ref{strong_right-a} and
\ref{strong_right-b} we immediately obtain that \(L\) is
bicompressible. Conversely, suppose that \(L\) is
bicompressible. Since $L$ is left compressible, then $L=L_{l} + L'$,
with $L' \cong L$ and $L_{l} \neq\emptyset$. Since $L'\cong L$ is
right compressible, we can write $L'= \widetilde{L} + L_{r}$, with
$\widetilde{L} \cong L'$ and $L_{r} \neq\emptyset$. Hence,
$L=L_{l} + \widetilde{L} + L_{r}$, with $\widetilde{L} \cong L' \cong L$ and $L_l,L_r\neq\emptyset$.
\end{proof}

We denote by $\LO_{r}\subseteq \LO$ the set of (codes for) right
compressible linear orders on \(\N\), and by $\LO_{l}\subseteq \LO$
the set of (codes for) left compressible linear orders on \(\N\).
Note that $\LO_{r}=\{L \in \LO \mid L^* \in \LO_{l}\}$, and vice
versa. Moreover, each of the four sets
\begin{equation}\label{subsets of LO}
\LO \setminus (\LO_{r}\cup \LO_{l}) \qquad \quad \LO_{l}\setminus \LO_{r} \qquad \quad  \LO_{r}\setminus \LO_{l} \qquad \quad \LO_{r}\cap \LO_{l}
\end{equation} 
is closed under isomorphism. The next proposition shows that they are also closed under \(\cvxeq\).

\begin{prop}\label{inv_wrt_cvxeq} 
If \( L \) is a right compressible linear order and \( L' \mathrel{\underline{\bowtie}} L \) (which implies \( |L'| = |L| \)), then \( L' \) is right compressible as well. Similarly, if \( L' \mathrel{\underline{\bowtie}} L \) and  \( L \) is left compressible (respectively: bicompressible, incompressible), then so is \( L' \).

In particular, the four subsets \(\LO \setminus (\LO_{r}\cup \LO_{l})\),
\(\LO_{l}\setminus \LO_{r},\LO_{r}\setminus \LO_{l}\) and
\(\LO_{r}\cap \LO_{l}\) are invariant with respect to
${\cvxeq}$.
\end{prop}

\begin{proof}
It is clearly enough to consider the case of right compressible linear orders.
Since \( L \) is right compressible, then
$L = \widetilde{L} + L_{r}$ for some $\widetilde{L} \cong L$ and $L_{r} \neq\emptyset$. 
Let \( f \colon L' \to \widetilde{L} \) and \(g \colon L \to L' \) be convex embeddings witnessing \( L' \trianglelefteq \widetilde{L} \) and \( L \trianglelefteq L' \), respectively, so that \( \widetilde{L} = \widetilde{L}_l + f(L') + \widetilde{L}_r \) and \( L' = L'_l + g(L) + L'_r \). Then 
\begin{align*}
L' & = L'_l + g(L) + L'_r \\
& = L'_l + g(\widetilde{L}) + g(L_r) + L'_r \\
& = L'_l + g(\widetilde{L}_l) + g(f(L')) + g(\widetilde{L}_r) + g(L_r) + L'_r.
\end{align*}
Since \( g(f(L')) \cong L' \) and 
\( g(\widetilde{L}_r) + g(L_r) + L'_r \supseteq g(L_r) \neq \emptyset \),
by
Lemma \ref{strong_right} we have $L' \in \LO_r$, as desired.
\end{proof}



We are now ready to go back to the study of the complexity of convex
biembeddability. We can prove that the
restrictions of \(\cvxeq\) to each of the four sets in \eqref{subsets
of LO}, which we denote by 
$\underline{\bowtie}_{\tiny \LO \setminus (\LO_{r}\cup \LO_l)}$, 
$\underline{\bowtie}_{\tiny \LO_{l}\setminus \LO_r}$,
$\underline{\bowtie}_{\tiny \LO_{r}\setminus \LO_l}$ and 
$\underline{\bowtie}_{\tiny \LO_{r}\cap \LO_l}$, respectively,
are Borel bireducible with \(\iso\).

To this aim, we first observe that the map \( \varphi_0 = \varphi \) from Theorem~\ref{thm:red_iso_cvxeq} reduces isomorphism to convex biembeddability restricted to incompressible linear orders, and that suitable variations of it do the same job but with left compressible (respectively, right compressible, bicompressible) linear orders.

\begin{prop} \label{prop:thm 3.2.3}
Given a linear order \( L \), set
\begin{align*}
\varphi_0(L) & = \boldsymbol{1}+\z L + \boldsymbol{1} \\
\varphi_1(L) & = \eta + \z L + \boldsymbol{1}\\
\varphi_2(L) & = \boldsymbol{1} + \z L + \eta \\
\varphi_3(L) & = \eta + \z L + \eta.
\end{align*}
Then \( \varphi_0(L) \) is incompressible, \( \varphi_1(L) \) is left compressible but not right compressible, \( \varphi_2(L) \) is right compressible but not left compressible, and \( \varphi_3(L) \) is bicompressible. 

Moreover, Theorem~\ref{thm:red_iso_cvxeq} is still true when \( \varphi \) is replaced by any of the above \( \varphi_i \)'s.
\end{prop}

\begin{proof}
As the minimum and
the maximum of \(\varphi_0(L)\) are the only elements without
immediate predecessor and successor, respectively, we have that \(\varphi_0(L)\) is
not isomorphic to any of its proper convex subsets, i.e.\ it is incompressible.
Hence we are done with \( \varphi_0 \) by Theorem~\ref{thm:red_iso_cvxeq}.

Using a similar argument, one easily sees that \( \varphi_1(L) \) is not right compressible.
Indeed, any convex embedding \( f \) of \( \varphi_1(L) \) into itself cannot send \( \max \varphi_1(L) \) into \( \zeta L \) (by the argument in the previous paragraph)
and cannot send it into \( \eta \) either (because otherwise \( f(\zeta L) \subseteq \eta \), which is clearly impossible). 
On the other hand, \( \varphi_1(L) \) is trivially left compressible
because one can map \( \eta \) onto any of its (proper) final segments. 
Obviously \( |\varphi_1(L)| = \max \{ \aleph_0, |L| \} \),
\( {L \cong L'} \Rightarrow {\varphi_1(L) \cong \varphi_1(L')} \), \( {\varphi_1(L) \cong \varphi_1(L')} \Rightarrow {\varphi_1(L) \mathrel{\underline{\bowtie}} \varphi_1(L') } \), and also \( {\varphi_1(L) \mathrel{\underline{\bowtie}} \varphi_1(L') } \Rightarrow {\varphi_1(L) \trianglelefteq \varphi_1(L')} \), so it remains to prove that if \( {\varphi_1(L) \trianglelefteq \varphi_1(L') } \) then \( {L \cong L'} \). Let \( f \colon \varphi_1(L) \to \varphi_1(L') \) be a convex embedding. Since the elements of \( \eta \) are the unique non-maximal points without immediate predecessor and immediate successor (both in \( \varphi_1(L) \) and \( \varphi_1(L') \)), then \( f(\eta) \subseteq \eta \). Similarly, since the elements of \( \z L \) and \( \z L' \) are the only elements having both an immediate predecessor and an immediate successor, then \( f(\z L) \subseteq \z L' \). Moreover, the maximal element \( \boldsymbol{1} \) has no immediate predecessor, which forbids \( f(\boldsymbol{1}) \in \z L \), and we cannot have \( f(\boldsymbol{1}) \in \eta \) because otherwise \( f(\z L) \subseteq \eta \): thus \( f(\boldsymbol{1}) = \boldsymbol{1} \). Since the range of \( f \) is convex, it then follows that \( f (\z L) = \z L' \), hence \( \z L \cong \z L' \) and thus \( L \cong L' \) by Lemma~\ref{lem:isom_zetaL}.

The cases of \( \varphi_2(L)\) and \( \varphi_3(L) \) are similar.
\end{proof}

When restricted to \( \LO \), the functions \( \varphi_i \) are clearly Borel, thus we obtain:

\begin{cor}\label{thm 3.2.3}
The isomorphism relation \( \iso \) is Borel reducible to any of 
\( \underline{\bowtie}_{\tiny \LO \setminus (\LO_{r}\cup \LO_l)} \), 
\( \underline{\bowtie}_{\tiny \LO_{l}\backslash \LO_r} \),
\( \underline{\bowtie}_{\tiny \LO_{r}\backslash \LO_l} \), and
\( \underline{\bowtie}_{\tiny \LO_{r}\cap \LO_l} \).
\end{cor}

Notice that the ranges of the four reductions used in the proof of Corollary~\ref{thm 3.2.3} are all Borel, and that on such ranges isomorphism and convex biembeddability coincide.

\begin{thm}\label{red_cvxeq_jump_iso}
\begin{enumerate-(a)}
\item\label{red_cvxeq_jump_iso-i}
On the set $\LO \setminus (\LO_{r}\cup \LO_{l})$ the relations $\cvxeq$ and $\iso$ coincide, so that $\underline{\bowtie}_{\tiny \LO \backslash (\LO_{r}\cup \LO_l)}$ is Borel reducible to $\iso$ via the identity map.

\item \label{red_cvxeq_jump_iso-ii}
Each of $\underline{\bowtie}_{\tiny \LO_{r}\cap \LO_l}$, $\underline{\bowtie}_{\tiny \LO_{l}\setminus \LO_r}$, and
$\underline{\bowtie}_{\tiny \LO_{r}\setminus \LO_l}$
is Borel reducible to $(\iso)^+$, and thus to \( \iso \).
\end{enumerate-(a)}
\end{thm}

\begin{proof}
\ref{red_cvxeq_jump_iso-i} 
Let $L,L' \in \LO \setminus (\LO_{r}\cup \LO_{l})$. It is obvious that if $L \iso L'$ then $L \cvxeq L'$. For the other direction, assume \(L \cvxeq L'\) and let $f$ and $g$ be convex embeddings
witnessing $L \trianglelefteq_{\LO} L'$ and $L' \trianglelefteq_{\LO} L$, respectively. Then
$L' = L'_{l} + f(L) + L'_{r}$ and $L= L_{l} + g(L'_l)+ g(f(L))+ g(L'_r) + L_{r}$. Since $L$ is incompressible and \( g(f(L)) \cong L \) we have $L_{l} + g(L'_l) = g(L'_r) + L_{r} = \emptyset$ and hence $L'_{l} = L'_{r} = \emptyset$, showing that $f$ is an isomorphism.

\ref{red_cvxeq_jump_iso-ii}
We start by considering the case of \( \underline{\bowtie}_{\tiny \LO_{r}\cap \LO_l} \). 
Let $\varphi_{r+l}\colon  \LO \setminus [\z,\omega,\omega^*]_{\cong}
\to \LO^\N$ be a Borel map such that
$\varphi_{r+l}(L)$ is an enumeration (possibly with repetitions) of all the \emph{infinite}
subsets of \(L\) of the form $ [n,m]_L$. Since we are omitting the isomorphism types of \( \z \), \( \omega \), and \( \omega^* \) the map is well-defined, i.e.\ for each \( L \) in its domain there is at least one infinite interval \( [n,m]_L \), and clearly \( \LO_l \cap \LO_r \subseteq \dom(\varphi_{r+l}) \). By the same reason, its domain is Borel because we are omitting finitely many \( \iso \)-classes, which are Borel themselves. We claim that for all \( L , L' \in \LO_l \cap \LO_r \)
\[  
L \cvxeq L' \iff \varphi_{r+l}(L) \mathrel{(\iso)^+} \varphi_{r+l}(L'),
\]
so that any Borel extension of \( \varphi_{r+l} \) to \( \LO \) witnesses \( {\underline{\bowtie}_{\tiny \LO_{r}\cap \LO_l}} \leq_B {(\iso)^+} \), and hence \( {\underline{\bowtie}_{\tiny \LO_{r}\cap \LO_l}} \leq_B {\iso} \) by Theorem~\ref{cor 2}.

Assume first that
\(L \cvxeq L'\), and let $f$ be a convex embedding witnessing
$L\trianglelefteq_{\LO} L'$. Given any infinite $[n,m]_L$, we have
$[n,m]_L \iso [f(n),f(m)]_{L'}$, so that in particular the latter is infinite and appears among the linear orders in \( \varphi_{r+l}(L') \). Symmetrically, if $g$ is a convex embedding witnessing
$L' \trianglelefteq_{\LO} L$, then for every infinite $[n,m]_{L'}$ we have
$[n,m]_{L'} \cong [g(n),g(m)]_L$. It follows that
$\varphi_{r+l}(L) \mathrel{(\iso)^+} \varphi_{r+l}(L')$.

Conversely, observe that since $L \in \LO_{l}\cap \LO_{r}$ then by
Lemma~\ref{strong_right} we have $L = L_l+ \widetilde{L} + L_r$,
with \(\widetilde{L}\cong L\) and both \(L_l\) and \(L_r\)
nonempty. Fix $k \in L_l$ and $m \in L_r$. Then
\(\widetilde{L} \csube [k,m]_L\), and hence
$L \trianglelefteq [k,m]_L$ and $[k,m]_L$ is infinite. Thus if
\( \varphi_{r+l}(L) \mathrel{(\iso)^+} \varphi_{r+l}(L')\),
there are $k', m' \in L'$ such that $[k,m]_L\cong
[k',m']_{L'}$. But then \( L \cvx L' \) because
$L \trianglelefteq [k,m]_L \cong [k',m']_{L'}\trianglelefteq L'$.
The argument to show that if \( \varphi_{r+l}(L) \mathrel{(\iso)^+} \varphi_{r+l}(L')\) then $L'\trianglelefteq_{\LO} L$ is symmetric.

We now move to the case of $\underline{\bowtie}_{\tiny \LO_{l}\setminus \LO_r}$. Let
$\varphi_{l}\colon \LO\setminus [\omega^*]_{\cong} \to
\LO^\N$ be a Borel map such that $\varphi_l(L)$ is an
enumeration of all the \emph{infinite} subsets of \(L\) of the
form $[n,+\infty)_L$, which is well-defined on all \( L \not\cong \o^* \) and such that \( \LO_l \setminus \LO_r \subseteq \dom(\varphi_l)\). Arguing as above, it is enough to show that for all $L,L' \in \LO_{l}\setminus \LO_{r}$
\[  
L \cvxeq L' \iff \varphi_{l}(L) \mathrel{(\iso)^+} \varphi_{l}(L').
\]

For the forward direction, let
$f$ and \(g\) be convex embeddings witnessing $L\trianglelefteq_{\LO} L'$ and
$L'\trianglelefteq_{\LO} L$,  respectively.
We first show that \( f(L) \) is a final segment of \( L' \). Since \( f \) is a convex embedding,
$L'=L'_{l} + f(L) + L'_{r}$ with
\(L'_l\) and $L'_{r}$ possibly empty. Then $L= L_l + g(f(L))+L_r$
with $L_{r}\supseteq g(L'_r)$. Since $g(f(L))\cong L$ and
$L \notin \LO_{r}$, we have \(L_r = \emptyset\) and hence
\(L'_r = \emptyset\), i.e.\ \( L' = L'_l+f(L) \). 
Thus if \( [n, \infty)_L \) is infinite, then \( f([n,\infty)_L) = [f(n),\infty)_{L'} \), so that, being infinite, the latter appears in \( \varphi_l(L') \) and \( [n,\infty)_L \cong [f(n),\infty)_{L'} \). 
Analogously, \( g(L') \) is a final segment of \( L \) because \( L' \notin \LO_r \),  hence for every infinite
$[n,+\infty)_{L'}$, we have
$[n,+\infty)_{L'} \cong [g(n),+\infty)_L$. It follows that
$\varphi_{l}(L) \mathrel{(\iso)^+} \varphi_{l}(L')$.

Conversely, assume that $\varphi_{l}(L) \mathrel{(\iso)^+} \varphi_{l}(L')$. Using \( L \in \LO_l \), let  \(L= L_l + \widetilde{L}\) with
\(L_l \neq \emptyset \) and \( \widetilde{L} \cong L \), and fix
any $m \in L_l$. Then \(\widetilde{L} \csube [m,+\infty)_L\), and
thus the latter, being infinite, appears in \( \varphi_l(L) \) and $L \trianglelefteq [m,+\infty)_L$. Let 
$m' \in L'$ be such that
$[m,+\infty)_L\iso [m',+\infty)_{L'}$: then
$L \trianglelefteq_{\LO} [m,+\infty)_L \iso
[m',+\infty)_{L'}\trianglelefteq_{\LO} L'$. Reversing the role of \( L \) and \( L' \) we get $L'\trianglelefteq_{\LO} L$ and we are done.

The case of $\underline{\bowtie}_{\tiny \LO_{r}\setminus \LO_l}$ is symmetric, with the desired Borel reduction be given by any Borel map $\varphi_{r}\colon  \LO\setminus [\omega]_{\cong} \to
\LO^\N$ such that $\varphi_r(L)$ is an
enumeration of all the \emph{infinite} subsets of \(L\) of the
form $(-\infty,n]_L$. 
\end{proof}

\begin{remark}
The statement and proof of Theorem~\ref{red_cvxeq_jump_iso} can easily be adapted to deal with uncountable linear orders of a given cardinality \( \kappa \). However, since we have no use for this in the present paper, for the sake of simplicity we decided to stick to the countable case.
\end{remark}

If $\LO_{r}$ and $\LO_{l}$ were  Borel subsets of $\LO$, then we could glue the reductions from the proof of Theorem~\ref{red_cvxeq_jump_iso} and
obtain a Borel reduction from the whole \( \cvxeq \) to \( \iso \).
Unfortunately, this is not the case: none of the subclasses of \( \LO \) involved  in Theorem~\ref{red_cvxeq_jump_iso} is Borel.
To prove this, we need the following lemmas.

\begin{lem}\label{incomp_z_cl1}
Let \( \alpha > 0 \).
For any \( z \in \Z^{\a} \) and \( \beta < \a \) there exists
\( \g \) such that \( \beta \leq \g <\a \) and
\( {\Z^{\g} \o} \cvx {[z, +\infty)_{\Z^\a}} \).
\end{lem}

\begin{proof}
We consider the isomorphic copy of \(\Z^\a\) given by Proposition \ref{power_Z}:
\[
\bigg(\sum_{\beta \leq \g < \a} \Z^{\g} \o\bigg)^* + \bigg(\sum_{\beta \leq \g < \a} \Z^{\g} \o\bigg).
\]   
Without loss of generality we can assume \(z \in \sum_{\beta \leq \g < \a} \Z^{\g} \omega \), so that there exists \(\g\) with \(\beta \leq \g < \alpha\) such that \( z \in \Z^{\g} \omega\). Since \( \Z^{\g} \omega \cvx [z, +\infty)_{\Z^{\g} \omega} \cvx [z, +\infty)_{\Z^{\a}} \), this $\g$ works.
\end{proof}

\begin{lem}\label{incomp_z}
For every ordinal \(\alpha>0\), \(\mathbb{Z}^{\alpha}\) is incompressible.
\end{lem}

\begin{proof}
By induction on \(\a>0\). We have already noticed that \(\Z^1 \iso \z\) is
incompressible. Fix \(\a>1 \) and assume that \(\Z^\beta\)
is incompressible for every \(\beta < \a\).

We consider the isomorphic copy of \(\Z^\a\) given by Proposition \ref{power_Z}
with \(\beta=0\):
\[
\bigg(\sum_{\g < \a} \Z^{\g} \o\bigg)^* + \bigg(\sum_{\g < \a} \Z^{\g} \o\bigg).
\] 

We just prove that \( \Z^\a \notin \LO_r \), as \( \Z^\a \notin \LO_l\) can  be proved in a symmetric way. Suppose, towards a contradiction, that \(\Z^\a \in \LO_{r}\) and
let \(f\) be a convex embedding of \(\Z^\a\) into a proper initial
segment of \(\Z^\a\). 
Assume first that \( f(\Z^\alpha) \cap \big(\sum_{\g < \a} \Z^{\g} \o \big) \neq \emptyset \).
Let \(\beta<\a\) be least such that
\(\Z^{\g} \omega \nsubseteq f(\Z^\a) \) for every $\g \geq \beta$. (Such a \( \beta \) exists by the choice of \( f \).) 

\begin{claim} \label{claim:Z^alpha}
\(f(\Z^\a) \cap \Z^{\g} \omega = \emptyset\) for every $\g \geq \beta$, so that \( f(\Z^\alpha) \) is a final segment of
\[
\bigg(\sum_{\gamma < \alpha} \Z^\gamma
\o\bigg)^* + \bigg( \sum_{\gamma < \beta} \Z^\gamma \o \bigg).
\]
\end{claim}

\begin{proof}[Proof of the Claim]
If $\g>\beta$ the convexity of $f(\Z^\a)$ implies immediately \(f(\Z^\a) \cap \Z^{\g} \omega = \emptyset\), so we only need to consider the case $\g=\beta$. Towards a contradiction, assume that \(f(\Z^\alpha)\)
intersects \(\Z^\beta \o\), and using \( f(\Z^\a) \nsupseteq \Z^\beta \omega \) let \(n\) be maximum such that
\(f(\Z^\a) \cap (\Z^\beta \times \{ n \}) \neq \emptyset\). Pick \(z \in \Z^\a\)
such that \(f(z)\in \Z^\beta \times \{ n \}\). By Lemma~\ref{incomp_z_cl1}
there exists \( \beta \leq \gamma < \alpha \) such that
\( \Z^{\g} \omega \cvx [z, +\infty)_{\Z^\a} \) and hence
\( \Z^{\g} \omega \cvx [f(z), +\infty)_{\Z^{\beta} \times \{ n \}}
\). But then \( \Z^{\g} \o \cvx \Z^\beta\), and since $\Z^\beta \cvx \Z^{\g} \cong \Z^\gamma \times \{ 0 \} $ by \( \beta \leq \gamma \) (see Definition~\ref{def:Z^alpha}) this shows that $\Z^\beta$ is right compressible, against the induction hypothesis.
\end{proof}

Using Proposition \ref{power_Z} again, we have 
\begin{align*}
\bigg(\sum_{\gamma < \alpha} \Z^\gamma \o\bigg)^* + \bigg(\sum_{\gamma < \beta} \Z^\gamma \o \bigg) & = \bigg(\sum_{\beta \leq \gamma < \alpha} \Z^\gamma \o\bigg)^* + \bigg(\sum_{\gamma < \beta} \Z^\gamma \o\bigg)^* + \bigg( \sum_{\gamma < \beta} \Z^\gamma \o \bigg) \\
& \cong \bigg(\sum_{\beta \leq \gamma < \alpha} \Z^\gamma \o\bigg)^* + \Z^\beta
\end{align*}
Let \(g\) be the isomorphism between the first and last element of
this chain. Choose \(z \in \Z^\a\) such that
\(g(f(z)) \in \Z^\beta\) --- such a \(z\) exists because \(f(\Z^\a)\) is
cofinal in
\(\big(\sum_{\gamma < \alpha} \Z^\gamma \o\big)^* + \big( \sum_{\gamma <
\beta} \Z^\gamma \o \big)\) by Claim~\ref{claim:Z^alpha}. Arguing as before, 
\( \Z^{\g} \omega \cvx [g(f(z)), +\infty)_{ \Z^{\beta}} \) for some \(\beta \leq \gamma < \alpha \),
contradicting again the incompressibility of $\Z^\beta$.

Finally,
assume that \( f(\Z^\alpha) \cap \big(\sum_{\g < \a} \Z^{\g} \o \big) = \emptyset \), i.e.\ \( f(\Z^\alpha) \subseteq \big( \sum_{\g < \a} \Z^{\g}\o \big)^* \). Let \( \beta < \a \) be smallest such that \( f(\Z^\alpha) \cap (\Z^\beta \omega)^* \neq \emptyset \), and let \( n  \) be smallest such that \( f(\Z^\alpha) \cap (\Z^\beta \times \{ n \})^* \neq \emptyset \).
Pick \( z \in \Z^\alpha \) such that \( f(z) \in (\Z^\beta \times \{ n \})^* \).
Arguing as before, there is \( \beta \leq \gamma < \alpha \) such that \( \Z^\gamma \omega \cvx (\Z^\beta \times \{ n \})^* \).
Since \( (\Z^\beta \times \{ n \})^* \cong (\Z^\beta)^* \cong \Z^\beta \) by Proposition~\ref{power_Z^L}, this would mean that \( \Z^\gamma \omega \cvx \Z^\beta \), contradicting again the incompressibility of the latter.
\end{proof}

\begin{thm}\label{compl_subsets}
\begin{enumerate-(a)}
\item \label{compl_subsets-a}
\(\LO_{l}\), \(\LO_{r}\) and $\LO_{r} \cap \LO_l$ are
$\mathbf{\Sigma}_1^1$-complete.

\item \label{compl_subsets-b}
$\LO \setminus (\LO_r \cup \LO_l)$ is
$\mathbf{\Pi}_1^1$-complete.

\item \label{compl_subsets-c}
$\LO_{r}\setminus \LO_l$ and $\LO_{l}\setminus \LO_r$
are $D_2(\mathbf{\Pi}_1^1)$-complete.
\end{enumerate-(a)}
\end{thm}

\begin{proof}
\ref{compl_subsets-a} 
First, we check that $\LO_{r}$ is
$\mathbf{\Sigma}_1^1$. Indeed, \( L \in \LO_{r} \) if and only
if
\begin{equation*}
\begin{aligned} 
\exists f\colon \N \to \N \big[ & \forall n,m \; (n< _L m \rightarrow f(n)<_L f(m)) \; \wedge\\
& \forall n,m,k \; (f(n)\leq_L k \leq_L f(m)\rightarrow \exists k'(f(k')=k)) \; \wedge\\
& \exists n\; \forall m \; (f(m) <_L n)\big].
\end{aligned}  
\end{equation*}
In a similar way, one can prove that $\LO_{l}$ (and hence also
$\LO_l \cap \LO_r$) is $\mathbf{\Sigma}_1^1$.

We now show that $\LO_{l}$, \(\LO_{r}\) and $\LO_{r} \cap \LO_l$ are
$\mathbf{\Sigma}_1^1$-hard by continuously reducing the \( \mathbf{\Sigma}^1_1 \)-complete set \( \LO \setminus \WO \) to each of them. We can actually use the continuous
function $L \mapsto \Z^L$ for all three sets.
Indeed, if $L \notin \WO$, by Proposition \ref{power_Z^L} we have
\(\Z^L \cong \mathbb{Z}^\alpha \eta\) for some ordinal \(\alpha\), and hence $\Z^L$ is obviously bicompressible.  If $L \in \WO$, then
$\Z^L$ is incompressible by Lemma \ref{incomp_z}.

\ref{compl_subsets-b} is immediate from the proof of \ref{compl_subsets-a}.

\ref{compl_subsets-c} By \ref{compl_subsets-a} it follows that \(\LO_r \setminus \LO_l\) and
\(\LO_l \setminus \LO_r\) are $D_2(\mathbf{\Pi}_1^1)$. Consider
now the set
$A = \{(L,L') \in \LO \times \LO \mid L \notin \WO \mbox{ and } L' \in
\WO \}$ and recall that it is
$D_2(\mathbf{\Pi}_1^1)$-complete. Define the continuous map
\(\psi \colon \LO \times \LO \to \LO\) by
\(\psi(L, L') = \Z^{1+L'} + \eta + \Z^{1+L}\).

We claim that \( \psi(L,L') \) is left compressible if and only if \( L' \notin \WO \). One direction is obvious: if \( L' \notin \WO \), then \( \Z^{1+L'} \cong \Z^\alpha \eta \) for some \( \alpha \geq 1 \), and thus it has a convex self-embedding onto a proper final segment of it, which can then be naturally extended to a witness of \( \psi(L,L') \in \LO_l \). For the other direction, we use the fact that every convex subset of \( \eta  \) consists of points which have neither an immediate predecessor nor an immediate successors, while convex subsets of \( \Z^{1+L} \) and \( \Z^{1+L'} \) with at least two points always contain elements with both an immediate predecessor and an immediate successor in the given linear order. (Here we use again the fact that \( \Z^{1+L} \) and \( \Z^{1+L'} \) are either of the form \( \Z^\alpha \) or \( \Z^\alpha \eta \) for some $\a\geq 1$, depending on whether \( L \) and \( L' \) are well-orders or not.) Thus if \( f \colon \psi(L,L') \to \psi(L,L') \) is a convex embedding we must have \( f(\eta) = \eta \), and hence \(f(\Z^{1+L'}) \subseteq \Z^{1+L'} \). Thus if \( L' \in \WO \) then \( \Z^{1+L'} \notin \LO_l \) by Lemma \ref{incomp_z}, which implies  \( f(\Z^{1+L'}) = \Z^{1+L'} \): since \( f \) was arbitrary, this shows that \( \psi(L,L') \notin \LO_l \). 

Analogously, one can check that \( \psi(L,L') \) is right compressible if and only if \( L \notin \WO \). Using these facts, 
it is then easy to prove that
\((L,L') \in A\) if and only if
\(\psi(L,L') \in \LO_{r}\setminus \LO_l\), hence \( \psi \) witnesses that
$\LO_{r}\setminus \LO_l$ is $D_2(\mathbf{\Pi}_1^1)$-hard. 

For
$\LO_{l}\setminus \LO_r$ it suffices to switch the positions of
$\Z^{1+L}$ and $\Z^{1+L'}$ in the definition of $\psi$.
\end{proof}

Even if they are not Borel, the sets \(\LO \setminus (\LO_{r}\cup \LO_{l})\), \(\LO_{l}\setminus \LO_{r}\), \(\LO_{r}\setminus \LO_{l}\) and \(\LO_{r}\cap \LO_{l}\) belong to the Boolean algebra generated by the analytic sets, and hence have the Baire property and are universally measurable. 
By Theorem \ref{red_cvxeq_jump_iso} and Proposition \ref{Baire_red_for_qclas_count_str} we thus obtain the following result.

\begin{cor}\label{cvxeq_baire_iso}
The equivalence relation
\(\cvxeq\) is \(\sigma\)-classifiable by countable structures, and
therefore
\(\cvxeq \leq_{\text{\scriptsize
\textit{Baire}}} {\iso}\).
\end{cor}

Notice that, since the partition of $\LO$ given by \eqref{subsets of LO} is finite, we actually have that the preimages of Borel sets via the reduction of $\cvxeq$ to $\iso$ are Boolean combinations of analytic sets. 
The problem of whether \(\cvxeq\) is Borel reducible to \(\iso\) remains open. 
However, from the reductions above we can derive some more information about the complexity of \(\cvxeq\), showing that it shares some important properties with $\iso$.

\begin{cor}\label{turb}
If \(X\) is a Polish space on which the action of a Polish group \(G\) is turbulent, then \( E^X_G \nleq_B \cvxeq \).
\end{cor}

\begin{proof}
If \(E^X_G \leq_B \cvxeq\), then by Corollary \ref{cvxeq_baire_iso}
we would have that
\(E^X_G \leq_{\text{\scriptsize \textit{Baire}}}
{\iso}\), against Theorem \ref{turb_act}.
\end{proof}

In Proposition~\ref{cor 2} we observed that \( (\iso)^+ \leq_{B} {\iso} \). Replacing Borel reducibility with Baire reducibility, we get an analogous result for \( \cvxeq \).

\begin{cor}
\({(\cvxeq)^+} \leq_{\text{\scriptsize \textit{Baire}}} {\cvxeq}\).
\end{cor}

\begin{proof}
Since \({\cvxeq} \leq_{\text{\scriptsize \textit{Baire}}} {\iso}\), we
have that
\({(\cvxeq)^+} \leq_{\text{\scriptsize \textit{Baire}}} {(\iso)^+}\),
but \({(\iso)^+} \leq_B {\iso} \leq_B {\cvxeq}\), so
\({(\cvxeq)^+} \leq_{\text{\scriptsize \textit{Baire}}} {\cvxeq}\).
\end{proof}



\begin{cor}\label{E1_cvxeq}
\({E_1} \nleq_{\text{\scriptsize \textit{Baire}}} {\cvxeq}\).
\end{cor}

\begin{proof}
If \(E_1 \leq_{\text{\scriptsize \textit{Baire}}} {\cvxeq}\), by
Corollary \ref{cvxeq_baire_iso} we would have
$E_1 \leq_{\text{\scriptsize \textit{Baire}}} {\iso}$, contradicting Theorem \ref{E1_orbit}.
\end{proof}


Each one of Corollaries \ref{turb} and \ref{E1_cvxeq} implies that $\cvxeq$ is not
complete for analytic equivalence relations, thus by Proposition
\ref{compl_qo} we obtain:

\begin{cor} \label{cor:notcompletepapillon} $\trianglelefteq_{\LO}$ is
not complete for analytic quasi-orders.
\end{cor}

Recall that by \(\Cvx(\R)\) we denote the set of the open intervals of
\(\R\). We can naturally equip \(\Cvx(\R)\) with a Polish topology: indeed, if we extend the usual order on $\R$ to $\R \cup \{\pm \infty\}$ in the obvious way, then $\Cvx(\R)$ is the open subset $\{(x,y) \mid x < y\}$ of the Polish space $(\R\cup \{\pm \infty\})^2$.  The inclusion relation on $\Cvx(\R)$ is then closed. Notice now that the embedding
from \((\Cvx(\R), \subseteq)\) to \((\LO,\cvx)\) defined in the proof
of Lemma \ref{open_int} is actually a Borel reduction.
Thus we have
the following corollary.

\begin{cor} \label{cor:intervalsleqBcvx}
\((\Cvx(\R),\subseteq) \leq_B (\LO,\cvx)\).
\end{cor}

\subsection{Convex embeddability between countable circular orders}\label{sec:cvx_co}

Our goal in this section is to define a relation of convex
embeddability among circular orders. We first recall the definition of
convex subset of a circular order as given by Kulpeshov and Macpherson
(\cite{KM05}).

\begin{defn} \label{def:en}
Let \(C\) be a circular order. The set \(A \subseteq C\) is said to
be \textbf{convex} in \(C\), in symbols \(A \csube C\), if for any distinct
\(x,y \in A\) one of the following holds:
\begin{enumerate-(i)}
\item\label{en1} for every \(c \in C\) with \(C(x,c,y)\) we have
\(c \in A\);
\item\label{en2} for every \(c \in C\) with \(C(y,c,x)\) we have
\(c \in A\).
\end{enumerate-(i)}
If $A$ is a proper subset of $C$ we write $A \csub C$.
\end{defn}

Note that if \(A \csub C\) then exactly one of \ref{en1} and
\ref{en2} holds for each pair of distinct \(x,y \in A\).

The following propositions collect some basic properties of convex subsets of circular orders.

\begin{prop}\label{pro:complconvex}
If \(C\) is a circular order and \(A \csube C\) then \(C \setminus A\) is a convex subset of $A$ as well.  
\end{prop}

\begin{proof}
If \( C \setminus A \) is empty or a singleton the result is trivial, so we can assume that \( C \setminus A \) contains at least two points.
Toward a contradiction, suppose \(x, y \in C \setminus A\) are distinct and such that:
\begin{enumerate-(1)}
\item there exists \(c_0 \in A\) with \(C(x,c_0,y)\), and
\item there exists \(c_1 \in A\) with \(C(y,c_1,x)\).
\end{enumerate-(1)}
By cyclicity and transitivity we obtain \(C(c_0,y,c_1)\) and
\(C(c_1,x,c_0)\), and since \(A\) is convex we would have that at
least one of \(x\) and \(y\) belongs to \(A\), a contradiction.
\end{proof}

The previous proposition highlights a major difference between convex subsets of circular and linear orders: the complement of a convex subset of a linear order is not in general convex. On the other hand, convex subsets of linear orders are closed under intersections, while this is not the case for circular orders: consider the circular order $C[\boldsymbol{4}]$ and its convex subsets $\{0,1,2\}$ and $\{2,3,0\}$. However the intersection of two convex subsets of a circular order is not convex only in some circumstances.

\begin{prop}\label{prop:intconvex}
Let \(C\) be a circular order. If \(A,B\csube C\)
then \(A \cap B\) is the union of two convex subsets of \(C\). Moreover, if $A \cap B$ is not convex then $A \cup B = C$.
\end{prop}

\begin{proof}
If \(B \subseteq A\) or \(A \subseteq B\), the result is trivial. So, suppose there exist \(w \in A \setminus B\) and \(z \in B\setminus A\), and consider the partition of \(A \cap B\) given by the sets
\[
A_1 = \{x \in A \cap B \mid C(w,x,z)\} \quad \text{and} \quad
A_2 = \{x \in A \cap B \mid C(z,x,w)\}.
\]

We claim that \(A_1 \csube C\). Let \(x,y\in A_1\) be distinct: since $x,y \in A$, without loss of generality we can assume that $u \in A$ for every $u \in C$ such that $C(x,u,y)$. Since $z \notin A$ we have that $C(x,z,y)$ fails and, by totality and cyclicity, we have \(C(x,y,z)\). Using cyclicity, transitivity and $C(w,x,z)$ we obtain $C(y,w,x)$. Since $w \notin B$ and $B$ is convex this implies that $u \in B$ for every $u \in C$ such that $C(x,u,y)$. If now $u$ is such that $C(x,u,y)$ we already showed that $u \in A \cap B$. From $C(y,w,x)$ and $C(x,u,y)$ it follows that $C(y,w,u)$ which, combined with $C(w,y,z)$, yields $C(w,u,z)$ and hence $u \in A_1$.
The proof that \(A_2\) is convex is symmetric.

Now assume that $A \cap B$ is not convex, and hence both $A_1$ and $A_2$ are non empty. Fix $x \in A_1$ and $y \in A_2$. From $C(w,x,z)$ and $C(z,y,w)$ it follows that we have $C(x,z,y)$ and $C(y,w,x)$. Since $x,y \in A$ but $z \notin A$ we must have that $C(y,u,x)$ implies $u \in A$. Similarly we obtain that $C(x,u,y)$ implies $u \in B$. By totality it follows that $A \cup B = C$.
\end{proof}

\begin{prop}\label{prop:partconvex}
Let \(C\) be a circular order. Let $\{A_i \mid i \in I\}$ and $\{B_j \mid j \in J\}$ be two collections of pairwise disjoint convex subsets of $C$. Then there exists at most one pair $(i,j) \in I \times J$ such that $A_i \cap B_j$ is not convex.
\end{prop}

\begin{proof}
Suppose that $A_i \cap B_j$ is not convex. By the second part of Proposition \ref{prop:intconvex} we have $A_i \cup B_j = C$. Hence for every  $i' \neq i$ and $j' \neq j$ we have that $A_{i'} \subseteq B_j$ and $B_{j'} \subseteq A_i$. Therefore $A_{i'} \cap B_j =A_{i'}$, $A_i \cap B_{j'} =B_{j'}$ and $A_{i'} \cap B_{j'} \subseteq A_{i'} \cap A_i = \emptyset$ are all convex.
\end{proof}

If \( f \) is an embedding between linear orders \( L \) and \( L' \) and \( f(L) \csube L' \), then \( f(B) \csube L' \) for every \( B \csube L \). This ceases to be true for circular orders, as shown by the following example. The identity map between \( C = C[\zeta] \) and \( C' = C[\zeta + \boldsymbol{1}] \) has convex range, but the image of the convex set \( B = C \setminus \{ 0 \} \csube C \) is no longer convex in \( C' \). The following proposition gives a weakening of the above property which is however sufficient for the ensuing proofs.

\begin{prop} \label{prop:embeddingandconvex}
Let \( f \) be an embedding between the circular orders \( C \) into \( C' \). If \( A' \csube C' \), then \( f^{-1}(A') \csube C \). Conversely, if \( A \csub C \) is such that \( f(A) \csube C' \), then \( f(B) \csube C' \) for all \( B \csube C \) with \( B \subseteq A \). 
\end{prop}

\begin{proof}
%
The first part is obvious, so let us consider \( A \csub C \) with \( f(A) \csube C' \), and fix any \( B \csube C \) contained in \( A \).
Pick distinct points \( f(x),f(y) \in f(B) \subseteq f(A) \), so that \( x,y \in B \) and \( x \neq y \) because \( f \) is injective. 
Since \( A \csub C \), without loss of generality we might assume that \( c \in A \) for all \( c \in C \) with \( C(x,c,y) \) and
that there is \( d \in C \) with \( C(y,d,x) \) and \( d \notin A \), so that the same is true with \( A\) replaced by \( B \) because \( B \subseteq A \) is convex. Since \( f \) is an embedding, \( f(d) \) is such that \( C'(f(y),f(d),f(x)) \) but \( f(d) \notin f(A) \). Since \( f(A) \csube C' \) by hypothesis, this means that \( c' \in f(A) \) for all \( c' \in C \) such that \( C'(f(x),c',f(y)) \). 
So for such a \( c' \in C' \) there is \( c \in A \) such that \( c' = f(c)\): then \( C(x,c,y) \) because \( f \) is an embedding, and so \( c \in B \) and \( f(c)= c' \in f(B) \). This shows that \( f(B) \) satisfies~\ref{en1} of Definition~\ref{def:en} with respect to $x$ and $y$. Hence \( f(B) \csube C' \).
\end{proof}

The first natural attempt to define convex embeddability between circular orders is the following.

\begin{defn}
Let \(C\) and \(C'\) be circular orders. We say that \(C\) is \textbf{convex embeddable} into \(C'\), and write \(C \ccvx C'\), if there exists an embedding \( f\) from \( C \) to \( C' \) such that \( f(C) \csube C'\).
\end{defn} 

However, \(\ccvx\) is \textbf{not} transitive, as witnessed by
\(C[\z] \ccvx C[\z+\boldsymbol{1}]\),
\(C[\z+\boldsymbol{1}] \ccvx C[\o+\boldsymbol{1}+\o^*+\eta]\) (because
\(C[\z+\boldsymbol{1}] \cong_{c} C[\o+\boldsymbol{1}+\o^*]\)), and
\(C[\z] \nccvx C[\o+\boldsymbol{1}+\o^*+\eta]\). Nevertheless, notice that if
we partition \(C[\z]\) into the \emph{two} convex subsets \(\o^*\) and
\(\o\) then they are isomorphic to the two convex subsets \(\o^*\) and
\(\o\) of~\(C[\o+\boldsymbol{1}+\o^*+\eta]\).

By taking the transitive closure of $\ccvx$ (i.e.\ the smallest binary relation containing $\ccvx$) we are naturally led to the following definition. We call \textbf{finite convex partition} of the circular order \( C \) any finite partition \( \{ C_i \mid i < n \} \) of \( C \) such that 
\begin{itemizenew}
\item
\( C_i \csube C \) for all \( i < n \), and 
\item
for all \( x,y,z \in C \), if \( C(x,y,z) \) then \( C[\boldsymbol{n}](i,j,k) \) for the unique \( i,j,k < n \) such that \( x \in C_i \), \( y \in C_j \), and \( z \in C_k \).
\end{itemizenew}
Notice that this implies that the \( C_i \)'s are ordered as \( C[\boldsymbol{n}] \), that is: if \( i,j,k < n \) are distinct and \( C[\boldsymbol{n}](i,j,k) \) then \( C(x,y,z) \) for every \( x \in C_i \), \( y \in C_j \), and \( z \in C_k \). Also, the convexity of the \( C_i \)'s follows from the second condition if \( n \geq 3 \).

\begin{defn}\label{cvx_co}
Let \(C\) and \(C'\) be circular orders. We say that \(C\) is \textbf{piecewise convex embeddable} into \(C'\), and write \(C \pcvx C'\), if there are a finite convex partition \( \{ C_i \mid i < n \} \) of \( C \) and an embedding \( f \) of \( C \) into \( C'\) such that \( f(C_i) \csube C' \) for all \( i < n \). 

We denote by \(\pccvx\) the restriction of \(\pcvx\) to the set \(\CO\) of (codes for) circular orders on \(\N\).
\end{defn}

Clearly, $C \ccvx C'$ implies $C \pcvx C'$.
Notice also that when \( C \) has at least two elements and \( C \pccvx C' \) as witnessed by \( \{ C_i \mid i < n \} \) and \( f \), without loss of generality we can assume that \( n > 1 \) and hence \(C_i \csub C \). (If not, split \( f(C_0) \csube C' \) into two nonempty convex subsets \( A, B \) of \( C' \), and consider the finite convex partition \( \{ f^{-1}(A), f^{-1}(B) \} \) of \( C \) together with the same embedding \( f \).)

\begin{prop}\label{prop:pccvx_is_trans}
$\pcvx$ is transitive.
\end{prop}

\begin{proof}
Suppose that $C \pcvx C'$, as witnessed by the embedding \( f \) and the finite convex partition \( \{ C_i \mid i<n \}\) of \(C\),  
and that $C' \pcvx C''$ via the embedding \( g \) and the finite convex partition \( \{ C'_j \mid j < m \} \) of \(C'\). 
If $C'$ has only one element than so does $C$ and \( C \pcvx C'' \) is immediate. Thus, without loss of generality, we can assume that \(  m  > 1 \),
so that \( C'_j \csub C' \) for all \( j < m \).
Notice that $\{f(C_i) \mid i < n\}$ and $\{C'_j \mid j<m\}$ are two collections of pairwise disjoint convex subsets of $C'$. We distinguish two cases. 

If
\(
C'_{i,j} = f(C_i) \cap C'_j
\)
is a convex subset of \( C' \) for every \( i < n \) and \( j < m \), 
then we can order the family of pairwise disjoint convex sets 
\[
\left\{ C'_{i,j} \mid (i,j) \in n \times m \wedge C'_{i,j} \neq \emptyset \right\} 
\] 
following the circular order of \( C'\). In this way we obtain a family
\( \{ D'_k \mid k < \ell \} \), for the suitable \( \ell \leq n \cdot m\), such that if \( x_0 \in D'_{k_0} \), \( x_1 \in D'_{k_1} \), and \( x_2 \in D'_{k_2} \) satisfy \( C'(x_0,x_1,x_2) \) then
\( C[\boldsymbol{\ell}](k_0,k_1,k_2)\).
Then \( \{ f^{-1}(D'_k) \mid k < \ell \} \) is a finite convex partition of \( C \) and \( g \circ f\) is an embedding of \( C \) into \( C'' \). Moreover, for every \( k < \ell \) we have \( (g \circ f)(f^{-1}(D'_k)) = g(D'_k) \csube C''\) because \( D'_k \subseteq C'_j \csub C' \) for some \( j< m\) (Proposition~\ref{prop:embeddingandconvex}).
Thus \( C \pcvx C'' \).

Suppose now that \( C'_{i,j} = f(C_i) \cap C'_j \) is not convex for some \( (i,j) \).   
By Proposition \ref{prop:partconvex} there is at most one such pair \((\bar{\imath},\bar{\jmath})\). By Proposition \ref{prop:intconvex}, \(C'_{\bar{\imath}, \bar{\jmath}}\) is the union of two disjoint convex subsets \(A_0\) and \(A_1\) of \(C'\). Then we can argue as in the previous paragraph but starting with the family 
\[
\left\{ C'_{i,j} \mid (i,j) \in n \times m \wedge (i,j) \neq (\bar{\imath},\bar{\jmath}) \wedge C'_{i,j} \neq \emptyset \right\} \cup \{ A_0,A_1 \}. \qedhere
\]
\end{proof}

Thus \( \pcvx \) is a quasi-order, and  it is easy to see that its restriction \( \pccvx \) to the Polish space  \( \CO \) is analytic.

We first show that \(\pccvx \)
satisfies combinatorial properties similar to those proved
for \(\cvx\) in Section \ref{sec:comb_prop}.
A key point is that it still makes sense to talk about (finite) condensation in the realm of circular orders. Indeed, given a circular order \( C \) the \textbf{condensation class} \( c_F^C(\ell) \) of \( \ell \) is the collection of those \( m \) such that either \( \{ k \mid C(\ell,k,m) \} \) or \( \{ k \mid C(m,k,\ell) \} \) is finite. Each \( c_F^C(\ell) \) is convex in \( C \), and it again holds that the condensation classes form a partition of \( C \). This allows us to define the (\textbf{finite}) \textbf{condensation} \( C_F \) of \( C \) in the obvious way.
The crucial observation is that we can substitute Proposition~\ref{cond_classes} with the following lemma.

\begin{lem} \label{lem:condensationCO}
Let \( f \) be an embedding between the circular orders \( C \) and \( C' \). Fix any \( \ell \in C \), and let \( A \csub C \) and \( a,b \in A \setminus c_F^C(\ell) \) be such that \( f(A) \csube C' \), \( c_F^C(\ell) \subseteq A \) and \( C(a,\ell',b) \) for all \( \ell' \in c_F^C(\ell) \).
Then the restriction of \( f \) to \( c_F^C(\ell) \) is an isomorphism between \( c_F^C(\ell) \) and \( c_F^{C'}(f(\ell)) \), and thus 
\( |c_F^C(\ell)|=|c_F^{C'}(f(\ell))| \).
\end{lem}

\begin{proof}
By Proposition~\ref{prop:embeddingandconvex}, we have \( f(c_F^C(\ell)) \csube C' \), which easily implies \( f(c^C_F(\ell)) \subseteq c_F^{C'}(f(\ell)) \).
Conversely, pick any \( d' \in c_F^{C'}(f(\ell)) \) distinct from \( f(\ell) \), and first assume that \( \{ k \mid C'(d',k,f(\ell)) \} \) is finite. 
Consider the set \( B = \{ k \in C \mid C(a,k,\ell)\} \subseteq A\). 
Since \( a \notin c_F^C(\ell) \) and \( B \csube C \), the set \( f(B) \) is infinite and by Proposition \ref{prop:embeddingandconvex} \( f(B) \csube C' \). 
We cannot have \( C' (d',f(a),f(\ell)) \), otherwise \( \{ k \mid C'(d',k,f(\ell)) \} \supseteq f(B) \) and the former would be infinite. 
Thus \( C'(f(a),d',f(\ell)) \), and so \( \{ k \mid C'(d',k,f(\ell)) \} \subseteq f(B) \) because \( f(B) \csube C' \).
This easily implies that \( d' = f(d) \) for some \( d \in c_F^C(\ell) \), and we are done.
(When \( \{ k \mid C'(f(\ell),k,d') \} \) is finite, we work symmetrically on the other side of \( \ell \) and use \( b \) instead of \( a \).)
\end{proof}

\begin{prop} \label{prop:antichainsCO}
\begin{enumerate-(a)}
\item \label{prop:antichainsCO-a}
There is an embedding from the partial order
\((\Cvx(\R),\subseteq)\) into $\pccvx$, and indeed \( {(\Cvx(\R),\subseteq)} \leq_B {\pccvx} \).
\item \label{prop:antichainsCO-b}
$\pccvx$ has chains of order type \( (\R, {<} ) \), as well as antichains of size $2^{\aleph_0}$.
\end{enumerate-(a)}
\end{prop}

\begin{proof}
Given an interval \( (x,y) \in \Cvx(\R)\), consider the circular order \( C_{(x,y)} \in \CO \) defined by \( C_{(x,y)} =  C[L_{(x,y)}]  \), where \( L_{(x,y)} \) is as in the proof of Lemma~\ref{open_int}. We claim that the map \( (x,y) \mapsto C_{(x,y)}\) witnesses~\ref{prop:antichainsCO-a}.
Pick two intervals \( (x,y) , (x',y') \in \Cvx(\R) \). If \( (x,y) \subseteq (x',y') \) then the identity map witnesses \( C_{(x,y)} \ccvx C_{(x',y')}\), hence \( C_{(x,y)} \pccvx C_{(x',y')}\). Suppose now that \( (x,y) \not\subseteq (x',y') \), and for the sake of definiteness assume that \( x < x' \).
Towards a contradiction suppose that there are a finite convex partition \( \{ C_i \mid i < n \} \) of \( C_{(x,y)} \) and an embedding \( f \) witnessing \( C_{(x,y)} \pccvx C_{(x',y')}\).
As usual, we can assume \( n > 1\), so that \( C_i \csub C_{(x,y)}\) for all \( i < n \). Since there are infinitely many rationals between \( x \) and \( x' \) and all condensation classes of \( C_{(x,y)} \) are finite, we can find \( i < n \) and \( q,q',r \in \Q \) such that \( x < q < r < q' < x' \) and the hypothesis of Lemma~\ref{lem:condensationCO} are satisfied with \( A = C_i \), \( a = (0,q) \), \( b = (0,q') \) and \( \ell = (0,r) \). Thus the condensation class of \( f(0,r) \) has the same size of the condensation class of \( (0,r) \), which by construction can happen only if \( r \in (x',y')\), a contradiction.

Part~\ref{prop:antichainsCO-b} is derived from~\ref{prop:antichainsCO-a} as in Proposition~\ref{prop1}.
\end{proof}

\begin{prop} \label{prop:boundingCO}
\( \mathfrak{b}(\pccvx) = \aleph_1\), and indeed every \( C \in \CO \) is the bottom of a strictly increasing \( \pccvx \)-unbounded chain of length \( \omega_1\).
\end{prop}

\begin{proof}
For \( C \in \CO \) and \( \ell \in C \), let \( \alpha_{\ell,C} \) be the sup of those \( \omega \leq \alpha < \omega_1\) such that \( C[\boldsymbol{\alpha}] \ccvx C \) via some \( f \) satisfying \( f(0) = \ell \). Since \( \alpha_{\ell,C} \) is attained by definition of convexity, the ordinal \( \alpha_C = \sup_{\ell \in C} \alpha_{\ell,C} \) is countable, and by construction \( C[\boldsymbol{\alpha}_C+\boldsymbol{1}] \not\trianglelefteq_c C \). Let \( \alpha \) be an additively indecomposable%
\footnote{An ordinal \( \alpha\) is additively indecomposable if \( \beta+ \gamma < \alpha \) for all \( \beta, \gamma < \alpha \). Additively indecomposable ordinals are precisely those of the form \( \omega^\delta \) for some ordinal \( \delta \).}
countable ordinal above \( \alpha_C+1\): we claim that \( C[\boldsymbol{\alpha}] \npccvx C\). 
Suppose towards a contradiction that \( \{ C_i \mid i < n \} \) is a finite convex partition and \( f \) an embedding witnessing \( C[\boldsymbol{\alpha}] \pccvx C \). As usual we can assume \( n> 1\). 
Then there are \( i < n \) and \( \gamma < \alpha \) such that \( A'_\gamma =  \{ \beta \in C[\boldsymbol{\alpha}] \mid \beta \geq \gamma \} \) is contained in \( C_i \). 
Since \( \alpha \) is additively indecomposable, the linear order determined by \( A_\gamma \) has order type \( \alpha \geq \alpha_C+1 \), thus we can consider the set \( A_\gamma = \{ \beta \in A'_\gamma \mid \beta < \gamma+\alpha_C+1 \} \), which has order type \( \alpha_C+1 \). Since \( A_\gamma \csube C[\boldsymbol{\alpha}] \) and \( A_\gamma \subseteq C_i \csub C[\boldsymbol{\alpha}] \), by Proposition~\ref{prop:embeddingandconvex} the restriction of \( f \) to \( A_\gamma \) witnesses \( C[\boldsymbol{\alpha}_C+\boldsymbol{1}] \ccvx C\), a contradiction.

This shows that the family \( \{ C[\boldsymbol{\alpha}] \mid \omega \leq \alpha < \omega_1 \} \) is \( \pccvx \)-unbounded in \( \CO \). Since \( C[\boldsymbol{\alpha}] \ccvx C[\boldsymbol{\beta}] \) when \( \alpha \leq \beta \), we can extract from it a strictly increasing chain witnessing \( \mathfrak{b}(\pccvx) \leq \aleph_1 \). 
To show that \( \mathfrak{b}(\pccvx) > \aleph_0\), consider a countable family \( \{ C_i \mid i \in \N \} \subseteq \CO \). For each \( i \in \N \) pick an arbitrary \( \ell_i \in C_i \) and define \( L_i \in \LO\) by setting \( x \leq_{L_i} y \) iff \( C_i(\ell_i,x,y) \). Then the circular order \( C = C \left[  \sum_{i \in \boldsymbol{\N}} L_i \right] \in \CO \) is such that \( C_i \ccvx C \) for all \( i \in \N \), and thus the given family is \( \pccvx \)-bounded.

For the second part, pick \( \ell \in C \) and let \( L \in \LO\) be defined by \( x \leq_L y \) iff \( C(\ell,x,y) \). Consider the \( \pccvx \)-nondecreasing sequence \( (C_\alpha)_{\alpha < \omega_1}\) of circular orders defined by \( C_0 = C  \) and \( C_\alpha = C[L+\boldsymbol{\omega}+\boldsymbol{\alpha}] \) when \( \alpha > 0 \). Since \( C[\boldsymbol{\omega}+\boldsymbol{\alpha}] \ccvx C_\alpha \) for all \( \alpha \neq 0 \), such a sequence is \( \pccvx \)-unbounded. Thus we can extract from it a strictly \( \pccvx \)-increasing subsequence of length \( \omega_1 \) with \( C_0 \) as first element: being cofinal in the original sequence, it will be \( \pccvx \)-unbounded too, as required.
\end{proof}

\begin{prop} \label{prop:basisCO}
\begin{enumerate-(a)}
\item \label{prop:basisCO-1}
There are \( 2^{\aleph_0}\)-many \( \pccvx \)-incomparable \( \pccvx \)-minimal elements in \( \CO \). In particular, all bases for \( \pccvx \) are of maximal size.

\item \label{prop:basisCO-2}
There is a \( \pccvx\)-decreasing \( \omega \)-sequence in \( \CO \) which is not \( \pccvx \)-bounded from below.
\end{enumerate-(a)}
\end{prop}

\begin{proof}
\ref{prop:basisCO-1}
Given an infinite \( S \subseteq \N \), let \( C_S = C[\eta_{f_S}]\) where \( \eta_{f_S} \) is as in the proof of Proposition~\ref{prop:basisforcvx}\ref{prop:basisforcvx-1}.
If \( A \csube C_S \) is infinite, then there exist \( q,q' \in \Q \) with \( q < q' \) such that \( \{ (\ell,q'') \in C_S \mid q \leq q'' \leq q' \} \subseteq A \),
and thus \( C_S \) is convex embeddable into (the circular order determined by) \( A \).

Let \( C \pccvx C_S \), as witnessed by the finite convex partition \( \{ C_i \mid i < n \} \) and the embedding \( f \).
Then there is \( i < n \) such that \( C_i \), and hence also \( f(C_i) \) is infinite. Setting \( A = f(C_i) \) in the previous paragraph, we get that \( C_S \ccvx f(C_i) \cong C_i \csube C \), hence \( C_S \pccvx C \). This shows that \( C_S \) is \( \pccvx \)-minimal.

Assume now that \( C_S \pccvx C_{S'} \) for some infinite \( S,S' \subseteq \N \), and let \( \{ C_i \mid i < n \} \) be a finite convex partition of \( C_S \) and \( f \colon C_S \to C_{S'} \) be an embedding witnessing this. As usual, we may assume \( n > 1 \), so that \( C_i \csub C_S \) and \( f(C_i) \csub C_{S'} \). Fix any \( i < n \) such that \( C_i \) is infinite. By the first paragraph, there are \( q < q' \) such that \( \{ (\ell,q'') \in C_S \mid q \leq q'' \leq q' \} \subseteq C_i \). Given an arbitrary \( m \in S \), pick \( q'' \in \Q \) such that \( q < q'' < q' \) and \( f_S(q'') = m \). 
Then the hypotheses of Lemma~\ref{lem:condensationCO} are satisfied when we set \( A = C_i \), \( a = (0,q) \), \( b = (0,q') \), and \( \ell = (0,q'') \). Thus \( C_{S'} \) must contain a condensation class of size \( m \), which is possible only if \( m \in S'\). This shows that \( S \subseteq S' \). Conversely, given \( m \in S' \) we work with the infinite set \( f(C_i) \csub C_{S'} \) and pick \( q,q' \in \Q \) such that \( \{ (\ell,q'') \in C_{S'} \mid q \leq q'' \leq q' \} \subseteq f(C_i) \). Then we pick \( q'' \in \Q \) such that \( q < q'' < q' \) and \( f_{S'}(q'') = m \).
Applying (the proof of) Lemma~\ref{lem:condensationCO} we get that the condensation class of \( f^{-1}(0,q'')\) has size \( m \), hence \( m \in S\). Since \( m \in S' \) was arbitrary, \( S' \subseteq S \), and thus \( S = S' \). This shows that \( \{ C_S \mid S \subseteq \N \wedge S \text{ is infinite} \}\) is a \( \pccvx \)-antichain and we are done.

\ref{prop:basisCO-2}
Consider the family \( \{ C_{(m,+\infty)} \mid m \in \N \} \), where \( C_{(m,+\infty)} \) is as in the proof of Proposition~\ref{prop:antichainsCO}\ref{prop:antichainsCO-a}.
It is a strictly \( \pccvx \)-decreasing chain, so we only need to show that it is \( \pccvx \)-unbounded from below. Let \( C \in \CO \) be such that \( C \pccvx C_{(0,+\infty)} \), as witnessed by the finite convex partition \( \{ C_i \mid i < n \} \) (for some \( n > 1 \)) and the embedding \( f \).
Then there is \( i < n \) such that \( C_i \) is infinite, which means that \( \{ (\ell,q'') \in C_{(0,+\infty)} \mid q < q'' < q' \} \subseteq f(C_i) \csub C_{(0,+\infty)} \) for some rational numbers \( 0 \leq q < q' \).
Thus \( C \) contains a convex subset isomorphic to \( C_{(q,q')} \) by Lemma~\ref{prop:embeddingandconvex}. Pick \( m \in \N \) with \( m > q' \). Then \( C \npccvx C_{(m,+\infty)}\) because otherwise \( C_{(q,q')} \ccvx C \pccvx C_{(m,+\infty)} \), contradicting (the proof of) Proposition~\ref{prop:antichainsCO}\ref{prop:antichainsCO-a}.
\end{proof}

\begin{prop} \label{prop:antichains2CO}
Every \( \pccvx \)-antichain is contained in a \( \pccvx \)-antichain of size \( 2^{\aleph_0}\). In particular, there are no maximal \( \pccvx \)-antichains of size smaller than \( 2^{\aleph_0} \), and every \( C \in \CO \) belongs to a \( \pccvx \)-antichain of size \( 2^{\aleph_0}\).
\end{prop}

\begin{proof} 
Following the proof of Proposition~\ref{prop:everylobelongstoantichain}, we only need to verify that for every \( C \in \CO \) the set
$\{ S \subseteq \N \mid S \text{ is infinite} \wedge C_S \pccvx C \}$ is countable, where the $C_S$'s are defined in the proof Proposition \ref{prop:basisCO}\ref{prop:basisCO-1}. 

First observe that arguing as at the beginning of that proof and using Proposition~\ref{prop:embeddingandconvex} one can prove that if \( C_S \pccvx C\) then \( C_S \ccvx C \). 
Indeed, let \( C_S \pccvx C\) be witnessed by the finite convex partition \( \{ C_i \mid i < n \} \) (for some \( n > 1 \)) of \( C_S \) and the embedding \( f \colon C_S \to C \). 
Then some \( C_i \) must be infinite, so there is an embedding \( g \colon C_S \to C_S \) such that \( \Ima g \csube C_S \) and \( \Ima g \subseteq C_i \csub C_S \). 
Hence \( f(\Ima g) \csube C\), and so \( f \circ g\) witnesses \( C_S \ccvx C \). 

Suppose that \( S,S' \subseteq \N \) are distinct infinite sets such that \( C_S \pccvx C\) and \( C_{S'} \pccvx C \) via corresponding embeddings \( f \) and \( g \), respectively. 
Without loss of generality, we may assume that \( \Ima f \neq C\) and \( \Ima g \neq C\), as otherwise \( C_{S'} \pccvx  C_S\) or \( C_S \pccvx C_{S'} \), contradicting (the proof of) Proposition~\ref{prop:basisCO}\ref{prop:basisCO-1}.
If \( \Ima f \cap \Ima g \neq \emptyset \), then by Proposition~\ref{prop:intconvex} such intersection is the union of (at most) two proper convex subsets \( A_0, A_1 \) of \( C \), each of which must be infinite by definition of \( C_S \) and \( C_{S'}\). 
Thus \( f^{-1}(A_0) \) is an infinite convex proper subset of \( C_S \), and so \( C_S \ccvx f^{-1}(A_0) \), which in turn implies \( C_S \ccvx A_0 \) and \( C_S \ccvx C_{S'} \), a contradiction. Thus \( \Ima f \cap \Ima g = \emptyset \). Since \( C \) is countable, there can be at most countably many infinite \( S \subseteq \N \) such that \( C_S \pccvx C \) and the claim follows.
\end{proof}

Once we know that for every \( C \in \CO \) there are at most countably many infinite sets \( S \in \N \) such that \( C_S \pccvx C \), arguing as in Proposition~\ref{prop:dom_fam} we easily get

\begin{prop} \label{prop:dominatingCO}
\( \mathfrak{d}(\pccvx) = 2^{\aleph_0} \).
\end{prop}

We now move to the study of
the 
(analytic) equivalence relation  $\underline{\bowtie}^{<\omega}_{\CO}$ induced by \( \pccvx \). 
Obviously if $C \cong_{\CO} C'$ then
we also have $C \mathrel{\underline{\bowtie}^{<\omega}_{\CO}} C'$.

\begin{thm} \label{thm:COpiecewise}
${\iso} \leq_B {\underline{\bowtie}^{<\omega}_{\CO}}$.
\end{thm}

\begin{proof}
Consider the Borel map $\varphi\colon \LO \to \CO$ defined by
\[
\varphi(L)= C[(\boldsymbol{1}+\z L)\omega].
\]  
We claim that $\varphi$ is a reduction. Clearly, if \( L \iso L' \) then \( \varphi(L) \cong_{\CO} \varphi(L') \) and hence \( \varphi(L) \mathrel{\underline{\bowtie}^{<\omega}_{\CO}} \varphi(L') \). For the converse, let the finite convex partition \( \{ C_i \mid i<n \}\) of \(\varphi(L)\) and the embedding \( f \) of \( \varphi(L) \) into \(\varphi(L')\) witness $\varphi(L) \pccvx \varphi(L')$. 
Without loss of generality \( n > 1\), so that \( C_i \csub \varphi(L) \) for all \( i < n \).
Since \( n \) is finite, there exists some \(j < n\) such that $C_j$ contains at least two copies of $\boldsymbol{1}+\z L$, so we can consider a convex set of the form \(\boldsymbol{1}+\z L+\boldsymbol{1} \subseteq C_j\), so that \( f(\boldsymbol{1}+\z L+\boldsymbol{1}) \csube \varphi(L') \) by Proposition~\ref{prop:embeddingandconvex}. 
Since the $\boldsymbol{1}$'s are the only elements which do not have immediate predecessor and successor both in 
\( \varphi(L) \) and in $\varphi(L')$, and since \( f(\boldsymbol{1}+\z L+\boldsymbol{1}) \) is convex, we have that the images via $f$ of the two $\boldsymbol{1}$'s in $\boldsymbol{1}+\z L+\boldsymbol{1} \subseteq C_j$ are two necessarily ``consecutive'' $\boldsymbol{1}$'s in $\varphi(L')$. It follows that \( \boldsymbol{1}+\z L+\boldsymbol{1} \subseteq C_j\) is isomorphic to a copy of $\boldsymbol{1}+\z L'+\boldsymbol{1}$ in $\varphi(L')$. We thus obtain $\z L\iso \z L'$, hence $L\iso L'$ by Lemma~\ref{lem:isom_zetaL}.
\end{proof}

The next results contrasts with Corollary \ref{E1_cvxeq}. To simplify the notation, we let \( \vec{x} \) and \( \vec{y} \) denote the sequences \( (x_n)_{n \in \N} , (y_n)_{n \in \N} \in \R^\N \), respectively.

\begin{thm}\label{no_action_group}
\({E_1} \leq_B {\underline{\bowtie}^{<\omega}_{\CO}}\).
\end{thm}

\begin{proof}
By Proposition \ref{E1} it suffices to define a Borel reduction from \(E_1^t\) to \({\underline{\bowtie}^{<\omega}_{\CO}}\). To this end, fix an injective map \(f \colon \Q \to \{\boldsymbol{n} \mid n \in \N \setminus \{ 0,1 \} \}\) and, as in the proofs of Lemma \ref{open_int} and Proposition~\ref{prop1}, consider the linear orders \(\eta_f\) and \(L_{(x,x+1)}\), with \(x \in \R\). By Lemma \ref{open_int}, \(L_{(x,x+1)}\) and \(L_{(x',x'+1)}\) are isomorphic if and only if \(x=x'\). Consider the Borel map that sends \( \vec{x} = (x_n)_{n \in \N} \in \R^\N \) to the linear order
\[  
L(\vec{x}) = \sum_{n \in \Z} \overline{L}_n,
\]
where \( \overline{L}_n = \eta_f + \eta \) if \( n < 0 \) and \( \overline{L}_n = L_{(x_n,x_n+1)}+\eta \) if \( n \geq 0 \).
We claim that the Borel map \( \varphi \colon \R^\N \to \CO \) defined by $\varphi(\vec{x}) = C[L(\vec{x})]$ is a reduction from \(E_1^t\) to
\({\underline{\bowtie}^{<\omega}_{\CO}}\).

First suppose that \({\vec{x}} \mathrel{E_1^t} {\vec{y}}\), i.e.\ that there are \(\bar{n},\bar{m} \in \N\) such that \(x_{\bar{n}+k}=y_{\bar{m}+k}\) for all \(k \in \N \). 
Consider the finite convex partition \( \{ C_i \mid i < 2 \bar{n}+2 \} \) of \( \varphi(\vec{x})\) given by setting for \( 0 \leq j < \bar{n} \)
\begin{align*}
C_0 & = \sum_{n \in \Z \setminus \N} \overline{L}_n = \{ (\ell,n) \in L(\vec{x}) \mid n < 0 \} \\
C_{2j+1} & = L_{(x_j,x_j+1)} \times \{ j \} \\
C_{2j+2} & = \eta \times \{ j \} \\
C_{2 \bar{n} +1} & = \sum_{n \geq \bar{n}} \overline{L}_{n} = \{ (\ell,n) \in L(\vec{x}) \mid n \geq \bar{n} \}.
\end{align*}
Consider the embedding \( f \) of \( \varphi(\vec{x}) \) into \( \varphi(\vec{y}) \) defined by
\[  
f(\ell,n) = 
\begin{cases}
(\ell, n - \bar{n}) & \text{if } n < \bar{n} \\
(\ell, \bar{m} + (n - \bar{n})) & \text{if } n \geq \bar{n}.
\end{cases}
\]
By choice of \( \bar{n}, \bar{m} \in \N \) and
since \( L_{(x,x+1)} \csube \eta_f\) for all \( x \in \R \), it is easy to verify that \( f \) is well-defined  and that \( f(C_i) \csube \varphi(\vec{y}) \) for all \( i < 2 \bar{n} + 2\). This witnesses 
\( \varphi(\vec{x}) \mathrel{\trianglelefteq^{<\omega}_{\CO}} \varphi(\vec{y})\), and since \(\varphi(\vec{y}) \mathrel{\trianglelefteq^{<\omega}_{\CO}} \varphi(\vec{x})\) can be proved symmetrically, we obtain ${\varphi(\vec{x})} \mathrel{\underline{\bowtie}^{<\omega}_{\CO}}  {\varphi(\vec{y})}$.

Suppose now that%
\footnote{Our proof actually shows that 
\(\varphi(\vec{x}) \mathrel{\trianglelefteq^{<\omega}_{\CO}} \varphi(\vec{y})\) already suffices to obtain \( \vec{x} \mathrel{E^t_1} \vec{y} \), so that in particular we get
\( \varphi(\vec{x}) \mathrel{\trianglelefteq^{<\omega}_{\CO}} \varphi(\vec{y}) \iff {\varphi(\vec{x})} \mathrel{\underline{\bowtie}^{<\omega}_{\CO}} {\varphi(\vec{y})}  \)
for all \( \vec{x},\vec{y} \in \R^\N\).} 
\({\varphi(\vec{x})} \mathrel{\underline{\bowtie}^{<\omega}_{\CO}} {\varphi(\vec{y})}\). Let \( \{ C_i \mid i<b \}\) with 
\( b \in \N \setminus \{ 0 \} \) be a finite convex partition of \( \varphi(\vec{x}) \) and \( f \) be an embedding of $\varphi(\vec{x})$ into \(\varphi(\vec{y})\) witnessing \(\varphi(\vec{x}) \mathrel{\trianglelefteq^{<\omega}_{\CO}} \varphi(\vec{y})\). (As usual, we can assume \( b > 1 \), so that Proposition \ref{prop:embeddingandconvex} can be applied when necessary.)
Since \( b \) is finite, for some \( i < b \) and \( \bar{n} \in \N \setminus \{ 0 \} \) we must have \( \sum_{n\geq \bar{n}-1} \overline{L}_{n} \subseteq C_i \). 
Notice that for every \( n \geq \bar{n}-1 \) and \( q \in \eta \), the point \( (q,n) \in \eta \times \{ n \} \csub \varphi(\vec{x}) \) has no immediate predecessor and immediate successor, while points of the form \( (\ell,m) \) for \( \ell \in L_{(y_m,y_m+1)} \) and \( m \in \N \) or \( \ell \in \eta_f \) and \( m \in \Z \setminus \N \) have an immediate predecessor or an immediate successor (or both): thus \( f(q,n) \in \eta \times m\) for some \( m \in \Z \). 
By a similar argument, \( f(L_{(x_n,x_n+1)} \times \{ n \}) \subseteq L_{(y_m,y_{m+1})} \times \{ m \}\) or \( f(L_{(x_n,x_n+1)} \times \{ n \}) \subseteq \eta_f \times \{ m \} \) for a suitable \( m \in \Z \). This two facts together with the convexity of \( f(C_i) \) and the fact that, by the proof of Lemma~\ref{open_int}, the only convex subset of \( \eta_f \) isomorphic to \( L_{(x,x+1)}\) is \( L_{(x,x+1)}\) itself, imply that \( f(L_{(x_{\bar{n}},x_{\bar{n}+1})} \times \{ \bar{n} \}) = L_{(x_{\bar{m}}, x_{\bar{m}+1})} \times \{ \bar{m} \}\) for some \( \bar{m} \in \N \), and in turn \( f(L_{(x_{\bar{n}+k},x_{\bar{n}+k+1})} \times \{ \bar{n}+k \}) = L_{(x_{\bar{m}+k}, x_{\bar{m}+k+1})} \times \{ \bar{m}+k \}\) for all \( k \in \N \). But by Lemma~\ref{open_int} again, this means that
\(x_{\bar{n} + k}=y_{\bar{m}+k}\) for all \( k \in \N \), hence \( \vec{x}  \mathrel{E^t_1} \vec{y} \).
\end{proof}

\begin{cor} \label{cor:piecewiseiscomplex}
\({\iso} <_B {\underline{\bowtie}^{<\omega}_{\CO}}\) and \(\cvxeq <_{\text{\scriptsize \textit{Baire}}} {\pccvxeq}\). Moreover \({\underline{\bowtie}^{<\omega}_{\CO}}\) is not Baire reducible to
an orbit equivalence relation..
\end{cor}

\begin{proof}
All the statements follow from Theorem \ref{no_action_group} and some of the previous results. The first two statements need Theorem \ref{thm:COpiecewise} and Corollaries \ref{E1_cvxeq} and \ref{cvxeq_baire_iso}; the last one follows from Theorem \ref{E1_orbit}.
\end{proof}

\section{Anti-classification results in knot theory}

In this section we recall the basic notions of knot theory and the
relation between proper arcs/knots and linear orders established in
\cite{Kul17}. We analyse the relation of subarc among proper arcs 
(Definition~\ref{def 2.1.2})
using the results and methods developed for convex embeddability -- its
counterpart in the realm of linear orders. We then move to knots and
define a natural notion of subknot (Definition~\ref{defn0 4}), uncovering a natural 
connection with circular orders. This allows us to show that the
equivalence relation associated to the latter quasi-order is not induced by a Borel 
action of a Polish group, in strong contrast with knot equivalence.

\subsection{Knots and proper arcs: definitions and basic facts} \label{sec:introknotsarcs}

In mathematics, there are essentially two ways to formalize the intuitive concept of a knot: a (mathematical) \emph{knot} is obtained from a real-life knot by joining its ends so that it cannot be undone, while a \emph{proper arc} is obtained by embedding the real-life knot in a closed 3-ball and sticking its ends to the border of the ball, so that again it cannot be undone. The two concepts are strictly related, although not equivalent as there exist knots that cannot be \lq\lq cut\rq\rq\ to obtain a proper arc (\cite{Bi56}).  Let us recall the main definitions and related concepts.

Depending on the situation, we think of \( S^1 \) as either the unit circle in $\R^2$ or the one-point compactification of \( \mathbb{R}\) obtained by adding $\infty$ to the space. 
Similarly, \( S^3 \) can be viewed as the one-point compactification of \( \mathbb{R}^3\).

\begin{defn}
A \textbf{knot} \( K \) is a homeomorphic image of $S^1$ in $S^3$, that is, a subspace of \( S^3 \) of the form \( K = \Ima f\) for some topological embedding \( f \colon S^1 \to S^3 \).
\end{defn}

The
collection of all knots is denoted by~$\Kn$;
as shown in \cite{Kul17}, it can be construed as a standard Borel space.
Obviously, if \( K \subseteq S^3 \) is a knot and \( \varphi \colon S^3 \to S^3 \) is an embedding, then \( \varphi(K) \) is a knot as well.

\begin{remark}\label{rem:circ_ord_on_knots}
Knots can be naturally endowed with a circular order induced by the standard circular order $C_{S^1}$ defined on $S^1$ (see Section \ref{sec:circular orders}). More precisely, let $f \colon S^1 \to S^3$ be an embedding and $K = \Ima f$ be the knot induced by \( f \). 
Then for every $x,y,z \in K$ we can set 
\[  
C_f(x,y,z) \iff C_{S^1}(f^{-1}(x),f^{-1}(y),f^{-1}(z)).
\]
If \( f , f' \colon S^1 \to S^3 \) are two embeddings giving rise to the same knot \( K = \Ima f = \Ima f' \), then \( f^{-1} \circ f' \colon S^1 \to S^1 \) is a homeomorphism, and thus it is either order-preserving or order-reversing with respect to \( C_{S^1} \). It follows that either \( C_f = C_{f'} \) or \( C_f = C_{f'}^* \). 
Thus a knot \( K \) can be endowed with exactly two circular orders, corresponding to the two possible orientations of $K$ sometimes used in knot theory, which are one the reverse of the other one and depend on the specific embedding used to witness \( K \in \Kn \).
We speak of \textbf{oriented} knot \( K \) when we single out one specific orientation between the two possibilities.
\end{remark}

Two knots $K, K' \in \Kn$ are \textbf{equivalent}, in
symbols
\(
{K} \equiv_{\Kn} {K'}, 
\)
if there exists a homeomorphism
$\varphi\colon S^3\to S^3$ such that $\varphi(K) = K'$.
%
The relation \( \equiv_{\Kn} \) is an analytic equivalence
relation on \( \Kn \). A knot is \textbf{trivial} if it is equivalent to the unit circle  $I_{\Kn} = \{(x,y,z) \in S^3 \mid x^2+y^2=1 \wedge z=0\}$.

\begin{remark}
In knot theory it is more common to consider the oriented version of \( \equiv_{\Kn} \), according to which two knots \( K \) and \( K' \) are equivalent if there is an \emph{orientation-preserving} homeomorphism \( \varphi \colon S^3 \to S^3 \) such that \( \varphi(K) = K' \) or, equivalently, an ambient isotopy sending \( K \) to \( K' \). Nevertheless, we are mostly going to prove anti-classification results, and thus they become even stronger if we consider the coarser equivalence relation \( \equiv_{\Kn} \). For the interested reader,  however, we point out that all our results remain true if we stick to common practice and replace all the relevant equivalence relations and quasi-orders with their oriented versions. Similar considerations apply to the ensuing definitions and results concerning proper arcs.
\end{remark}

We now move to proper arcs. A proper arc is a special case of an $n$-tangle,
namely it is a $1$-tangle. Contrary to most expositions \cite{Conway1970AnEO,Khovanov2002}, however, we
do not assume that the tangles (arcs) are tame. 
Given \( x \in \R^3 \) and a positive \( r \in \R \), the closed ball with center \( x \) and radius \( r \) is denoted by \( \bar{B}(x,r) \). The origin \( (0,0,0) \) of \( \R^3 \) is sometimes denoted by \( \bar{0} \). 
To avoid repetitions, we convene that from now $\bar{B}$, possibly with subscripts and/or superscripts, is always a closed topological 3-ball, i.e.\ a homeomorphic copy of a closed ball in $\mathbb{R}^3$. Recall that by compactness of \( \bar{B} \)
and the invariance of domain theorem, the notion  of boundary of \( \bar{B} \) as a topological subspace of \( \R^3 \) and the notion of boundary of \( \bar{B}\) as a topological \( 3 \)-manifold coincide. Thus we can unambiguously denote by \(\partial \bar{B} \) the \textbf{boundary} of \( \bar{B}\),
and set \( \Int \bar{B} = \bar{B} \setminus \partial \bar{B} \). Notice also that by the same reasons, if \( \varphi \colon \bar{B} \to \R^3\) is an embedding, then \( \varphi(\partial \bar{B}) = \partial \, \varphi(\bar{B}) \) and \( \varphi (\Int \bar{B}) = \Int \varphi(\bar{B}) \).

\begin{defn} \label{def:arc}
Given a topological embedding $f \colon [0,1] \to \bar{B}$ we say that the pair $(\bar{B},\Ima f)$ is an \textbf{proper arc} if $f(x) \in \partial \bar{B} \iff {x = 0} \vee {x =1} $. With an abuse of notation which is standard in knot theory, when there is no danger of confusion \emph{we identify $f$ with its image \( \Ima f\)} and write e.g.\ $\arc{B}{f}$ in place of $(\bar{B},\Ima f)$. 
\end{defn}

Any proper arc \( \arc{B}{f} \) can be canonically turned  (up to knot equivalence) into a knot \( K_{\arc{B}{f}} \) by joining its ends \( f(0) \) and \( f(1) \) with a simple curve running on the boundary \( \partial \bar B\) of its ambient space, assuming, without loss of generality,
that $\bar B=\bar B(\bar 0,1)$.
The collection of proper arcs is denoted by~\(\Ar\), and
 can be construed as a standard Borel subspace of the product \( K(\R^3) \times K(\R^3)\) of the Vietoris space \( K(\R^3) \) over \( \R^3\).
 This follows as an application of a theorem by Ryll-Nardzewski~\cite{Ryll65}
 (see \cite{Kul17} for the analogous construction of the coding space \( \Kn\) of knots). 
Notice that if \( \arc{B}{f}\) is a proper arc and \( \varphi \colon \bar{B} \to \R^3 \) is an embedding, then \( \barc{\varphi(\bar{B})}{\varphi(f)} =  \barc{\varphi(\bar{B})}{\varphi(\Ima f)}\) is a proper arc, as witnessed by the embedding \( \varphi \circ f \colon [0,1] \to \varphi(\bar{B})\).

\begin{remark} \label{rmk:orientation}
\begin{enumerate-(1)}
\item \label{rmk:orientation-1}
Every specific embedding \( f \) giving rise to an arc \( (\bar{B},\Ima f) \) induces an orientation on it, namely, the linear order \( \leq_f \) on \( \Ima f \) defined by 
\[
b_0 \leq_f b_1 \iff f^{-1}(b_0) \leq f^{-1}(b_1).
\]
If \( f,f' \colon [0,1] \to \bar{B} \) are two topological embeddings inducing the same proper arc (that is, \( \Ima f = \Ima f' \)), then \( f^{-1} \circ f' \colon [0,1] \to [0,1] \) is a homeomorphism, and thus it is either order-preserving or order-reversing.
It follows that every proper arc has exactly two orientations. Moreover,  the minimum and the maximum of \( \leq_f \) always exist; they can be identified, independently of \( f \), as the only points of \( \Ima f \) belonging to \( \partial \bar{B}\). We speak of \textbf{oriented} proper arc \( \arc{B}{f}\) when we equip it with the specific orientation given by the displayed \( f \).

\item \label{rmk:orientation-2}
If \( \arc{B}{f} \) and \( \arc{B'}{g}\) are proper arcs and \( \varphi \colon \bar{B} \to \bar{B}'\) is a topological embedding such that \( \varphi(\Ima f) \subseteq \Ima g\), then \( h = g^{-1} \circ \varphi \circ f \colon [0,1] \to [0,1] \) is a topological embedding. 
It follows that when \( \arc{B}{f} \) and \( \arc{B'}{g}\) are construed as oriented proper arcs, then \( \varphi \) is either order-preserving (that is, \( \varphi(b_0) \leq_g \varphi(b_1)\) for all \( b_0,b_1 \in \bar{B}\) with \( b_0 \leq_f b_1\)) or order-reversing (that is, \( \varphi(b_1) \leq_g \varphi(b_0)\) for all \( b_0,b_1 \in \bar{B}\) with \( b_0 \leq_f b_1\)).
\end{enumerate-(1)}
\end{remark}

Two proper arcs $\arc{B}{f}$ and $ \arc{B'}{g} $ are
\textbf{equivalent}, in symbols
\(
\arc{B}{f}\equiv_{\Ar}\arc{B'}{g},
\)
if there exists a homeomorphism
$\varphi\colon \bar{B}\to \bar{B}'$ such that $\varphi(\Ima f)= \Ima g$.
%
The relation \( \equiv_{\Ar} \) is an analytic equivalence relation on the standard Borel space \( \Ar \). A proper arc \(\arc{B}{f}\) is \textbf{trivial} if it is equivalent to  $I_{\Ar} = (\bar{B}(\bar{0}, 1), [-1, 1] \times \{(0, 0)\})$.

An important dividing line among knots (respectively, proper arcs) is given by tameness, i.e.\ the absence of singular points. Given a knot \( K \in \Kn \), 
a \textbf{subarc} of \( K \) is any proper arc of the form 
\( \arc{B}{K \cap \bar{B}} \).
A point \( x \in K \) is called \textbf{singular}, or a \textbf{singularity}, 
of \( K \) if there is no \( \bar B \) such that \( x \in \Int \bar B\) and 
\( \arc{B}{K \cap \bar B} \) is a trivial proper subarc of \( K \). 
The space of singularities of \( K \) is denoted by \( \Sigma_K \).
An \textbf{isolated} singular point
of \( K \) is an isolated point of the topological
space \( \Sigma_K \), and
the (sub)space of isolated singular points of \( K \) is denoted by \( I\Sigma_K \).
Finally, a knot \( K \) is \textbf{tame}%
\footnote{Our definition of tame knot is equivalent to the classical one, according to which a knot is tame if it is equivalent to a finite polygon (see \cite[Definition 1.3]{BZ03}). }
if it has no singular
points, and \textbf{wild} otherwise.
Notice also that if \( x \in K \) is not a singularity of \( K \), then there are arbitrarily small closed topological \( 3 \)-balls \( \bar B \) witnessing this.

The previous definitions can be naturally adapted to proper arcs.
%
Let \( \arc{B}{f} \in \Ar\).
A point \( x \in \Ima f \) is called
\textbf{singular}, or a \textbf{singularity}, of \( \arc{B}{f}\) 
if it belongs to \( \Sigma_{K_{\arc{B}{f}}}\), 
while an \textbf{isolated} singular point of \( \arc{B}{f} \) is an an element of \( I\Sigma_{K_{\arc{B}{f}}} \).
Accordingly, the space of singularities of \( \arc{B}{f} \) is denoted by \( \Sigma_{\arc{B}{f}}\), while the space of isolated singular points is denoted by \( I\Sigma_{\arc{B}{f}}\).
An arc \( \arc{B}{f} \) is \textbf{tame} if \( \Sigma_{\arc{B}{f}} = \emptyset \) (equivalently, if \( K_{\arc{B}{f}} \) is tame), and \textbf{wild} otherwise.
%
Notice that if \( x \in \Ima f \cap \Int \bar B\), then \( x \notin \Sigma_{\arc{B}{f}}\)
if and only if  there is \( \bar{B}' \subseteq \bar{B}\) such that \( x \in \Int \bar{B}'  \) and 
\( \arc{B'}{f \cap \bar{B}'} \) is a trivial proper arc. For points on the boundary \( \partial \bar B \), instead, it is not enough to consider closed topological \( 3 \)-balls \( \bar B' \subseteq \bar B\), as we necessarily need to consider a ``trivial prolungation'' of the curve \( \Ima f \) beyond its extreme points in order to determine whether they are singular or not.

We also introduce the notion of circularization of a proper arc, which generates a knot and gives a characterization of tame knots.

\begin{defn}\label{def:circularization_of_arcs}
Let $\arc{B}{f} \in \Ar$. Up to equivalence, we can assume that $\bar B =[-1,1]^3$, $f(0)=(-1,0,0)$ and $f(1)=(1,0,0)$. Consider the equivalence relation obtained setting $(-1,y,z) \sim (1,y,z)$ for all $(y,z) \in [-1,1]^2$, so that in the quotient space $T = [-1,1]^3 / {\sim}$ the two lateral faces of the cube $\bar B$ are glued and we have a solid torus.
Given a topological embedding $h$ of $T$ into $S^3$, we call \textbf{circularitazion of} $\arc{B}{f}$, denoted by $C^h[\arc{B}{f}]$, the knot which is obtained as the image of $\Ima f / {\sim}$ via $h$.    
\end{defn}

Notice that the circularization of a proper arc depends on the topological embedding of the solid torus into $S^3$, hence it is not unique. Moreover, we have that $K \in \Kn$ is tame if and only if $K=C^h[I_{\Ar}]$ for some topological embedding $h \colon T \to S^3$.

A substantial part of the analysis of tame knots relies on their prime factorization, which is in turn based on the classical notion of sum (see \cite[Definition 2.7]{BZ03}, where the sum is actually called product). 
We will work with the corresponding sum for proper arcs which 
is akin to the tangle sum~\cite{Conway1970AnEO} (note that we allow wild arcs).

\begin{defn} \label{def:sumarc}
Let \( \arc{B_0}{f_0} \) and \( \arc{B_1}{f_1} \) be \emph{oriented} proper arcs. Up to equivalence, we may assume that \( \bar{B}_0 = [-1,0] \times [-1,1]^2\), \( \bar{B}_1 = [0,1] \times [-1,1]^2\), \( f_0(0) = (-1,0,0)\), \( f_0(1) = f_1(0) = (0,0,0)\), and \( f_1(1) = (1,0,0) \).
The \textbf{sum} 
\(
\arc{B_0}{f_0} \oplus \arc{B_1}{f_1} 
\) 
is the proper arc \( \arc{B}{f} \) where \( \bar{B} = [-1,1]^3 \) and \( f \colon [0,1] \to \bar{B} \)
is defined by \( f(x) = f_0(2x)\) if \( x \leq \frac{1}{2}\) and \( f(x) = f_1(2x-1) \) if \( x \geq \frac{1}{2} \).
\end{defn}

By induction on \( n \in \N \), one can then define finite sums of proper arcs
\(
\arc{B_0}{f_0} \oplus \dots \oplus \arc{B_n}{f_n} ,
\)
abbreviated by \( \bigoplus_{i \leq n} \, \arc{B_i}{f_i} \).

\begin{remark} 
Although the sum of two oriented proper arcs is again oriented, in this paper we will tacitly consider it as an unoriented proper arc. Also, we will often sum unoriented proper arcs: what we mean in this case is that the arcs are summed using the natural orientation coming from the way we present them.
\end{remark}


Tame knots and arcs are in canonical one-to-one correspondence up to
equivalence. Given a tame knot $K$ in $S^3$, pick 
an arbitrary $x\in K$
and a neighbourhood $B_0$ of $x$ such that $(\bar B_0,K\cap B_0)$ is a trivial 
arc. Then
let $\bar B_K=S^3\setminus B_0$ and $f_K=f\setminus B_0$. Then
$K_{(\bar B_K,f_K)}$ is equivalent to~$K$. Vice versa, for 
all equivalent arcs $\arc{B}{f}$ and $\arc{B'}{f'}$ it is easy to see 
that $K_{\arc{B}{f}} \equiv_{\Kn} K_{\arc{B'}{f'}}$. This transformation commutes with the 
knot sum \( \# \) of tame knots~\cite[Definition~7.1]{BZ03}. The latter
can actually be defined using~\(\oplus\). 
Given two oriented tame knots we then have that, up to equivalence, 
\[K_0 \mathrel{\#} K_1 = K_{\arc{B_{K_0}}{f_{K_0}} \oplus \arc{B_{K_1}}{f_{K_1}}}.\] 
Notice also that if \( \arc{B_0}{f_0} \) and \( \arc{B_1}{f_1} \) are (oriented) tame proper arcs, then 
\[K_{\arc{B_0}{f_0} \oplus \arc{B_1}{f_1}} \equiv_{\Kn} K_{\arc{B_0}{f_0}} \mathrel{\#} K_{\arc{B_1}{f_1}}.\]

Recall that a nontrivial tame knot is \textbf{prime} if it cannot be written as a sum of nontrivial knots.
The prime knot decomposition 
theorem  \cite[Theorem 7.12]{BZ03} states that every
tame knot can be expressed as a finite knot sum of prime knots in a unique way up to knot equivalence.
We will later consider also prime arcs which are defined in a similar way.

\subsection{Proper arcs and their classification}\label{sec:arcs}

The following notion is  equivalent to~\cite[Definition 2.10]{Kul17}.



\begin{defn}{\label{def 2.1.2}}
Let $\arc{B}{f}, \arc{B'}{g} \in \Ar$. We say that
\( \arc{B}{f} \) is a \textbf{subarc} of $\arc{B'}{g}$, or that
$\arc{B'}{g}$ has $\arc{B}{f}$ as a subarc,
if there exists a topological embedding
$\varphi\colon \bar{B}\to \bar{B}'$ such that
$\varphi(f) = g \cap \Ima\ \varphi$. 
In this case we write 
\begin{equation*}    
\arc{B}{f}\precsim_{\Ar} \arc{B'}{g}.
\end{equation*} 
(Notice that we automatically have that \( \barc{\varphi(\bar{B})}{g \cap \Ima \varphi)} \) is a proper arc.)
\end{defn}



Clearly, the subarc relation $\precsim_{\Ar}$ is
an analytic quasi-order on the standard Borel space $\Ar$. 
We denote by \( \prec_{\Ar}\) the strict part of \( \precsim_{\Ar}\), i.e.\ 
\[
\arc{B}{f} \prec_{\Ar} \arc{B'}{g} \iff {\arc{B}{f} \precsim_{\Ar} \arc{B'}{g}} \wedge {\arc{B'}{g} \not\precsim_{\Ar} \arc{B}{f}} .
\]
The
analytic equivalence relation associated to $\precsim_{\Ar}$
is denoted by $\approx_{\Ar}$, and we say that two proper arcs
$\arc{B}{f}$ and $\arc{B'}{g}$ are \textbf{mutual subarcs} if
\[{\label{eq}}
\arc{B}{f} \approx_{\Ar} \arc{B'}{g}.
\] 
This may be interpreted as asserting that the two
arcs have the ``same complexity'' because each of
them is a subarc of the other one. Notice also that $\arc{B}{f}\equiv_{\Ar} \arc{B'}{g}$ trivially implies $\arc{B}{f}\approx_{\Ar} \arc{B'}{g}$.


If \( ( \bar{B},f), \arc{B'}{g} \in \Ar\) and 
\( \varphi \) witnesses \( \arc{B}{f} \equiv_{\Ar} \arc{B'}{g} \), then \( \varphi \) induces a homeomorphism between the spaces \( \Sigma_{\arc{B}{f}} \) and \( \Sigma_{\arc{B'}{g}}\), and hence also a homeomorphism between \( I\Sigma_{\arc{B}{f}}\) and \( I\Sigma_{\arc{B'}{g}} \).
If instead
\( \varphi \colon \bar{B} \to \bar{B}'\) is just an embedding witnessing \( \arc{B}{f} \precsim_{\Ar} \arc{B'}{g}\), then we still have that \( \varphi \) induces an embedding of \( \Sigma_{\arc{B}{f}}\)into \( \Sigma_{\arc{B'}{g}}\), but needs not send isolated singular points into isolated singular points: if \( x \in I\Sigma_{\arc{B}{f}} \cap \partial \bar{B} \), then it might happen that \( \varphi(x) \in \Sigma_{\arc{B'}{g}} \setminus I\Sigma_{\arc{B'}{g}} \). However, this is the only exception.

\begin{lem} \label{lem:isolatedthroughembeddings}
\begin{enumerate-(a)}
\item \label{lem:isolatedthroughembeddings-1}
Let \( \arc{B}{f} \in \Ar\) and \( \bar{B}' \subseteq \bar{B} \) be such that \( \arc{B'}{f \cap \bar{B}'} \in \Ar \). Then 
\( \Sigma_{\arc{B'}{f \cap \bar B'}} \subseteq \Sigma_{\arc{B}{f}}\), and \( \Sigma_{\arc{B'}{f \cap \bar{B}'}} \cap \Int \bar{B}' = \Sigma_{\arc{B}{f}} \cap \Int \bar{B}' \).

\item \label{lem:isolatedthroughembeddings-2}
Let \( ( \bar{B},f), \arc{B'}{g} \in \Ar\), and let \( \varphi \colon \bar{B} \to \bar{B}'\) witness \( \arc{B}{f} \precsim_{\Ar} \arc{B'}{g}\).
If \( x \in I\Sigma_{\arc{B}{f}} \cap \Int \bar{B} \), then \( \varphi(x) \in I\Sigma_{\arc{B'}{g}}\).
\end{enumerate-(a)}
\end{lem}

\begin{proof}
\ref{lem:isolatedthroughembeddings-1} 
The first part is easy and is left to the reader. 
For the nontrivial inclusion of the second part,
%
assume that \( x \in \Int \bar B' \) (so that \( x \in \Int \bar B \) as well because \( \bar B' \subseteq \bar B \)) and \( x \notin \Sigma_{\arc{B'}{f \cap \bar{B}'}}\). Let \( \bar{B}'' \subseteq \bar{B}'\) be a witness of this: then \( \bar{B}''\) also witnesses \( x \notin \Sigma_{\arc{B}{f}}\).

\ref{lem:isolatedthroughembeddings-2}
By hypothesis and the fact that \( \varphi \colon \bar{B} \to \varphi(\bar{B}) \) is a homeomorphism, \( \varphi(x) \in I\Sigma_{\barc{\varphi(\bar{B})}{g \cap \Ima \varphi}} \cap \Int \varphi(\bar{B}) \). By part~\ref{lem:isolatedthroughembeddings-1}, this implies that \( \varphi(x) \in \Sigma_{\arc{B'}{g}}\). Using \( \varphi(x) \in \Int \varphi(\bar{B}) \), pick a small enough open set \( U \subseteq \Int \varphi(\bar{B})  \) such that \( U \cap \Sigma_{\barc{\varphi(\bar{B})}{g \cap \Ima \varphi}}  = \{ \varphi(x) \} \): then by part~\ref{lem:isolatedthroughembeddings-1} again 
\( U \cap \Sigma_{\arc{B'}{g}}  = U \cap \Sigma_{\barc{\varphi(\bar{B})}{g \cap \Ima \varphi}} \), and thus \( U \cap \Sigma_{\arc{B'}{g}} \) witnesses \( \varphi(x) \in I\Sigma_{\arc{B'}{g}}\).
\end{proof}

We now define an infinitary version of the sum operation for (tame) proper arcs introduced in Definition~\ref{def:sumarc}.
Since the ambient space \( \bar{B}\) in the definition of a proper arc is a compact space, in order to define such infinitary sums we need the summands to accumulate towards a point  \( b \in \bar{B} \), which thus becomes a singularity when infinitely many summands are not trivial. 

\begin{defn} \label{def:infinitesumarc}
Let \( \arc{B_i}{f_i} \) be oriented proper arcs, for \( i \in \N \).%
\footnote{When summing unoriented proper arcs, if not specified otherwise we use the natural orientation coming from their presentation.} 
The (\textbf{infinite}) \textbf{sum with limit \( b \in \bar{B} \)}, denoted by
\( \bigoplus_{i \in \N}^b \, \arc{B_i}{f_i} \), 
is defined up to equivalence as follows. Without loss of generality, we may assume that \( b \) is of the form \( (b',0,0)\) for some \( b' \) with \( 0 < b' \leq 1  \). Up to equivalence, we may also assume that \( \bar{B}_i = [b' - 2^{-i}, b'-2^{-(i+1)}] \times [-2^{-i},2^{-i}]^2\), and that \( f_i(0) = (b'-2^{-i},0,0) \) and \( f_i(1) = (b'-2^{-(i+1)},0,0) \) for all \( i \in \N \). Then \( \bigoplus_{i \in \N}^b \, \arc{B_i}{f_i}\) is the arc \( \arc{B}{f} \) where \( \bar{B} = [-1,1]^3 \) and \( \Ima f \) is 
the union of \( \bigcup_{i \in \N } \Ima f_i\) together with \( [-1,b'-1] \times \{ (0,0) \} \) and \( [b',1] \times \{ (0,0) \}\) (the latter might reduce to the point \( (1,0,0) \) if \( b' = 1\) or, equivalently, if \( b \in \partial \bar{B} \)).
\end{defn}

Trivially, \( \arc{B_j}{f_j} \precsim_{\Ar} \bigoplus_{i \leq n} \, \arc{B_i}{f_i} \) for all \( n \geq j \) and
\( \arc{B_j}{f_j} \precsim_{\Ar} \bigoplus_{i \in \N}^b \arc{B_i}{f_i} \) for all \( b \in \bar{B} \).
Notice that, up to \( \equiv_{\Ar}\), Definition~\ref{def:infinitesumarc} gives rise to precisely two non-equivalent proper arcs, depending on whether \( b \in \partial \bar{B} \) or not---besides this dividing line the actual choice of the limit point \( b \in \bar{B} \) is completely irrelevant. Therefore we can simplify the notation by denoting with \( \bigoplus_{i \in \N} \, \arc{B_i}{f_i}\) the infinite sum \( \bigoplus^b_{i \in \N} \, \arc{B_i}{f_i}\) for some/any \( b \in \Int \bar{B} \), and with \( \bigoplus^\partial_{i \in \N} \, \arc{B_i}{f_i}\) the infinite sum \( \bigoplus^b_{i \in \N} \, \arc{B_i}{f_i}\) for some/any \( b \in \partial \bar{B} \). It is not hard to see that \( \bigoplus^\partial_{i \in \N} \, \arc{B_i}{f_i} \prec_{\Ar} \bigoplus_{i \in \N} \, \arc{B_i}{f_i} \).
Finally, if all the proper arcs \( \arc{B_i}{f_i}\) are equivalent to the same arc \( \arc{B'}{g}\), the two possible infinite sums will be denoted by \( \bigoplus_{\N} \, \arc{B'}{g}\) and \( \bigoplus^\partial_{\N} \, \arc{B'}{g} \), respectively. Obviously, we can also replace \( \N \) with any infinite \( A \subseteq \N \) and write \( \bigoplus_{j \in A}^{(\partial)} \, \arc{B_j}{f_j} \) to denote \( \bigoplus_{i \in \N}^{(\partial)} \, \arc{B_{r(j)}}{f_{r(j)}}\), where \( r \colon \N \to A \) is the increasing enumeration of \( A \); similarly for \( \bigoplus_A^{(\partial)} \, \arc{B'}{g} \).

\begin{figure}[ht]
\centering \includegraphics[width=0.9\textwidth]{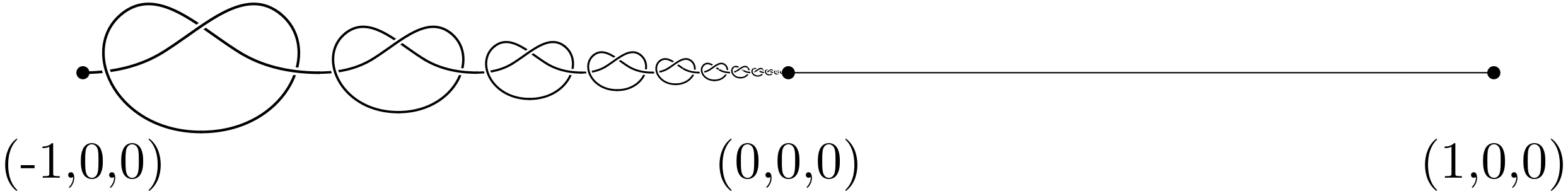}
\caption[LoF entry]{Infinite sum of trefoils, with limit point internal to the ambient space \( \bar{B} =  [-1,1]^3 \).}
\label{trefoils}
\end{figure}

Figure~\ref{trefoils} presents the arc \( \bigoplus_{\N} \, \arc{B'}{g}\) where \( \arc{B'}{g} \) is the trefoil; its variant \( \bigoplus^\partial_{\N} \arc{B'}{g} \) would be obtained my moving the current limit point \( (0,0,0)\) to the point \( (1,0,0) \) on \( \partial \bar{B} \).

In \cite[Theorem 3.1]{Kul17} it is shown that the isomorphism \( \cong_{\LO} \) on countable linear orders Borel reduces to equivalence \( \equiv_{\Kn} \) on knots. Employing the same construction, we establish a similar connection between convex embeddability  \( \cvx \) on linear orders and the subarc relation \( \precsim_{\Ar} \) on proper arcs.

Fix a proper arc \( \arc{B^*}{f^*} \) of the form \( \bigoplus_\N \, \arc{B_i}{f_i} \) with all the proper arcs \(\arc{B_i}{f_i}\) tame and not trivial. (For the sake of definiteness,  one can e.g.\ assume that \( \arc{B^*}{f^*}\)
is the sum of infinitely many trefoils depicted in Figure~\ref{trefoils}.) An important feature of such a \( \arc{B^*}{f^*} \) is that
\begin{equation} \tag{\( \dagger\)} \label{eq:not-inverting}
\begin{minipage}{11cm}
Any embedding \( \varphi \colon \bar{B^*} \to \bar{B^*}\) with \( \varphi(f^*) = f^* \cap \Ima\varphi \) preserves the (natural) orientation of the arc.
\end{minipage}
\end{equation}
Notice also that the only singularity of \( \arc{B^*}{f^*} \), which is trivially isolated, belongs to \( \Int \bar{B^*} \).

We first define a Borel map that given $L \in \LO$ produces an 
order-embedding $h_L$ of \( L \) into $(\Q, \leq)$ and a function $r_L \colon L \to \Q$ such that:
\begin{enumerate-(a)}
\item\label{Vpwdisj} the open intervals $V_n^L=(h_L(n)-2r_L(n), h_L(n)+2r_L(n))$ are included in $[-1,1]$ and pairwise disjoint;
\item\label{dense} 
$\bigcup_{n \in \N} V_n^L$ is dense in $[-1,1]$.
\end{enumerate-(a)}
To this end, we first establish in a Borel way whether $L$ has extrema, what are they are, and when one element of the linear  order is the immediate successor of another. So suppose that the maps $L\mapsto (h_L,r_L)$ are defined.

Notice that $\lim_{n \to \infty} r_L(n) = 0$ and that we can assume that $r_L(n+1)<r_L(n)$ for every $n\in \N$. Let $U_n^L = [h_L(n)-r_L(n), h_L(n)+r_L(n)]$.
Thinking of $[-1,1]$ as lying on the $x$-axis, we replace $U_n^L$ with the cube $\bar{B}_n^L = U_n^L \times [-r_L(n),r_L(n)]^2$. Let $\arc{B_n^L}{ f_n^L}$ be equivalent to $\arc{B^*}{f^*}$ and such that \( f_n^L(0) < f_n^L(1) \) and both belong to the \( x \)-axis, and set $f_r^L = ([-1,1] \setminus \bigcup_{n \in \N} U_n^L) \times \{(0,0)\}$. Then we define the map
\begin{equation} \label{eq:FfromLOtoA} 
F \colon \LO \to \Ar, \qquad L \mapsto \arc{B}{f_L}
\end{equation}
by letting
\(\bar{B}=[-1,1]^3\) and \(f_L = f_r^L \cup \bigcup_{n \in \N} f_n^L\).

By construction, every $(h_L(n),0,0)$ is singular and isolated in \( \Sigma_{F(L)} \) by (the trace of) $\bar{B}_n^L$, and every other member of \( \Sigma_{F(L)} \) is a limit of these singular points. Thus $\Sigma_{F(L)}$ is contained in the $x$-axis and $I\Sigma_{F(L)}=\{(h_L(n),0,0) \mid n \in \N\}$. The latter is naturally ordered by considering the restriction of \( \leq_{f_L} \) to \( I\Sigma_{F(L)} \), or equivalently, by considering first coordinates ordered as elements of \( \R \). Then the map \( n \mapsto (h_L(n),0,0) \) is an isomorphism between the linear orders \( L \) and \( I\Sigma_{F(L)} \).

Since the entire construction really depends on the proper arc \( \arc{B^*}{f^*}\), when relevant we will add this information to the notation and write e.g.\ \( F_{\arc{B^*}{f^*}}(L)\). For future reference, we also notice that by construction \( I\Sigma_{F(L)} \subseteq \Int \bar{B} \).

\begin{thm}\label{thm:reduc_lo_arcs}
The map \( F \) from~\eqref{eq:FfromLOtoA} simultaneously witnesses
${\trianglelefteq_{\LO}} \leq_B {\precsim_{\Ar}}$ (hence also ${\underline{\bowtie}}_{\LO} \leq_B {\approx_{\Ar}}$)
and \( {\cong_{\LO}} \leq_B {\equiv_{\Ar}}\). 
\end{thm}

The lower bound \( {\cong_{\LO}} \leq_B {\equiv_{\Ar}}\) for the relation \( \equiv_{\Ar} \) is implicit in (the proof of)~\cite[Theorem 3.1]{Kul17}. Notice however that our proof is more natural, as it avoids reducing first $\iso$ to its restriction to linear orders with minimum and without maximum, and then the latter to the relations on arcs and knots, as it is done instead in \cite{Kul17}. 

\begin{proof}
In order to check that $F$ is a Borel function between the Polish space \( \LO \) and the standard Borel space \( \Ar \) one can argue as in \cite{Kul17}, so we only need to prove that \( F \) is a reduction.
 
Assume first that \( L,L' \in \LO \) are such that $L \cvx L'$, and let $g\colon L\to L'$ witness this. 
For every $n \in L$ the proper arcs $\arc{B_n^L}{f_n^L}$ and $\arc{B_{g(n)}^{L'}}{f_{g(n)}^{L'}}$ are both equivalent to $\arc{B^*}{f^*}$, and hence we can consider a homeomorphism $\varphi_1^n\colon \bar{B}_n^L \to \bar{B}_{g(n)}^{L'}$ witnessing this. Notice that \( \varphi_1^n\) is necessarily order-preserving
by~\eqref{eq:not-inverting}. Let $\varphi_1 = \bigcup_{n \in \N} \varphi_1^n$ and notice that $\varphi_1 (h_L(n),0,0) = (h_{L'}(g(n)),0,0)$ for every $n \in L$, so that the restriction of $\varphi_1$ to $I\Sigma_{F(L)}$ is order-preserving into $I\Sigma_{F(L')}$.

For every $n \in L$ let $M_n= \max \{2r_L(n), 2r_{L'}(g(n)) \}$ and let $\varphi_2^n \colon \overline{V_n^{L}} \times [-M_n,M_n]^2 \to \overline{V_{g(n)}^{L'}} \times [-M_n,M_n]^2$ be a homeomorphism which extends $\varphi_1^n$ and has the following properties:

\begin{enumerate}[label={\upshape (\roman*)}, leftmargin=2pc]
\item for all $(y,z) \in [-M_n,M_n]^2$, $\varphi_2^n(h_L(n) \pm 2r_L(n), y,z)=(h_{L'}(g(n)) \pm 2r_{L'}(g(n)), y,z)$;
\item for all $(y,z) \in [-M_n,M_n]^2$ with $\max \{|y|,|z|\}=M_n$ and for all $t \in [-1,1]$ we have
\[
\varphi_2^n(h_L(n)+ 2r_L(n)t, y,z)=(h_{L'}(g(n)) + 2r_{L'}(g(n))t, y,z)
\]
(this condition is missing in \cite{Kul17}).
\end{enumerate}

Let $W_n^{L}= \overline{V_n^{L}} \times [-1,1]^2$ and $W_{g(n)}^{L'} = \overline{V_{g(n)}^{L'}}\times [-1,1]^2$. We can then define a homeomorphism $\varphi_3^n\colon W_n^L\to W_{g(n)}^{L'}$ which extends $\varphi_2^n$ and is such that:
\begin{enumerate}[resume*]
\item\label{lines} for every $(y,z)\in [-1,1]^2$ such that $\max \{|y|, |z|\} \geq M_n$ we have
\[
\varphi_3^n (h_L(n)+ 2r_L(n)t,y,z) = (h_{L'}(g(n)) + 2r_{L'}(g(n))t,y,z),
\]
so that outside $\overline{V_n^{L}} \times [-M_n,M_n]^2$ the lines parallel to the $x$-axis are mapped into themselves.
\end{enumerate}
Then $\varphi_3 = \bigcup_{n \in \N} \varphi_3^n$ is a homeomorphism between $\bigcup_{n \in \N} W_n^L$ and $\bigcup_{n \in \N} W_{g(n)}^{L'}$.

We finally extend $\varphi_3$ to $\varphi \colon \bar{B} \to \bar{B}$ by looking at each $x_0 \in [-1,1] \setminus \bigcup_{n \in \N} \overline{V_n^L}$ (which is a cluster point of $\Ima h_L$) and setting $\varphi (x_0,y,z) = (x'_0,y,z)$ for every $(y,z) \in [-1,1]^2$, where $x_0= \lim_{i \to \infty} h_L(n_i)$ and $x'_0 = \lim_{i \to \infty} h_{L'}(g(n_i))$. Condition \ref{lines} ensures that $\varphi$ is continuous and indeed a homeomorphism.
It is immediate that $\varphi$ witnesses $F(L) \precsim_{\Ar} F(L')$, and that if \( g \colon L \to L' \) was actually an isomorphism, then \( \varphi \) witnesses
$F(L) \equiv_{\Ar} F(L')$.

Conversely, suppose that $\varphi \colon  \bar{B}\to \bar{B}$ is an embedding witnessing $F(L)\precsim_{\Ar} F(L')$.
Since all isolated points of \( F(L) \) belong to \( \Int \bar{B}\), by Lemma~\ref{lem:isolatedthroughembeddings} the map $\varphi \restriction I\Sigma_L$ embeds \( I\Sigma_{F(L)} \) into $I\Sigma_{F(L')}$. Furthermore, as explained in \cite{Kul17}, the embedding
\( \varphi \) preserves the betweenness relation. 
By \eqref{eq:not-inverting}, for any $n \in L$ the restriction of $\varphi$ to the arc
$\arc{B_n^L}{f_L \cap \bar{B}_n^L}$, which maps it to 
$\barc{\varphi(\bar{B}_n^L)}{f_{L'} \cap \varphi(\bar{B}_n^L)}$, is order-preserving and hence $\varphi \restriction I\Sigma_{F(L)}$ is order-preserving too. Moreover, since $\varphi$ is continuous and
\( f_L \) is connected we get that also \( \varphi(f_L) \) is
connected: it follows that $\varphi(I\Sigma_{F(L)})$ is a convex subset of
$I\Sigma_{F(L')}$. Summing up, $\varphi\restriction I\Sigma_{F(L)}$
witnesses that $I\Sigma_{F(L)}\trianglelefteq I\Sigma_{F(L')}$, hence
\(L \trianglelefteq_{\LO} L' \) because
$L\cong I\Sigma_{F(L)} \trianglelefteq I\Sigma_{F(L')}\cong L'$. Obviously, if \( \varphi \colon \bar{B} \to \bar{B} \) was actually a homeomorphism, then \( \varphi \restriction I\Sigma_{F(L)} \) would be onto \( I\Sigma_{F(L')} \), and thus it would witness $I\Sigma_{F(L)}\cong I\Sigma_{F(L')}$, which in turn implies \( L \iso L'\).
\end{proof}

By Theorem~\ref{thm:reduc_lo_arcs},
the reduction \( F \colon \LO \to \Ar \) 
 allows us to transfer some
combinatorial properties of ${\trianglelefteq_{\LO}}$
discussed in Section~\ref{sec:comb_prop} 
to the quasi-order $\precsim_{\Ar}$ (cfr.\ Lemma~\ref{open_int}, Proposition~\ref{prop1}, and Corollary~\ref{cor:intervalsleqBcvx}).

\begin{cor} \label{cor:sumuparcs}
\begin{enumerate-(a)}
\item \label{cor:sumuparcs-a}
There is an embedding from the partial order
\((\Cvx(\R),\subseteq)\) into $\precsim_{\Ar}$, and indeed \( {(\Cvx(\R),\subseteq)} \leq_B {\precsim_{\Ar}} \).
\item \label{cor:sumuparcs-b}
$\precsim_{\Ar}$ has chains of order type \( (\R, {<} ) \), as well as antichains of size $2^{\aleph_0}$.
\end{enumerate-(a)}
\end{cor}


In contrast, the combinatorial properties uncovered in Propositions~\ref{prop_WO}, \ref{prop:basisforcvx}, and \ref{prop:dom_fam}, being universal statements, do not transfer through the reduction \( F \). We overcome some of these difficulties by using the following construction.

Using the orientation induced by \( f \), when \( \arc{B}{f} \) is a proper arc the set \( I\Sigma_{\arc{B}{f}}\) can naturally be viewed as a linear order \( L_{\arc{B}{f}} = (I\Sigma_{\arc{B}{f}}, \leq_f) \). Since \( \Sigma_{\arc{B}{f}}\), being a subspace of the Polish space \( \bar{B} \), is second-countable, the set \( I\Sigma_{\arc{B}{f}}\) is (at most) countable and thus up to isomorphism \( L_{\arc{B}{f}} \) is an element of \( \Lin \). 
We remark that the linear order \( L_{\arc{B}{f}} \) really depends on the topological embedding \( f \) (or, more precisely, on the orientation it induces) rather than its image. However, if \( f \) and \( f' \) are two topological embeddings giving rise to the same arc, then either \( L_{\arc{B}{f}} = L_{\arc{B}{f'}} \) or \( L_{\arc{B}{f}} = (L_{\arc{B}{f'}})^* \) --- indeed the two linear orders correspond to the two possible orientations of the arc \( (\bar{B},\Ima f) \).
Recall that by construction, for proper arcs of the form%
\footnote{If not specified otherwise, we always choose the natural orientation of \( F(L) \).}
\( F(L) = \arc{B}{f_L} \)  we have
\( I\Sigma_{F(L)} \cong L \). 

\begin{lem} \label{lem:froarcstoLin}
  Let \( \arc{B}{f}, \arc{B'}{g} \in \Ar \)
  be such that \( \arc{B}{f}\precsim_{\Ar} \arc{B'}{g} \), and let 
  \( K = (I\Sigma_{\arc{B}{f}} \cap \Int \bar{B},{\leq_f}) \). Then either 
  \( K \trianglelefteq L_{\arc{B'}{g}}\) or \( K \trianglelefteq (L_{\arc{B'}{g}})^*\).
\end{lem}
\begin{proof} 
Let \( \varphi \colon \bar{B} \to \bar{B}' \) be an embedding witnessing \( \arc{B}{f}\precsim_{\Ar} \arc{B'}{g} \). By
Lemma~\ref{lem:isolatedthroughembeddings}\ref{lem:isolatedthroughembeddings-2}, 
\( \varphi \restriction (I\Sigma_{\arc{B}{f}} \cap \Int \bar{B})\) is an embedding of \( I\Sigma_{\arc{B}{f}} \cap \Int \bar{B}\) into \( I\Sigma_{\arc{B'}{g}}\), and arguing as in the proof of Theorem~\ref{thm:reduc_lo_arcs} we can observe that
\( \varphi(I\Sigma_{\arc{B}{f}} \cap \Int \bar{B})\) is a convex subset of \( I\Sigma_{\arc{B'}{g}}\) with respect to \( \leq_g \) because \( \Ima f \cap \Int \bar{B} \), which is homeomorphic to \( (0,1) \), is connected. 
As already noticed, the embedding \( \varphi\) is either order-preserving or order-reversing (Remark~\ref{rmk:orientation}\ref{rmk:orientation-2}): thus \( \varphi \restriction (I\Sigma_{\arc{B}{f}} \cap \Int \bar{B}) \) witnesses \( K \trianglelefteq L_{\arc{B'}{g}} \) in the former case, and \( K \trianglelefteq  (L_{\arc{B'}{g}})^*\) in the latter.
\end{proof}

\begin{remark}\label{rmk:froarcstoLin}
In the special case where \( \arc{B}{f}\) is of the form \( F(L)\) for some \( L \in \LO\), then all its isolated singular points belong to \(\Int \bar{B}\) and \( I\Sigma_{F(L)} \cong L\). Thus in this case Lemma~\ref{lem:froarcstoLin} reads as follows: For every \( L \in \LO \) and \( \arc{B'}{g} \in \Ar\) with \( F(L) \precsim_{\Ar} \arc{B'}{g} \), either \( L \trianglelefteq L_{\arc{B'}{g}}\) or \( L \trianglelefteq (L_{\arc{B'}{g}})^*\).
\end{remark}

Lemma~\ref{lem:froarcstoLin} allows us to prove analogues of Propositions~\ref{prop_WO} and \ref{prop:dom_fam}.

\begin{thm} \label{thm:umboundedanddominatingforarcs}
\( \mathfrak{b} (\precsim_{\Ar}) = \aleph_1\) and 
\( \mathfrak{d} (\precsim_{\Ar}) = 2^{\aleph_0}\).
\end{thm}

\begin{proof} 
We begin with the unbounding number. First notice that \( \mathfrak{b} (\precsim_{\Ar}) > \aleph_0 \) because given a countable family of proper arcs \( \{ \arc{B_i}{f_i} \mid i \in \N \} \), their infinite sum \( \bigoplus_\N \, \arc{B_i}{f_i} \) is a \( \precsim_{\Ar}\)-upper bound for them. To show the existence of an \( \precsim_{\Ar}\)-unbounded family of arcs of size \( \aleph_1 \) we use Proposition~\ref{prop_WO} as follows. Let \( F \colon \LO \to \Ar \) be the reduction introduced in~\eqref{eq:FfromLOtoA},
and consider the family \( \{ F(\boldsymbol{\alpha}) \mid \omega \leq \alpha < \omega_1 \} \). It is strictly \( \precsim_{\Ar}\)-increasing
by Theorem~\ref{thm:reduc_lo_arcs}.
Suppose towards a contradiction that there is \( \arc{B}{f} \in \Ar \) such that \( F(\boldsymbol{\alpha}) \precsim_{\Ar} \arc{B}{f}\) for all \( \alpha < \omega_1\). Then \( I\Sigma_{\arc{B}{f}} \) would be infinite and thus the linear order \( L = L_{\arc{B}{f}}\) would be, up to isomorphism, an element of \( \LO \). By Lemma~\ref{lem:froarcstoLin} and Remark~\ref{rmk:froarcstoLin}, this would lead to the fact that \( L+L^*\) is a \( \cvx\)-upper bound for \( \WO \), contradicting Proposition~\ref{prop_WO}.

We now deal with the dominating number.
Consider once again the \( \cvx \)-antichain \( \mathcal{A} =  \{ L_S \mid S \subseteq \N \} \), where \( L_S = \eta_{f_S}\) is as in the proof of Proposition~\ref{prop:basisforcvx}\ref{prop:basisforcvx-1}, and notice that by the usual back-and-forth argument \( L_S \cong (L_S)^*\). We first prove the analogue of Claim~\ref{claim:everylobelongstoantichain}.

\begin{claim} \label{claim:everylobelongstoantichainforarcs}
For every proper arc \( \arc{B}{f} \in \Ar\), the collection
\[
\{ F(L_S) \mid  {F(L_S) \precsim_{\Ar} \arc{B}{f}} \}
\] 
is countable.
\end{claim}

\begin{proof}[Proof of the Claim]
By Lemma~\ref{lem:froarcstoLin}, Remark~\ref{rmk:froarcstoLin}, and \( L_S \cong (L_S)^*\) we have that 
\[ 
\{ L_S \in \mathcal{A} \mid F(L_S) \precsim_{\Ar} \arc{B}{f} \} \subseteq  \{ L_S \in \mathcal{A} \mid {L_S \trianglelefteq L_{\arc{B}{f}}} \} . 
\]
If \( L_{\arc{B}{f}}\) is finite, then the latter set is empty and so is the set in the claim; if instead  \( L_{\arc{B}{f}} \) is infinite then, up to isomorphism, it is a member of \( \LO \), and thus the result easily follows from 
Claim~\ref{claim:everylobelongstoantichain}.
\end{proof}

The proof of the theorem can now be completed using the same argument as the one used in the proof of 
Proposition~\ref{prop:dom_fam}: every element of a \( \precsim_{\Ar}\)-dominating family has only countably many proper arcs of the form \( F(L_S)\) below it, and since by Theorem~\ref{thm:reduc_lo_arcs} there are \( 2^{\aleph_0}\)-many such arcs the dominating family must have size \( 2^{\aleph_0} \) too.
\end{proof}

As in the case of linear orders, one can then derive the following analogue of 
Corollary~\ref{cor_no_max}.
However, the proof is slightly more delicate.

\begin{cor} \label{cor:incomparablearcsandunboundedchains}
Every proper arc \( \arc{B}{f} \) is the bottom of an \( \precsim_{\Ar}\)-unbounded chain of length~\(\omega_1\).
\end{cor}

\begin{proof}
%
Consider the sequence of proper arcs \( (\arc{B_\alpha}{f_\alpha}))_{\alpha < \omega_1}\) where \( \arc{B_0}{f_0} = \arc{B}{f}\) and \( \arc{B_\alpha}{f_\alpha} = \arc{B}{f} \oplus F(\boldsymbol{\omega+\alpha}) \) if \(  \alpha \geq 1 \).
For every \( \alpha \leq \beta < \omega_1\) we have \( \arc{B_\alpha}{f_\alpha} \precsim_{\Ar} \arc{B_\beta}{f_\beta}\), and if \( \alpha > 0\) we also have \( F(\boldsymbol{\omega+\alpha}) \precsim_{\Ar} \arc{B_\alpha}{f_\alpha} \). By (the proof of) Theorem~\ref{thm:umboundedanddominatingforarcs}, this implies that the sequence is \( \precsim_{\Ar}\)-unbounded. Moreover, for every \( \alpha < \omega_1 \) there is \( \beta > \alpha \) such that \( \arc{B_\beta}{f_\beta} \not\precsim_{\Ar} \arc{B_\alpha}{f_\alpha}\) (and hence \( \arc{B_\alpha}{f_\alpha} \prec_{\Ar} \arc{B_\beta}{f_\beta} \)), as otherwise \( \arc{B_\alpha}{f_\alpha}\) would be an upper bound for the sequence \( (\arc{B_\alpha}{f_\alpha})_{\alpha < \omega_1}\). It follows that we can extract from the latter a strictly \( \precsim_{\Ar}\)-increasing subsequence of length \( \omega_1 \) with \( \arc{B}{f} \) as a first element: since such a subsequence is \( \precsim_{\Ar} \)-cofinal in \( (\arc{B_\alpha}{f_\alpha})_{\alpha < \omega_1}\), it is \( \precsim_{\Ar} \)-unbounded as well and the proof is complete.
\end{proof}

We now move to the possible generalizations of Proposition~\ref{prop:basisforcvx}, i.e.\ we discuss minimal elements and bases for the relation \( \precsim_{\Ar} \).

If we consider only tame proper arcs, which form a \( \precsim_{\Ar} \)-downward closed subclass of the collection of all proper arcs, then the situation is pretty clear: the trivial arc \( I_{\Ar} \) is the \( \precsim_{\Ar}\)-minimum within this class. Call \textbf{prime arc} any proper arc of the form \( \arc{B_K}{f_K} \) for \( K \) a prime knot: 
then one can observe that prime arcs are \( \precsim_{\Ar} \)-minimal above \( I_{\Ar} \) (and are the unique such).
Indeed, assume that \( \arc{B}{f} \) is a prime arc and that \( \arc{B'}{g} \precsim_{\Ar} \arc{B}{f} \) for some \( \arc{B'}{g} \in \Ar \).
Since \( \arc{B}{f} \) is tame, without loss of generality we can assume that \(  \bar{B}' \subseteq \Int \bar{B}\) and that both $\bar B$ and $\bar B'$ are smoothly embedded as well as that the witness $\varphi$
for \( \arc{B'}{g} \precsim_{\Ar} \arc{B}{f} \) is smooth. Let \( K_1  = K_{\arc{B'}{g}} \), and let \( K_2 \) be the knot obtained from the remainder \( f \setminus \Int \varphi[\bar{B}']\) by connecting  \( g(0) \) and \( g(1) \) with a simple curve lying on \( \partial \varphi[\bar B'] \) and the extrema \( f(0) \) and \( f(1) \) with a simple curve on \( \partial \bar{B} \). By construction, the prime knot used to construct \( \arc{B}{f}\) is  the sum of \( K_1\) and \( K_2\), thus one of \( K_1 \) and \( K_2\) is trivial. In the former case \( \arc{B'}{g} \equiv_{\Ar}  I_{\Ar} \), while in the latter \( \arc{B'}{g} \equiv_{\Ar} \arc{B}{f}\).

Prime arcs play the same role in the realm of tame proper arcs as prime knots do in the realm of tame knots: every tame proper arc is of the form \( \bigoplus_{i \leq n} \, \arc{B_i^p}{f_i^p} \) for some (unique, up to permutations) sequence of prime arcs \( \arc{B_i^p}{f_i^p}\).
This has a number of consequences on the structure of \emph{nontrivial} tame proper arcs under \( \precsim_{\Ar}\): 
\begin{enumerate-(1)}
\item
There are no infinite descending chains.
\item
Since up to equivalence there are only countably many tame proper arcs, and since there are infinitely many prime arcs (consider e.g.\ the prime arcs obtained from the \( (p,q) \) torus knots, where \( p,q > 1\)), it follows that the collection of prime arcs constitutes a countably infinite antichain basis. In particular, there are no finite bases.
\item
If \( \arc{B_i^p}{f_i^p} \), for \( i \in \N\), is an enumeration of the prime arcs, then \( (\bigoplus_{i \leq n} \, \arc{B_i^p}{f_i^p})_{n \in \N} \) is an unbounded \( \omega \)-chain. In particular, there is no \( \precsim_{\Ar}\)-maximal tame proper arc, and the unbounding number of \( \precsim_{\Ar} \) restricted to tame proper arcs is \( \aleph_0 \).
\item
Every dominating family is infinite: below every tame proper arc there are only finitely many of the infinitely many pairwise \( \precsim_{\Ar}\)-incomparable prime arcs, thus no finite family can be dominating with respect to \( \precsim_{\Ar} \). Hence the dominating number of \( \precsim_{\Ar} \) restricted to tame proper arcs is \( \aleph_0\).
\end{enumerate-(1)}

Having obtained the desired information in the realm of tame proper arcs, 
it is now natural to move to the wild side and consider the restriction \( \precsim_{\WAr}\) of  \( \precsim_{\Ar} \) to the collection \( \WAr \) of wild arcs. By \( \cvx \)-minimality of \( \eta_{f_S} \), Lemma~\ref{lem:froarcstoLin} and Remark~\ref{rmk:froarcstoLin}, one may be tempted to conjecture that the proper arcs \( F(\eta_{f_S}) \) used in the proof of Theorem~\ref{thm:umboundedanddominatingforarcs} are \( \precsim_{\WAr}\)-minimal. That is not quite true, as the arc 
\( \arc{B^*}{f^*} = \bigoplus_\N \ \arc{B_i}{f_i}\) used to define the reduction \( F = F_{\arc{B^*}{f^*}} \colon \LO \to \Ar \) from~\eqref{eq:FfromLOtoA} is such that \( \arc{B^*}{f^*} \prec_{\WAr} F(\eta_{f_S})\),
and moreover the proper arc \( \arc{B^\partial}{f^\partial} = \bigoplus^\partial_\N \, \arc{B_i}{f_i}\) is such that
\( \arc{B^\partial}{f^\partial} \prec_{\WAr} \arc{B^*}{f^*} \prec_{\WAr} F(\eta_{f_S}) \).
However, the following lemma allows us to obtain useful information on the \( \precsim_{\WAr} \)-predecessors of \( F(\eta_{f_S})\). 

\begin{lem} \label{lem:fewpredecessors} 
Let \( \{ \arc{B_i^p}{f_i^p} \mid i \in \N \} \) be a family of (oriented) prime arcs,
and let \( \arc{B^*}{f^*} = \bigoplus_{i \in \N }\, \arc{B_i^p}{f_i^p} \) and \( \arc{B^\partial}{f^\partial} = \bigoplus^\partial_{i \in \N} \, \arc{B_i^p}{f_i^p} \).
We consider the proper arc \( F_{\arc{B^*}{f^*}}(L_S)\) for some \( S \subseteq \N\), where \( L_S = \eta_{f_S}\) is as in the proof of Proposition~\ref{prop:basisforcvx}\ref{prop:basisforcvx-1}.
\begin{enumerate}[label={\upshape (\alph*)}, leftmargin=2pc]
\item \label{lem:fewpredecessors-1}
If \( \arc{B'}{g}\) is a prime arc, then 
\( \arc{B'}{g} \precsim_{\Ar} 
\arc{B^\partial}{f^\partial} \)
if and only if there is \( \bar\imath \in \N \) such that \( \arc{B'}{g} \equiv_{\Ar} \arc{B_{\bar\imath}^p}{f_{\bar\imath}^p} \). 
The same is true if \( \arc{B^\partial}{f^\partial} \) is replaced by 
\( F_{\arc{B^*}{f^*}}(L_S) \).
\end{enumerate}
Let now \( \arc{B}{f} \) be an arbitrary wild proper arc. Then:
\begin{enumerate}[resume*]
\item \label{lem:fewpredecessors-2}
\( \arc{B}{f} \precsim_{\WAr} \arc{B^\partial}{f^\partial} \) if and only if \( \arc{B}{f} \equiv_{\Ar} \bigoplus_{j \in A}^\partial \, \arc{B_j^p}{f_j^p} \) for some infinite set \( A \subseteq \N \).


\item \label{lem:fewpredecessors-4}
If \( \arc{B}{f} \precsim_{\WAr} F_{\arc{B^*}{f^*}}(L_S)\), then there is
\( \bar B' \subseteq \bar B\) such that \( \arc{B'}{f \cap \bar B'} \in \WAr \)
and \( \arc{B'}{f \cap \bar B'} \precsim_{\WAr} \arc{B^\partial}{f^\partial} \). 
\end{enumerate}
\end{lem}

\begin{proof}
\ref{lem:fewpredecessors-1} 
One direction is obvious. For the other direction, assume that
\( \arc{B'}{g} \precsim_{\Ar} \arc{B^\partial}{f^\partial} \).
Recall that the ambient space \( \bar B^\partial\) of \(  \bigoplus_{i \in \N}^\partial \, \arc{B_i^p}{f_i^p} \) is the cube \( [-1,1]^3\), and that its only singularity is the point \( (1,0,0) \). By the way we defined infinite sums, without loss of generality we may assume that \( \bar{B}_i^p = [1-2^{-i},1-2^{-(i+1)}] \times [-1,1]^2\) and that 
\( f_i^p(0) = (1-2^{-i},0,0) \), and \( f_i^p(1) = (1-2^{-(i+1)},0,0) \) for all \( i \in \N \). 
Let \( \varphi \colon \bar{B}' \to [-1,1]^3 \) witness \( \arc{B'}{g} \precsim_{\Ar} \arc{B^\partial}{f^\partial} \), and notice that 
since \( \arc{B_0^p}{f_0^p} \) is tame we may assume \( \Ima \varphi \subseteq [0,1] \times [-1,1]^2\).
Let \( I = \{ i \in \N \mid \varphi(\bar{B}') \cap f_i^p \cap \Int \bar{B}_i^p \neq \emptyset \} \): it is convex (with respect to the usual order \( \leq \) on \( \N \)) because \( g \) is connected.
Without loss of generality, we may assume 
that  for all \( i \in I\) the space
\( \bar B'_i =  \varphi(\bar{B}') \cap \bar{B}_i^p \) is a 
closed topological 3-ball. 
Consider the (tame) proper arcs \( \arc{B'_i}{f_i \cap \bar{B}'_i} \), which by construction are such that either \( \arc{B'}{g} \equiv_{\Ar} \bigoplus_{i \in I} \, \arc{B'_i}{f_i \cap \bar B'_i} \) if \( I \) is finite, or 
\( \arc{B'}{g} \equiv_{\Ar} \bigoplus^\partial_{i \in I} \, \arc{B'_i}{f_i \cap \bar B'_i} \) if \( I \) is infinite.
Moreover each \( \arc{B'_i}{f_i \cap \bar{B}'_i} \) is either trivial or equivalent to the corresponding \( \arc{B_i^p}{f_i^p} \) because \( \arc{B'_i}{f_i \cap \bar{B}'_i} \precsim \arc{B_i^p}{f_i^p} \) and the latter is prime.
If all the \( \arc{B'_i}{f_i \cap \bar{B}'_i} \)'s were trivial, then \( \arc{B'}{g} \) would be trivial too, a contradiction. 
Let \( \bar \imath \in I \) be such that \( \arc{B'_{\bar\imath}}{f_{\bar\imath} \cap \bar{B}'_{\bar\imath}} \) is not trivial: then \( \arc{B'_{\bar\imath}}{f_{\bar\imath} \cap \bar{B}'_{\bar\imath}} \equiv_{\Ar} \arc{B_{\bar\imath}^p}{f_{\bar\imath}^p} \), and since $\arc{B'}{g}$ is prime and \( \varphi^{-1} \restriction \bar B'_{\bar\imath} \) witnesses \( \arc{B'_{\bar\imath}}{f_{\bar\imath} \cap \bar{B}'_{\bar\imath}} \precsim_{\Ar} \arc{B'}{g} \) it follows that \( \arc{B_{\bar\imath}^p}{f_{\bar\imath}^p} \equiv_{\Ar} \arc{B'}{g} \), as desired.

Suppose now that \( \arc{B'}{g} \precsim_{Ar} F_{\arc{B^*}{f^*}}(L_S) \) via some \( \varphi \).
Recall the notation used in the proof of Theorem~\ref{thm:reduc_lo_arcs} and in the discussion preceding it.
Since \( \arc{B'}{g}\) is tame and not trivial,  \( \varphi (g) \) is tame and cannot be contained in \( [h_{L_S}(m), h_{L_S}(m) +2r_{L_S}(m)] \times [-1,1]^2\) for any $m \in L_S$; therefore \( \varphi (g) \) must be contained either in \( [h_{L_S}(n)-2 r_{L_S}(n), h_{L_S}(n) - \varepsilon] \times [-1,1]^2\)
or in \( [h_{L_S}(m), h_{L_S}(n) - \varepsilon] \times [-1,1]^2\) for some consecutive \( m,n \in L_S\) and small enough \( \varepsilon > 0  \). However, since the part on the right of the singularity \( h_{L_S}(m) \) is trivial and \( \arc{B_0^p}{f_0^p} \) is tame, we can actually assume that we are always in the first case and that \( \Ima \varphi \subseteq [h_{L_S}(n)-2 r_{L_S}(n), h_{L_S}(n)] \times [-1,1]^2 \). Since the subarc of \( F_{\arc{B^*}{f^*}}(L_S) \) determined by the latter set is equivalent to \( \bigoplus_{i \in \N}^{\partial} \, \arc{B_i^p}{f_i^p} \), we are done by the first part.

\ref{lem:fewpredecessors-2}
Let \( \varphi \colon \bar{B} \to [-1,1]^3\) witness \( \arc{B}{f} \precsim_{\WAr} \arc{B^\partial}{f^\partial} \). 
Being wild and \( \precsim_{\WAr} \)-below an arc with only one singularity, the proper arc \( \arc{B}{f} \) has a unique singularity \(  x \in \bar B \): clearly,  \( \varphi(x) = (1,0,0) \) by Lemma~\ref{lem:isolatedthroughembeddings}\ref{lem:isolatedthroughembeddings-1} and thus, necessarily,  \( x \in \partial \bar{B} \).
As before, set \( I = \{ i \in \N \mid \varphi(\bar{B}) \cap f_i^p \cap \Int \bar{B}_i^p \neq \emptyset \} \): now, we know that \( I \) is a final segment of \( (\N, {\leq})  \) because \( \varphi(x) = (1,0,0) \). 
We can assume that \( \Ima \varphi \setminus (\{ 1 \} \times [-1,1]^2) \subseteq \bigcup_{i \in I} \bar{B}_i^p \) and that  for all \( i \in I\) the space
\( \bar B'_i =  \varphi(\bar{B}) \cap \bar{B}_i^p \) is a 
closed topological 3-ball. Then \( \varphi \) witnesses \( \arc{B}{f} \equiv_{\Ar} \bigoplus_{i \in I}^\partial \, \arc{B'_i}{f_i^p \cap \bar{B}'_i} \).
Each of the proper arcs \( \arc{B'_i}{f_i^p \cap \bar{B}'_i} \), being a subarc of the prime arc \( \arc{B_i^p}{f_i^p} \), is either trivial or equivalent to it: set
\[  
A = \{ i \in I \mid \arc{B'_i}{f_i^p \cap \bar{B}'_i} \equiv_{\Ar} \arc{B_i^p}{f_i^p} \}.
\]
The set \( A \) is infinite because otherwise \( \bigoplus_{i \in I}^\partial \, \arc{B'_i}{f_i^p \cap \bar{B}'_i} \), and hence also \( \arc{B}{f} \), would be tame. 
Moreover, each of the \( \arc{B'_i}{f_i^p \cap \bar{B}'_i} \) is tame, so the trivial arcs \( \arc{B'_i}{f_i^p \cap \bar B'_i}\) occurring in the sequence can be ``absorbed'' by the next \( \arc{B'_j}{f_j^p \cap \bar{B}'_j}\) with \( j \in A \). Therefore 
\[ 
\arc{B}{f} \equiv_{\Ar}
\bigoplus_{i \in I}^\partial \, \arc{B'_i}{f_i^p \cap \bar{B}'_i} \equiv_{\Ar} 
\bigoplus_{j \in A}^\partial \, \arc{B'_j}{f_j^p \cap \bar{B}'_j} \equiv_{\Ar}
\bigoplus_{j \in A}^\partial \, \arc{B_j^p}{f_j^p}.
\]

Conversely,
assume that \( \arc{B}{f} \equiv_{\Ar} \bigoplus^\partial_{j \in A} \, \arc{B_j^p}{f_j^p} \) for some infinite \( A \subseteq \N\). For each \( i \in \N\), set \( \arc{B'_i}{f'_i} = \arc{B_i^p}{f_i^p} \) if \( i \in A\) and \( \arc{B'_i}{f'_i} = I_{\Ar} \) if \( i \notin A \): since all the proper arcs \( \arc{B_i^p}{f_i^p}\) are tame, we get \( \bigoplus^\partial_{j \in A} \, \arc{B_j^p}{f_j^p} \equiv_{\Ar} \bigoplus^\partial_{i \in \N} \, \arc{B'_i}{f'_i} \), so it is enough to show that the latter is a subarc of \( \bigoplus^\partial_{i \in \N} \, \arc{B_i^p}{f_i^p} \).  Without loss of generality, \( \bar B_i = \bar B'_i = [1-2^{-i},1-2^{-(i+1)}] \times [-2^{-i},2^{-i}]^2 \), \( f_i(0) = f'_i(0) = (1-2^{-i},0,0) \), and \( f_i(1) = f'_i(1) = (1-2^{-(i+1)},0,0) \). We can further assume that the ambient space of \( \bigoplus^\partial_{i \in \N} \, \arc{B'_i}{f'_i} \) is the ``step pyramid'' \( ([-1,0] \times [-1,1]^2) \cup {\bigcup_{i \in \N} \bar{B}'_i} \cup \{ (1,0,0) \} \).
For each \( i \notin A\), fix a tubular neighborhood \( \bar{B}''_i \subseteq \bar{B}_i^p\) of \( f_i^p\), i.e.\ a ``cylinder'' of radius \( \varepsilon_i \) with rotation axis given by \( f_i^p \) itself --- this is possible because \( \arc{B_i^p}{f_i^p} \) is tame. Moreover, since each block of consecutive \( i \in \N \setminus A\) is finite, we can assume that \( \varepsilon_i = \varepsilon_{i+1} \) if \( i, i+1 \notin A \). 
For \( i \in A \) pick instead \( \bar{B}''_i \subseteq \bar{B}_i^p \) so that: \( \arc{B''_i}{f_i^p}\) is a proper arc; \( \bar{B}''_i \) intersects the left face \( \{ 1-2^{-i} \} \times [-2^{-i},2^{-i}]^2 \) of \( \bar B_i^p \) in a disc of radius \( \varepsilon_i\) centered in \( (1-2^{-i},0,0) \), where \( \varepsilon_i = \varepsilon_{i-1} \) if \( i > 0 \) and \( i -1 \notin A \) and \( \varepsilon_i = 2^{-i} \) otherwise; similarly, \( \bar{B}''_i \) intersects the right face \( \{ 1-2^{-(i+1)} \} \times [-2^{-i},2^{-i}]^2 \) of \( \bar B_i \) in a disc of radius \( \varepsilon_i\)  centered in \( (1-2^{-(i+1)},0,0) \), where \( \varepsilon_i = \varepsilon_{i+1} \) if \( i +1 \notin A \) and \( \varepsilon_i = 2^{-(i+1)} \) otherwise. 
Finally, let \( \bar B''_{-1} \subseteq [-1,0] \times [-1,1]^2\) be such that \( \arc{B''_{-1}}{[-1,0] \times \{ (0,0) \}}\) is a proper arc and \( \bar B''_{-1}\) intersects the left face \( \{ 0 \} \times [-1,1]^2\) of \( \bar B_0^p \) in a disc of radius \( 1 \) centered in the origin \( (0,0,0) \). By construction, \( \bar B''_{-1} \cup \bigcup_{i \in \N} \bar B''_i \cup \{ (1,0,0) \} \) is homeomorphic to a (closed) cone, and thus it is a closed topological 3-ball. Moreover, every \( \arc{B'_i}{f'_i} \) is equivalent to \( \arc{B''_i}{f_i^p}\) via some \( \varphi_i \colon \bar B'_i \to \bar B''_i \). Fix also a homeomorphism \( \varphi_{-1} \colon [-1,0] \times [-1,1]^2 \to \bar B''_{-1} \) fixing the interval \( [-1,0] \times \{ (0,0)\}\), and let \( \varphi_\infty \) be the identity on the singleton \( (1,0,0) \). Then \( \varphi = \varphi_{-1} \cup {\bigcup_{i  \in \N} \varphi_i} \cup \varphi_\infty \) is an embedding witnessing \( \bigoplus^\partial_{i \in \N} \, \arc{B'_i}{f'_i} \precsim_{\Ar} \bigoplus^\partial_{i \in \N} \, \arc{B_i^p}{f_i^p} \), as desired.

\ref{lem:fewpredecessors-4}
Let \( \varphi \colon \bar B \to [-1,1]^3\) witness \( \arc{B}{f} \precsim_{\Ar} F_{\arc{B^*}{f^*}}(L_S)\).
We claim that there is \( x \in \Sigma_{\arc{B}{f}} \) such that \( \varphi(x)  \in I\Sigma_{F_{\arc{B^*}{f^*}}(L_S)}\).
Pick any \( x \in \Sigma_{\arc{B}{f}}\), so that \( \varphi(x) \in \Sigma_{F_{\arc{B^*}{f^*}}(L_S)} \) as well by Lemma~\ref{lem:isolatedthroughembeddings}\ref{lem:isolatedthroughembeddings-1}. If \( \varphi(x) \notin I\Sigma_{F_{\arc{B^*}{f^*}}(L_S)} \), we use the fact that
by construction every singularity \( y \in \Sigma_{F_{\arc{B^*}{f^*}}(L_S)} \setminus I\Sigma_{F_{\arc{B^*}{f^*}}(L_S)} \) is a limit of isolated singularities \emph{from both sides} (unless \( y \in \partial \bar B^* \), in which case there is only one side available), and hence for all \(  \bar B' \subseteq [-1,1]^3 \) with \( y \in \bar B' \) and \( \arc{B'}{f_{L_S} \cap \bar B'} \in \Ar \) the set \( I\Sigma_{F_{\arc{B^*}{f^*}}(L_S)} \cap \bar B' \) is infinite. 
Applying this to \( \bar B' = \varphi(\bar B) \) and \( y = \varphi(x) \), we get that there is some (in fact, infinitely many) \( y' \in I\Sigma_{F_{\arc{B^*}{f^*}}(L_S)} \cap \Int \bar B'\): by Lemma~\ref{lem:isolatedthroughembeddings}, replacing \( x \) with \( \varphi^{-1}(y') \) we are done.
Using the same notation as in the proof of Theorem~\ref{thm:reduc_lo_arcs}, let \( n \in L_S \) be such that \( \varphi(x) = (h_{L_S}(n),0,0) \).
Without loss of generality, we may assume that \( \bar B'' =  \Ima \varphi \cap [h_{L_S}(n)-2 r_{L_S}(n), h_{L_S}(n)] \times [-1,1]^2 \) is a closed topological  3-ball. Moreover, \( \varphi(x) =  (h_{L_S}(n),0,0) \in \Sigma_{\arc{B''}{f_{L_S} \cap \bar B''}} \) because otherwise \( x \) would not be a singularity of \( \arc{B}{f} \) (here we use the fact that on the right of \( (h_{L_S}(n),0,0) \) there is a trivial arc),
thus \( \arc{B''}{f_{L_S} \cap \bar B''} \in \WAr\).
Since the subarc of \( F_{\arc{B^*}{f^*}}(L_S) \) determined by \( [h_{L_S}(n)-2 r_{L_S}(n), h_{L_S}(n)] \times [-1,1]^2 \) is equivalent to \( \arc{B^\partial}{f^\partial} \), setting \( B' = \varphi^{-1}(B'') \) we are done.
\end{proof}

We are not able to get a full analogue of Proposition~\ref{prop:basisforcvx},
but Lemma~\ref{lem:fewpredecessors} 
allows us to
get a similar, although slightly weaker, result.

\begin{thm} \label{thm:basisforarcs}
\begin{enumerate-(a)} 
\item \label{thm:basisforarcs-1}
There are infinitely many \( \precsim_{\WAr}\)-incomparable \( \precsim_{\WAr} \)-minimal elements in \( \WAr \).
\item \label{thm:basisforarcs-2}
There is a strictly \( \precsim_{\WAr} \)-decreasing \( \omega \)-sequence in \( \WAr \) which is not \( \precsim_{\WAr}\)-bounded from below.
\item \label{thm:basisforarcs-3}
No basis for \( \precsim_{\WAr} \) has size smaller than \( 2^{\aleph_0} \).
\end{enumerate-(a)}
\end{thm}

\begin{proof}
Fix an enumeration without repetitions \( \{ \arc{B_i^p}{f_i^p} \mid i \in \N \} \) of all prime arcs.

\ref{thm:basisforarcs-1}
For each \( k \in \N \) set \( \arc{B'_k}{g_k} = \bigoplus^\partial_\N \, \arc{B_k^p}{f_k^p} \). Every \( \arc{B'_k}{g_k} \) is \( \precsim_{\WAr} \)-minimal by Lemma~\ref{lem:fewpredecessors}\ref{lem:fewpredecessors-2} (and the fact that all arcs in the infinitary sum are the same), and if \( k \neq k' \) then \( \arc{B'_k}{g_k} \not\precsim_{\WAr} \arc{B'_{k'}}{g_{k'}} \) because \( \arc{B_k^p}{f_k^p} \precsim_{\Ar} \arc{B'_k}{g_k} \) but \( \arc{B_k^p}{f_k^p} \not\precsim_{\Ar} \arc{B'_{k'}}{g_{k'}} \) by Lemma~\ref{lem:fewpredecessors}\ref{lem:fewpredecessors-1}.

\ref{thm:basisforarcs-2}
Now let \( \arc{B'_k}{g_k} = \bigoplus_{i \geq k}^\partial \, \arc{B_i^p}{f_i^p} \). 
By parts~\ref{lem:fewpredecessors-1} and~\ref{lem:fewpredecessors-2} of Lemma~\ref{lem:fewpredecessors}, if \( k < k' \) then \( \arc{B'_{k'}}{g_{k'}} \prec_{\Ar} \arc{B'_{k}}{g_{k}} \). Moreover, by Lemma~\ref{lem:fewpredecessors}\ref{lem:fewpredecessors-2} if \( \arc{B}{f} \in \WAr\) is such that \( \arc{B}{f} \precsim_{\WAr} \arc{B'_0}{g_0} \) then \( \arc{B}{f} \equiv_{\Ar} \bigoplus_{j \in A}^\partial \, \arc{B_j^p}{f_j^p} \), for some infinite $A \subseteq \N$. Let \( k = \min A \): since $\arc{B_k^p}{f_k^p} \precsim_{\WAr} \bigoplus_{j \in A}^\partial \, \arc{B_j^p}{f_j^p}$ but by Lemma~\ref{lem:fewpredecessors}\ref{lem:fewpredecessors-1} $\arc{B_k^p}{f_k^p} \not\precsim_{\WAr} \arc{B'_{k+1}}{g_{k+1}}$, we have \( \arc{B}{f} \not\precsim_{\WAr} \arc{B'_{k+1}}{g_{k+1}}\). Thus the chain formed by the proper arcs \( \arc{B'_k}{g_k} \) is as required.

\ref{thm:basisforarcs-3}
Let \( \{ A_x \mid x \in 2^{\N} \} \) be a family of infinite sets \( A_x \subseteq \N \) such that \( A_x \cap A_y \) is finite for all distinct \( x,y \in 2^{\N}  \). (For the sake of definiteness, set \( A_x = \{ h(x \restriction n) \mid n \in \N \} \), where \( h \) is a bijection from all finite binary sequences to the natural numbers.) 
Fix a basis \( \mathcal{B} \) for \( \precsim_{\WAr} \).
Then for every \( \bigoplus_{j \in A_x}^\partial \, \arc{B_j^p}{f_j^p} \) there is some \( \arc{B_x}{f_x} \in \mathcal{B}\) such that \( \arc{B_x}{f_x} \precsim_{\WAr} \bigoplus_{j \in A_x}^\partial \, \arc{B_j^p}{f_j^p} \). By Lemma~\ref{lem:fewpredecessors}\ref{lem:fewpredecessors-2}, \( \arc{B_x}{f_x} \equiv_{\Ar} \bigoplus_{j \in A'_x}^\partial \, \arc{B_j^p}{f_j^p} \) for some infinite \( A'_x \subseteq A_x \). If there were distinct \( x,y \in 2^{\N} \) such that \( \arc{B_x}{f_x} \equiv_{\Ar} \arc{B_y}{f_y} \), then we would get \( A'_x= A'_y\) by Lemma~\ref{lem:fewpredecessors}\ref{lem:fewpredecessors-1}, and thus \( A'_x \subseteq A_x \cap A_y \), which is impossible because \( A'_x \) is infinite. Thus all the proper arcs \( \arc{B_x}{f_x} \in \mathcal{B} \) are distinct, and thus \( |\mathcal{B}| \geq 2^{\aleph_0}\), as desired.
\end{proof}

It is plausible that part \ref{thm:basisforarcs-1} of the previous theorem can be improved by showing the existence of $2^{\aleph_0}$-many \( \precsim_{\WAr}\)-incomparable \( \precsim_{\WAr} \)-minimal elements in \( \WAr \). 
To this end, one could consider proper arcs that are singular everywhere and hence are very different from the arcs constructed in this paper. 

Lemma~\ref{lem:fewpredecessors} is also sufficient to recover an analogue of Proposition~\ref{prop:everylobelongstoantichain} for proper arcs.

\begin{thm} \label{thm:antichainsforarcs} 
Every \( \precsim_{\WAr}\)-antichain is contained in a \( \precsim_{\WAr}\)-antichain of size \( 2^{\aleph_0} \). In particular, there are no maximal \( \precsim_{\WAr}\)-antichains of size smaller than \( 2^{\aleph_0}\), and every \( \arc{B}{f} \in \WAr \) belongs to a \( \precsim_{\WAr}\)-antichain of size \( 2^{\aleph_0} \).
\end{thm}

\begin{proof}
Let \( \mathcal{A} =  \{ \arc{B'_m}{g_m} \mid m < \kappa \} \) be a \( \precsim_{\WAr} \)-antichain, where \( \kappa < 2^{\aleph_0} \).
Let \( \arc{B_i^p}{f_i^p} \), for \( i \in \N \), be an enumeration without repetitions of all prime arcs, and for \( S \subseteq \N \) let \( L_S = \eta_{f_S}\) be the linear order from the proof of Proposition~\ref{prop:basisforcvx}\ref{prop:basisforcvx-1}. 
Let \( \{ A_S \mid S \subseteq \N \} \) be a family of sets \( A_S \subseteq \N \) such that \( A_S \cap A_{S'} \) is finite for all distinct \( S,S' \subseteq \N \). 
(Such a family can be constructed as in the proof of Theorem~\ref{thm:basisforarcs}\ref{thm:basisforarcs-3}.)
For each \( S \subseteq \N \), set \( \arc{B^*_S}{f^*_S} = \bigoplus_{j \in A_S} \, \arc{B_j^p}{f_j^p} \) and
\[  
\arc{B_S}{f_S} = F_{\arc{B^*_S}{f^*_S}}(L_S).
\]
Let \( \mathcal{B}\) be the collection of all proper arcs of the form \( \arc{B_S}{f_S} \) which are \( \precsim_{\WAr} \)-incomparable with every \( \arc{B'_m}{g_m} \in \mathcal{A} \).

\begin{claim}\label{claim:antichain_for_arcs}
\( |\mathcal{B}| = 2^{\aleph_0} \).
\end{claim}

\begin{proof}[Proof of the Claim]
By the proof of Claim~\ref{claim:everylobelongstoantichainforarcs} and \( \kappa < 2^{\aleph_0} \) there are \( 2^{\aleph_0} \)-many proper arcs \( \arc{B_S}{f_S}\) such that \( \arc{B_S}{f_S} \not\precsim_{\WAr} \arc{B'_m}{g_m} \) for all \( m < \kappa \). 
On the other hand, we claim that there are at most \( \kappa \)-many proper arcs \( \arc{B_S}{f_S} \) such that \( \arc{B'_m}{g_m} \precsim_{\WAr} \arc{B_S}{f_S} \) for some \( m < \kappa \), which suffices to prove the claim. Indeed, suppose that \( m < \kappa \) and \( S \subseteq \N \) are such that \( \arc{B'_m}{g_m} \precsim_{\WAr} \arc{B_S}{f_S} \).
Then by parts~\ref{lem:fewpredecessors-4} and~\ref{lem:fewpredecessors-2} of Lemma~\ref{lem:fewpredecessors} there is \( \bar B'' \subseteq \bar B'_m \) such that \( \arc{B''}{g_m \cap \bar B''} \equiv_{\Ar} \bigoplus_{j \in A}^\partial \, \arc{B_j^p}{f_j^p} \) for some infinite \( A \subseteq A_S \).
Since if \( S' \subseteq \N \) is different from \( S \) then \( A_S \cap A_{S'}\) is finite, there is \( \bar \jmath \in A \) such that \( \bar \jmath \notin A_{S'} \). If \( \arc{B'_m}{g_m} \precsim_{\WAr} \arc{B_{S'}}{f_{S'}} \) then
\( \arc{B''}{g_m \cap \bar B''} \precsim_{\WAr} \arc{B_{S'}}{f_{S'}} \), 
and thus \( \arc{B_{\bar \jmath}^p}{f_{\bar \jmath}^p} \precsim_{\WAr} \arc{B_{S'}}{f_{S'}} \), which is impossible by Lemma~\ref{lem:fewpredecessors}\ref{lem:fewpredecessors-1} and the choice of \( \bar \jmath \).
Thus for every \( m < \kappa\) there is at most one \( S \subseteq \N \) such that \( \arc{B'_m}{g_m} \precsim_{\WAr} \arc{B_S}{f_S} \) and we are done.
\end{proof}

By (the proof of) Proposition~\ref{prop:basisforcvx}, Lemma~\ref{lem:froarcstoLin} and Remark~\ref{rmk:froarcstoLin} (together with \( L_S \cong (L_S)^* \)), if \( S,S' \subseteq \N \) are distinct then \( \arc{B_S}{f_S} \) and \( \arc{B_{S'}}{f_{S'}} \) are \( \precsim_{\WAr} \)-incomparable. Thus \( \mathcal{A} \cup \mathcal{B} \) is a \( \precsim_{\WAr} \)-antichain of size \( 2^{\aleph_0} \) containing \( \mathcal{A} \), as desired. 
\end{proof}

\subsection{Knots and their classification}

In the proof of \cite[Theorem 3.1]{Kul17}, a function from $\LO$ to $\Kn$
is defined, call it here $G$, by setting $G(L) = K_{F(\mathbf{1}+L+\mathbf{2}+\eta)}$, where $F$ is the reduction from~\eqref{eq:FfromLOtoA}. 
It was claimed that $G$ was a reduction of $\iso$ to $\equiv_{\Kn}$, but this is not the case. Indeed, notice that if $M$ is a linear order, then  we have $G(\eta+\mathbf{1}+M) \equiv_{\Kn} G(M)$, essentially because  $C[\mathbf{1}+\eta+\mathbf{1}+M+\mathbf{2}+\eta] \cong_{\CO} C[\mathbf{1}+M+\mathbf{2}+\eta]$; however if $M$ is scattered (and in many other cases) $\eta+\mathbf{1}+M \ncong_{\LO} M$. 

One can easily fix this problem by replacing 
\( K_{F(\mathbf{1}+L+\mathbf{2}+\eta)} \) with \( K_{F(\mathbf{1}+L+\mathbf{2}+\eta)+\bigoplus_{\N} \, \arc{B^*}{f^*}} \), where \( \arc{B^*}{f^*} \) is a figure-eight arc. More precisely, we can derive \cite[Theorem 3.1]{Kul17} from (the proof of) Theorem~\ref{thm:reduc_lo_arcs} connecting the endpoints of each arc $F(L)$ with $\bigoplus_{\N} \, \arc{B^*}{f^*}$ and get:

\begin{cor}\label{cor4}
${\cong_{\LO}} \leq_B {\equiv_{\Kn}}$.
\end{cor}

The next result follows from
Theorem~\ref{thm:iso_co_bired_iso_lo} and Corollary~\ref{cor4}.
However, exploiting the obvious analogy between circular orders and knots one obtains a direct and more natural proof.
(The reduction \( F_{\Kn}\) will be used also in the proof of Theorem~\ref{thm:lowerforknots}).

\begin{thm}\label{thm:reduc_co_knots}
${\cong_{\CO}} \leq_B {\equiv_{\Kn}}$.
\end{thm}
\begin{proof}
We define a Borel reduction $F_{\Kn} \colon \CO \to \Kn$ similar to the reduction of the proof of Theorem~\ref{thm:reduc_lo_arcs}. 
Instead of embedding a linear order $L \in \LO$ into $[-1,1] \subseteq \R$, we embed $C \in \CO$ into $S^1 = \R \cup \{\infty\}$ by defining a sequence of intervals $(h_C(n)-2r_C(n),h_C(n)+2r_C(n))_{n \in \N}$ of \( \R \) denoted by $V_n^C$, satisfying conditions analogous to \ref{Vpwdisj}--\ref{dense} of the proof of Theorem \ref{thm:reduc_lo_arcs}.

As before, for every $n \in \N$ let $U_n^C = [h_C(n)-r_C(n),h_C(n)+r_C(n)]$, consider $\bar B_n^C = U_n^C \times [-r_C(n),r_C(n)]^2$ and define a proper arc $(\bar B_n^C,f_n^C)$ as in Figure \ref{trefoils}. 
Set $f^C=\{(x,0,0) \mid (x,0,0) \notin \bigcup_{n \in \N} \bar B_n^C\} \cup \{\infty\}$. 
Finally we consider the knot $F_{\Kn}(C)$ given by $\bigcup_{n \in \N} f_n^C \cup f^C$. 
The rest of the proof is an adaptation of the proof of Theorem \ref{thm:reduc_lo_arcs} to this case.
\end{proof}

\begin{remark}
Theorem 4.1 of \cite{Kul17} shows that a certain equivalence relation induced by a turbulent action is Borel reducible to ${\equiv_{\Kn}}$.
Therefore, since $\cong_{\LO}$ and $\cong_{\CO}$ are induced by actions of $S_\infty$ the reductions in Corollary~\ref{cor4} and Theorem~\ref{thm:reduc_co_knots} are actually strict by Theorem \ref{turb_act}.
\end{remark}

In order to extend to knots the analysis of $\precsim_{\Ar}$ previously developed, 
one may be tempted to transfer the subarc relation from proper arcs to knots. To this aim, we first introduce the following notion.

\begin{defn}
    Let $\{\arc{B_i}{f_i} \mid i \leq n\}$ be a collection of oriented proper arcs, $h\colon T \to S^3$ a topological embedding, and $K \in \Kn$.
    We say that $K$ is the \textbf{circular sum of} the $\arc{B_i}{f_i}$'s via $h$, if $K=C^h[\bigoplus_{i\leq n}\arc{B_i}{f_i}]$ (recall Definitions \ref{def:circularization_of_arcs} and \ref{def:sumarc}).
\end{defn}

\begin{remark}\label{rem:tame_equiv_trivial}
A topological embedding $h$ of $T$ in $S^3$ is canonical if the closure of $S^3 \setminus h(T)$ is a solid torus as well (recall the solid torus theorem, see e.g.\ \cite[p.\ 107]{Ro90}). 
For any $\arc{B}{f} \in \Ar$, the knot $K_{\arc{B}{f}}$ introduced after Definition~\ref{def:arc} is equivalent to $C^h[\arc{B}{f}\oplus I_{\Ar}]$ for any canonical $h$.
Moreover, when $\arc{B}{f}$ is tame we have $K_{\arc{B}{f}} \equiv_{\Kn} C^h[\arc{B}{f}]$ for every such $h$, i.e.\ the two operations of joining the endpoints of $\arc{B}{f}$ with a trivial arc and of circularization of $\arc{B}{f}$ yield the same knot (up to equivalence).   
\end{remark}

Intuitively, a knot \( K \) is a subknot of \( K' \) when \( K' \) can be split into two subarcs, one of which (when circularized) is equivalent to $K$. Obviously, when splitting \( K' \) we want that its singularities are still singularities of its subarcs, and vice versa. More formally, let \( h_0 \colon T \to S^3 \) be the canonical embedding.

\begin{quote}
Given two knots $K, K' \in \Kn$, we say that $K$ is a \textbf{subknot} of $K'$ if either%
\footnote{This condition is added to have at least a reflexive relation.}
\( K \equiv_{\Kn} K' \), or there are proper arcs \( \arc{B_0}{f_0}, \arc{B_1}{f_1} \in \Ar \) such that 
\( K' \equiv_{\Kn} C^{h_0}[\arc{B_0}{f_0} \oplus \arc{B_1}{f_1}] \) and \( K \equiv_{\Kn} C^{h_0}[\arc{B_0}{f_0}] \).
\end{quote}
It is easy to see that for tame knots this is equivalent to: $K$ is a subknot of $K'$ if either $K \equiv_{\Kn} K'$ or $K \equiv_{\Kn} K_{\arc{B'}{K' \cap \bar{B'}}}$ 
for some subarc \( \arc{B'}{K' \cap \bar{B}'} \) of \( K' \) (see Remark \ref{rem:tame_equiv_trivial}).
However, as in the case of convex embeddability for circular orders, this relation is not transitive on $\Kn$ (it is however transitive on the class of tame knots).  
Since we want to study wild knots, to obtain a transitive relation we define a piecewise version of the subknot relation, which is the analogue for knots of $\pccvx$ (recall Definition \ref{cvx_co}). 
Furthermore, in the following definition we allow arbitrary circularizations: in this way, we work modulo tame variations of the knots at hand (see Remark~\ref{rem:tame_part_irrelevant}), so that the anti-classification results we are going to prove are actually stronger than if we allowed only canonical embeddings of \( T \) into \( S^3 \).




\begin{defn}\label{defn0 4}
Let $K, K'\in \Kn$. Then $K$ is a \textbf{(finite)
piecewise subknot} of $K'$, in symbols
\[
K \skr K',
\] 
if and only if either $K \equiv_{\Kn} K'$, or \( K' \) can be written as the circular sum of finitely many oriented proper arcs \( \{ \arc{B_j}{f_j} \mid j \leq k' \} \) (via some embedding \( h' \colon T \to S^3 \)) in such a way that \( K \equiv_{\Kn} C^h[\bigoplus_{i \leq k} \arc{B_{j_i}}{f_{j_i}}] \) for some (increasingly ordered) subsequence \(  j_0, \dotsc, j_k  \) of \( 0, \dotsc , k' \) and some embedding \( h \colon T \to S^3 \).

The \textbf{(finite) piecewise mutual subknot relation} is the relation defined by $K \approx_{\Kn}^{<\o} K'$ if and only if $K \skr K'$ and $K' \skr K$.
\end{defn}

\begin{prop}\label{prop:subarc_transitivity}  
$\skr$ and $\approx_{\Kn}^{<\o}$ are an analytic quasi-order and an analytic equivalence relation on $\Kn$, respectively.
\end{prop}
\begin{proof}
It is easy to see that $\skr$ is reflexive and analytic.
To prove transitivity we can mostly mimic the proof of Proposition \ref{prop:pccvx_is_trans}.
\end{proof}

The quasi-order $\skr$ is fine enough to distinguish between tame and wild knots, as shown in the next proposition.

\begin{prop}\label{prop:tame_knots_for_star}
Let $K \in \Kn$, and recall that we denote by \( I_{\Kn } \) the trivial knot. Then the
following are equivalent:
\begin{enumerate-(1)}
\item \( K \) is tame;
\item $K \approx_{\Kn}^{<\o} I_{\Kn}$;
\item \( K \skr I_{\Kn} \).
\end{enumerate-(1)}
In particular, the \( \approx^{<\o}_{\Kn} \)-class of
the tame knots is minimal with respect to (the quotient order of)
\( \skr \).
\end{prop}
\begin{proof}
The proof is immediate using the facts that a knot is tame if and only if it is a circularization of the trivial arc, and that the trivial knot can be written as a circular sum only if all summands are trivial arcs and the embedding of the solid torus is canonical.
\end{proof}

\begin{remark}\label{rem:tame_part_irrelevant}
If $\arc{B}{f}$ is a proper arc with a tame subarc and $\arc{B'}{g}$ is a tame arc then it is easy to check that $C^h[\arc{B}{f}] \approx^{<\o}_{\Kn} C^{h'}[\arc{B}{f} \oplus \arc{B'}{g}]$ for any topological embeddings $h$ and $h'$.   
\end{remark}

Notice that the relations $\skr$ and $\approx^{<\o}_{\Kn}$ differ from $\equiv_{\Kn}$ only on the set of knots which are circularizations of proper arcs. 
For this reason we focus on the following subset of $\Kn$.

\begin{defn}\label{def:ckn}
   We denote by $\CKn$ and $\WCKn$, respectively, the set of knots which are a circularization of a proper arc (that is, up to knot equivalence, those of the form \( C^h[\arc{B}{f}]\) for some \( \arc{B}{f} \in \Ar \) and some embedding \( h \colon T \to S^3\)), and its subset consisting of wild knots.
   Let $\cskr$ and $\precsim_{\WCKn}^{<\o}$ be the restrictions of $\skr$ to these sets.
\end{defn}

Notice that $\CKn$ is a proper subset of $\Kn$: for example, the knot constructed by Bing in \cite{Bi56} cannot be \lq\lq cut\rq\rq\ at any point and thus it does not belong to $\CKn$.
However $\CKn$ is quite rich, as it includes any wild knot $K$ satisfying any of the following equivalent conditions: $K$ has at least one isolated singularity (i.e.\ \( I \Sigma_K \neq \emptyset\)); the set \( \Sigma_K \) of singularities of $K$ is not dense in \( K \); there exists a point of $K$ which is not a singularity (i.e.\ \( \Sigma_K \neq K \)). Moreover, the wild knots built by Artin and Fox in \cite{FA48} do not satisfy the previous conditions, yet they belong to $\CKn$. 
Further evidence of the complexity and richness of $\cskr$ is provided in the results below (see Proposition \ref{prop:chains_skr} and Theorems \ref{thm:unbounded_and_dominating_for_knots}--\ref{thm:antichains_for_knots}).

Since $C^h[\arc{B}{f}] \approx_{\Kn}^{<\o} C^{h'}[\arc{B}{f}]$ for any topological embeddings $h$ and $h'$, every $K \in \CKn$ can be assumed to be, up to $\approx_{\CKn}^{<\o}$, of the form $C^h[\arc{B}{f}]$ for some canonical embedding $h\colon T \to S^3$. 
To simplify the notation we write $C[\arc{B}{f}]$ in place of $C^h[\arc{B}{f}]$ when $h$ is canonical and we do not mention $h$ and $h'$ witnessing $K \skr K'$ when they are canonical, which can always be assumed to be the case.

The next theorem establishes a lower bound for the complexity of $\cskr$ with respect to Borel reducibility.

\begin{thm}\label{thm:lowerforknots}
${\pccvx} 
\leq_B
{\cskr}$.
\end{thm}

\begin{proof}
We claim that the Borel map
\( F_{\Kn} \colon \CO \to \Kn\)
from the proof of Theorem~\ref{thm:reduc_co_knots} is the desired reduction. 
First of all, notice that by construction $F_{\Kn}(C) \in \CKn$ for every \( C \in \CO \). 
Fix now \( C, C' \in \CO \).

Assume first that ${C}\pccvx {C'}$, and let the finite convex partition $(C_i)_{i \leq k}$ of $C$ and the embedding $g$ witness this.
As usual, we can assume $k>0$.
For every $i\leq k$, let $C'_{2i}=g(C_i)$, for every $i<k$ let $C'_{2i+1}$ be the convex subset of $C'$ given by the elements $n$ such that $C'(\ell, n, \ell')$ for every $\ell \in C'_{2i}$ and $\ell' \in C'_{2i+2}$; moreover, let $C'_{2k+1} \csube C'$ be the set of the $n$'s such that $C'(\ell, n, \ell')$ for every $\ell \in C'_{2k}$ and $\ell' \in C'_0$ (if any of these sets is empty just delete it and reindex the remaining convex sets appropriately). 
For every $j\leq 2k+1$, set $\bar B_j= [a_j,b_j] \times [-1,1]^2$, where $a_j=\inf \bigcup_{n \in C'_j}V^{C'}_n$ and $b_j=\sup\bigcup_{n \in C'_j}V^{C'}_n$, and \( f_j = F_{\Kn}(C') \cap \bar B_j\), so that \( (\bar B_j, f_j)\) is a proper arc.
Notice now that $F_{\Kn}(C')=C[\bigoplus_{j\leq 2k+1}(\bar B_j, f_j)]$. Moreover, since \( C \cong_{\CO} \sum_{i \leq k} C'_{2i} \) we have $F_{\Kn}(C)\equiv_{\Kn} C[\bigoplus_{i\leq k}(\bar B_{2i}, f_{2i})]$.
Hence, $F_{\Kn}(C) \skr F_{\Kn}(C')$.

Conversely, suppose that $F_{\Kn}(C)$ and $F_{\Kn}(C')$ (which are elements of $\CKn$) are such that $F_{\Kn}(C) \skr F_{\Kn}(C')$, and let $\{(\bar{B}_i,f_i) \mid i \leq k'\}$ and the subsequence $(j_i)_{i \leq k}$ of $0,\dots,k'$ witness this.
By definition of $F_{\Kn}(C')$, when $\bar B_i \cap \bar B_m$ contains a point $x \in I\Sigma_{F_{\Kn}(C')}$ then $x$ is a singular point of only one of $(\bar{B}_i,f_i)$ and $(\bar{B}_m,f_m)$; by reindexing the sequence $\{(\bar{B}_i,f_i) \mid i \leq k'\}$ we can assume this occurs always for the index which is the immediate predecessor of the other one in $C[\mathbf{k'+1}]$. Since $F_{\Kn}(C) \equiv_{\Kn} C[\bigoplus_{i\leq k}(\bar B_{j_i}, f_{j_i})]$, by an analogue of \eqref{eq:not-inverting} we have $I\Sigma_{F_{\Kn}(C)} \cong_{\CO} I\Sigma_{C[\bigoplus_{i\leq k}(\bar B_{j_i}, f_{j_i})]}$ via some map $g$.
  
Recall that $h_C$ is an isomorphism of circular orders between $C$ and $I\Sigma_{F_{\Kn}(C)}$, so that $g \circ h_C$ is an isomorphism of circular orders as well, and let $C_i = (g \circ h_C)^{-1} (I\Sigma_{F_{\Kn}(C')}\cap\bar{B}_{j_i} \setminus \bar{B}_{j_{i+1}})$.
Notice that $(C_i)_{i \leq k}$ is a finite convex partition of \( C\).  
Moreover, since $\arc{B_{j_i}}{f_{j_i}}$ is a subarc of $F_{\Kn}(C')$ for every $i\leq k$, we have that each $I\Sigma_{F_{\Kn}(C')} \cap \bar{B}'_{j_i}\setminus \bar{B'_m}$ (for $m$ the immediate predecessor of $j_i$ in $C[\mathbf{k'+1}]$) is convex in \( I\Sigma_{F_{\Kn}(C')} \). 
Finally, since $I\Sigma_{F_{\Kn}(C')}\cong C'$ via $h_{C'}^{-1}$, then ${C} \pccvx C'$, as desired.
\end{proof}

\begin{cor} \label{cor:above}
${\pccvxeq} \leq_B {\approx_{\CKn}^{<\omega}}$,
hence also ${\iso}\leq_B {\approx_{\CKn}^{<\omega}}$ and
${E_1} \leq_B {\approx_{\CKn}^{<\omega}}$.
\end{cor}

The fact that the isomorphism on linear orders is Borel reducible to $\approx_{\CKn}^{<\omega}$ implies that $\approx_{\CKn}^{<\omega}$ is proper analytic. 
Moreover, \({\approx_{\CKn}^{<\omega}}\) is not Baire reducible to an orbit equivalence relation because it Borel reduces \( E_1 \), in stark contrast
with knot equivalence \( \equiv_{\Kn} \); in particular we have that
\({\approx_{\CKn}^{<\omega}}\) is not Borel, or even Baire, reducible to
$\equiv_{\Kn}$.

Using Theorem \ref{thm:lowerforknots}, we can transfer the combinatorial properties of $\pccvx$ proved in Proposition \ref{prop:antichainsCO} to $\cskr$.

\begin{prop}\label{prop:chains_skr}
\begin{enumerate-(a)}
\item
There is an embedding from the partial order
\((\Cvx(\R),\subseteq)\) into $\cskr$, and indeed \( {(\Cvx(\R),\subseteq)} \leq_B {\cskr} \).
\item
$\cskr$ has chains of order type \( (\R, {<} ) \), as well as antichains of size $2^{\aleph_0}$.
\end{enumerate-(a)}
\end{prop}

To extend the other combinatorial properties of $\pccvx$ to  $\cskr$ we need an analogous of Lemma \ref{lem:froarcstoLin}.
When \( K \) is a knot and $f$ is such that $\Ima f =K$, the set \( I\Sigma_{K}\) can naturally be viewed as a circular order \( C^K_f = (I\Sigma_{K}, C_f) \). 
As it was the case for proper arcs, the set \( I\Sigma_{K}\) is (at most) countable and thus \( C^K_f \) is either a finite or a countable circular order.  
If \( f,f'\colon S^1 \to S^3 \) are topological embeddings giving rise to the same knot, then either \( C^K_f = C^K_{f'} \) or \( C^K_f = (C^K_{f'})^* \).
Recall that by construction, for knots of the form%
\footnote{If not specified otherwise, we always choose the natural orientation of \( F_{\Kn}(C) \), witnessed by $f$.}
\( F_{\Kn}(C)\)  we have \( C^{F_{\Kn}(C)}_{f} \cong_{\CO} C \). 

\begin{lem} \label{lem:from_knots_to_CO}
Let \( K, K' \in \CKn \) be such that \( K\skr K' \) and let $f$ and $f'$ be such that $\Ima f =K$ and $\Ima f'=K'$.
Then there exists a finite set $A \subseteq I\Sigma_K$ such that either \( C^K_f \setminus A \pcvx C^{K'}_{f'}\) or \( C^K_f \setminus A \pcvx {(C^{K'}_{f'})^*}\).
\end{lem}

\begin{proof} 
Let $\{\arc{B_i}{f_i} \mid i \leq k'\}$ and the subsequence $(j_i)_{i \leq k}$ of $0,\dots,k'$ witness $K \skr K'$, and denote by $g$ the homeomorphism witnessing $K \equiv_{\Kn} C[\bigoplus_{i \leq k}\arc{B_{j_i}}{f_{j_i}}]$.
Let $A = \{x \in I\Sigma_K \mid \exists i \leq k \, (g(x) \in \partial \bar{B}_{j_i})\}$, and notice that \( A \) contains at most $k+1$ points.
We can assume that $g$ agrees with the orientations induced on $K$ and $K'$ by $f$ and $f'$, in which case we show that \( C^K_f  \setminus A \pcvx C^{K'}_{f'}\) (if $g$ agrees with only one of the orientations we obtain \( C^K_f  \setminus A \pcvx (C^{K'}_{f'})^*\), and if it disagrees with both it suffices to reverse both orientations).

For every $i \leq k$ let $C_i = I\Sigma_K \cap g^{-1}(\Int \bar{B}_{j_i})$. 
Then each $C_i$ is convex, and \( \{ C_i \mid i \leq k \} \) is a finite convex partition of \(C^K_f \setminus A\) (some of the $C_i$'s might actually be empty, in which case we would obtain a convex partition with less than $k+1$ sets, but for notational ease we avoid keeping track of this). Since each $\arc{B_{j_i}}{f_{j_i}}$ is a subarc of $K'$, it is easy to check that $g \rest I\Sigma_K$ is an embedding of circular orders and that \( g(C_i) \csube C^{K'}_{f'} \) for all \( i \leq k \).
\end{proof}

Since in $\Kn \setminus \CKn$ the relation $\skr$ is $\equiv_{\Kn}$, it is easy to show that \( \mathfrak{b} (\skr) = 2\) and  \( \mathfrak{d} (\skr) = 2^{\aleph_0}\).
It is therefore more interesting to compute the unbounding and dominating number of $\cskr$.  

\begin{thm} \label{thm:unbounded_and_dominating_for_knots}
\( \mathfrak{b} (\cskr) = \aleph_1\) and 
\( \mathfrak{d} (\cskr) = 2^{\aleph_0}\).
\end{thm}

\begin{proof} 
We first show the existence of an \( \cskr\)-unbounded family of knots of size \( \aleph_1 \). 
Let \( F_{\Kn} \colon \CO \to \Kn \) be the map defined in the proof of Theorem \ref{thm:reduc_co_knots} and used also in the proof of Theorem \ref{thm:lowerforknots}.
By (the proof of) Proposition \ref{prop:boundingCO} there exists a strictly increasing sequence \( \{ C_\a \mid \a < \o_1 \} \subseteq \CO\) without upper bound with respect to $\pccvx$.
Notice moreover that each $C_\a$ has the property that $C_\a \cong_{\CO} C_\a \setminus A$ for any finite $A \subseteq C_\a$.
The sequence \( \{ F_{\Kn}(C_\a) \mid \a < \o_1 \} \subseteq \CKn\) is then strictly $\cskr$-increasing, and we claim that it is also unbounded in $\CKn$.
Suppose towards a contradiction that there is \( K \in \Kn \) such that \( F_{\Kn}(C_\a) \cskr K\) for all \( \a < \omega_1\). Then \( I\Sigma_{K} \) is infinite and thus the circular order \( C_f^K \) is, up to isomorphism, an element of \( \CO \). 
Pick now $\ell \in  C_f^K$ and define $L \in \LO$ by setting $x \leq_L y$ if and only if $C_f^K(\ell,x,y)$ and $x=\ell$ when $y=\ell$. 
Notice that $C_f^K = C[L]$.
Then the circular order $C=C[L+L^*]$ is such that $C_f^K \ccvx C$ and $(C_f^K)^* \ccvx C$. 
By Lemma~\ref{lem:from_knots_to_CO} and the fact that each \( C^{F_{\Kn}(C_\alpha)}_{f} \cong_{\CO} C_\alpha \), this would imply that for every $\a < \o_1$ there exists a finite $A_\a \subseteq C_\a$ such that $C_\a \setminus A_\a \pccvx C$.
As $C_\a \cong_{\CO} C_\a \setminus A$, the circular order \( C\) would be a \( \pccvx\)-upper bound for \( \{ C_\a \mid \a < \o_1 \}\), yielding the desired contradiction.

We now prove that \( \mathfrak{b} (\cskr) > \aleph_0 \). Let \( \{ K_i \mid i \in \N \} \subseteq \CKn \) be a countable family of knots. By definition of $\CKn$, each $K_i$ can be written as $C[\arc{B_i}{f_i}]$ for some proper arc $\arc{B_i}{f_i}$ (we are using a canonical embedding of the solid torus in $S^3$).
Then the knot \( C[\bigoplus_\N \arc{B_i}{f_i}] \) is a \( \cskr\)-upper bound for \( \{ K_i \mid i \in \N \} \).

To prove that $\mathfrak{d} ({\cskr}) \geq 2^{\aleph_0}$ we follow the same strategy of the proof of Theorem \ref{thm:umboundedanddominatingforarcs}. 
Consider the \( \pccvx \)-antichain \(  \{ C_S \mid S \subseteq \N \} \) defined in the proof of Proposition \ref{prop:basisCO}\ref{prop:basisCO-1} and, using  Lemma~\ref{lem:from_knots_to_CO} (removing finitely many elements from $C_S$ does not affect the argument) and the proof of Proposition \ref{prop:antichains2CO}, prove that for every knot \( K \in \CKn\), \(\{ F_{\Kn}(C_S) \mid  {F_{\Kn}(C_S) \cskr K} \}\) is countable.
The proof is then completed using Theorem~\ref{thm:reduc_co_knots}.
\end{proof}

\begin{cor} \label{cor:incomparable_knots_and_unbounded_chains}
Every knot \(K \in \CKn \) is the bottom of an \( \cskr\)-unbounded chain of length~\(\omega_1\).
\end{cor}

\begin{proof}
Given $K \in \CKn$, let $\arc{B}{f}$ be a proper arc such that $K=C[\arc{B}{f}]$.
As in the proof of Theorem~\ref{thm:unbounded_and_dominating_for_knots} let \( \{ C_\a \mid \a < \o_1 \} \subseteq \CO\) be an unbounded strictly $\pccvx$-increasing sequence in $\CO$, so that \( \{ F_{\Kn}(C_\a) \mid \a < \o_1 \} \subseteq \CKn\) is unbounded and strictly $\cskr$-increasing in $\CKn$.
For every $\a<\o_1$, let $\arc{B_\a}{f_\a}$ be a proper arc obtained by cutting $F_{\Kn}(C_\a)$ in a point which is not an isolated singularity, so that in particular $F_{\Kn}(C_\a) = C[\arc{B_\a}{f_\a}]$.
Let \( K_0 = K\) and, for $0<\a<\o_1$, \( K_\a = C[{\arc{B}{f} \oplus \arc{B_\alpha}{f_\alpha}}] \).
For every \( \alpha < \beta < \omega_1\) we have \( K_\alpha \cskr K_\beta\) (even though it might happen that $\arc{B_\a}{f_\a} \not\precsim_{\Ar} \arc{B_\beta}{f_\beta}$, in which case we need a circular sum of proper arcs with more than one element to witness \( K_\alpha \cskr K_\beta\))  and \( F_{\Kn}(C_\a) \cskr K_\a \).
Hence the sequence $(K_\a)_{\a<\o_1}$ is \( \cskr\)-unbounded. 
By the same argument used in the proof of Corollary \ref{cor:incomparablearcsandunboundedchains} we can extract from $(K_\a)_{\a<\o_1}$ a strictly \(\cskr\)-increasing subsequence of length \( \omega_1 \) starting with \( K \).
\end{proof}

We finally deal with minimal elements and bases with respect to $\cskr$. 
In contrast with the case of proper arcs, it is not interesting to consider the restriction of $\cskr$ to the collection of tame knots because by Proposition \ref{prop:tame_knots_for_star} tame knots are all $\approx_{\CKn}^{<\omega}$-equivalent.
Thus let us consider $\precsim_{\WCKn}^{<\o}$.

\begin{lem}\label{lem:knots_sum_prime_arcs}
Let \( \{ \arc{B^p_i}{f^p_i} \mid i \in \N \} \) be a family of (oriented) prime arcs, and let $K^*_S = C[\bigoplus_{i \in S} \, \arc{B^p_i}{f^p_i}]$ for some infinite $S \subseteq \N$. 
\begin{enumerate-(a)}
\item \label{knots_sum_prime_arcs-0}
If \( K^*_S = C^h[\bigoplus_{i \leq k} \arc{B_i}{f_i}] \) for some \( k  \in \N \) and \( h \colon T \to S^3 \), then there is a unique \( j \leq k \) such that \( \arc{B_j}{f_j}\) is wild; moreover, either \( \arc{B_j}{f_j} \equiv_{\Ar} \bigoplus_{i \in S'} \arc{B^p_i}{f^p_i} \) or \( \arc{B_j}{f_j} \equiv_{\Ar} \bigoplus^\partial_{i \in S'} \arc{B^p_i}{f^p_i} \) for some \( S' \subseteq S \) with \( S \setminus S' \) finite. 
\item\label{knots_sum_prime_arcs-1}
The knot $K^*_S$ is $ \precsim_{\WCKn}^{<\o}$-minimal in $\WCKn$.
\item\label{knots_sum_prime_arcs-2}
If $K^*_{S_0} \precsim_{\WCKn}^{<\o} K^*_{S_1}$ then $S_0 =^* S_1$, where $=^*$ is the identity modulo a finite set.
\end{enumerate-(a)} 
\end{lem}

\begin{proof}
\ref{knots_sum_prime_arcs-0}
Let \( j \leq k \) be such that \( \barc{h(\bar B_j)}{h \circ f_j} \) is wild and contains the unique singularity \( x \) of \( K^*_S \). There is at least one such \( j \) because otherwise \( K^*_S \) would be tame, and it is unique because the  singularity \( x \) is ``one-sided'', i.e.\ it is witnessed only on one side while the other side is tame. Thus \( \arc{B_j}{f_j} \), being equivalent to \( \barc{h(\bar B_j)}{h \circ f_j} \) via \( h \restriction \bar B_j \), is wild, while all other \( \arc{B_i}{f_i}\) with \( i \neq j \) are tame because so are the proper arcs \( \barc{h(\bar B_i)}{h \circ f_i} \). Moreover, by construction \( \barc{h(\bar B_j)}{h \circ f_j} \) is either of the form \( \bigoplus_{i \in S'} \arc{B^p_i}{f^p_i} \) (if \( x \in \Int h(\bar B_j)\)) or \( \bigoplus^\partial_{i \in S'} \arc{B^p_i}{f^p_i} \) (if \( x \in \partial \, h(\bar B_j) \)), for some \( S' \subseteq S \) omitting finitely many elements of \( S \): since \( \arc{B_j}{f_j} \equiv_{\Ar} \barc{h(\bar B_j)}{h \circ f_j} \) we are done.

\ref{knots_sum_prime_arcs-1} 
Suppose that $K \in \WCKn$ is such that $K\precsim_{\WCKn}^{<\o} K^*_S$ but \( K \not\equiv_{\Kn} K^*_S \) (otherwise we are done), and let $\{ \arc{B_j}{f_j} \mid j \leq k'\}$ and $(j_i)_{i \leq k}$ witness this. 
By part~\ref{knots_sum_prime_arcs-0} there is a unique \( \ell \) such that \( \arc{B_\ell}{f_\ell} \) is wild, and necessarily \( \ell \) is one of the $j_i$'s because otherwise \( K \) would be tame. Using Remark~\ref{rem:tame_part_irrelevant} one easily gets
\[  
K \equiv_{\Kn} C[\bigoplus_{i \leq k}\arc{B_{j_i}}{f_{j_i}}] \approx_{\WCKn}^{<\o} C[\arc{B_\ell}{f_\ell}] \approx_{\WCKn}^{<\o} K^*_S.
\]

\ref{knots_sum_prime_arcs-2} 
Let $\{ \arc{B_j}{f_j} \mid j \leq k'\}$ and the subsequence $(j_i)_{i \leq k}$ of $0, \dots, k'$ witness $K^*_{S_0} \precsim_{\WCKn}^{<\o} K^*_{S_1}$.
Apply part~\ref{knots_sum_prime_arcs-0} to \( K^*_{S_1}\) to isolate the unique \( \ell \leq k' \) such that the proper arc \( \arc{B_\ell}{f_\ell} \) is wild. 
Let also \( S'_1 \subseteq S_1 \) be such that 
\(\arc{B_\ell}{f_\ell} \equiv_{\Ar} \bigoplus\nolimits^{(\partial)}_{i \in S'_1} \arc{B^p_i}{f^p_i} \), with \( S_1 \setminus S'_1 \) is finite. 
Arguing as in \ref{knots_sum_prime_arcs-1}, we have that $\ell$ is one of the $j_i$'s such that $K^*_{S_0}
\equiv_{\Kn} C[\bigoplus_{i \leq k} \arc{B_{j_i}}{f_{j_i}}]$. Then there is a subarc $\arc{B^*}{f^*}$ of $K$ such that $\arc{B^*}{f^*} \equiv_{\Ar} \arc{B_\ell}{f_\ell}$, and since it is wild, by applying \ref{knots_sum_prime_arcs-0} we obtain that $\arc{B^*}{f^*}$ is of the form $\bigoplus\nolimits^{(\partial)}_{i \in S'_0} \arc{B^p_i}{f^p_i}$ with \( S_0 \setminus S'_0 \) finite. Thus, 
\( \bigoplus^{(\partial)}_{i \in S'_0} \arc{B^p_i}{f^p_i} \equiv_{\Ar} \bigoplus^{(\partial)}_{i \in S'_1} \arc{B^p_i}{f^p_i}\), whence it follows that \( S'_0 = S'_1 \). Hence \( S_0 =^* S_1 \).
\end{proof}

\begin{thm}\label{thm:basis_for_knots}
\begin{enumerate-(a)}
    \item\label{thm:basis_for_knots-1} There are \( 2^{\aleph_0}\)-many \( \precsim_{\WCKn}^{<\o}\)-incomparable \( \precsim_{\WCKn}^{<\o} \)-minimal elements in \( \WCKn \). In particular, all bases for \( \precsim_{\WCKn}^{<\o} \) are of maximal size.
    \item\label{thm:basis_for_knots-2}
    There is a strictly \( \precsim_{\WCKn}^{<\o} \)-decreasing \( \omega \)-sequence in \( \WCKn \) which is not \( \precsim_{\WCKn}^{<\o}\)-bounded from below.
\end{enumerate-(a)} 
\end{thm}

\begin{proof}
Fix an enumeration without repetitions $\{\arc{B_i^p}{f_i^p} \mid i \in \N\}$ of all prime arcs.
    
\ref{thm:basis_for_knots-1} 
As in the proof of Theorem \ref{thm:basisforarcs}\ref{thm:basisforarcs-3}, let $\mathcal{P}$ be a family of size \( 2^{\aleph_0} \) consisting of infinite subsets of $\N$ with pairwise finite intersections. For every $S \in \mathcal{P}$ consider the knot $K^*_S$ defined in Lemma~\ref{lem:knots_sum_prime_arcs}.
By Lemma \ref{lem:knots_sum_prime_arcs}\ref{knots_sum_prime_arcs-1} each $K^*_S$ is \(\precsim_{\WCKn}^{<\o} \)-minimal in $\WCKn$, and if $S, S' \in \mathcal{P}$ are distinct then $K^*_S$ and $K^*_{S'}$ are \(\precsim_{\WCKn}^{<\o} \)-incomparable by Lemma~\ref{lem:knots_sum_prime_arcs}\ref{knots_sum_prime_arcs-2}.
    
\ref{thm:basis_for_knots-2}
Let $K_n=C[{\bigoplus_{i \geq n} (\bigoplus_{\N}^\partial \, \arc{B_i^p}{f_i^p})}]$ (notice that each $i \geq n$ is associated to an element of $I\Sigma_{K_n}$). 
We prove that $\{K_n \mid n \in \N\}$ is the desired $\o$-chain. 
    
Let $n < n'$. 
Clearly, $K_{n'} \precsim_{\WCKn}^{<\o} K_n$. Suppose now, towards a contradiction, that $K_n \precsim_{\WCKn}^{<\o} K_{n'}$, as witnessed by $\{(\bar{B}_j,f_j) \mid j \leq k'\}$ and the subsequence $(j_i)_{i \leq k}$ of $0,\dots,k'$.  
Since $K \equiv_{\Kn} C[\bigoplus_{i \leq k} \arc{B_{j_i}}{f_{j_i}}]$,  there exist $i \leq k$ and $m \in \N$ such that $\arc{B_{j_i}}{f_{j_i}}$ contains (a proper arc equivalent to) the tail $\bigoplus_{t \geq m}^\partial\, \arc{B_n^p}{f_n^p}$ of $\bigoplus_{\N}^\partial\, \arc{B_n^p}{f_n^p}$. But $\arc{B_{j_i}}{f_{j_i}}$ is a subarc of $K'_{n'}$, a contradiction. 
Hence $K_{n'} \prec_{\WCKn}^{<\o} K_n$.

Suppose now that $K \in \WCKn$ bounds from below $\{K_n \mid n \in \N\}$. 
Notice that $K \precsim_{\WCKn}^{<\o} K_0$ implies $I\Sigma_K \neq \emptyset$, so that we can fix $x \in I\Sigma_K$. Let $\{(\bar{B}_j,f_j) \mid j \leq k'\}$, and the subsequence $(j_i)_{i \leq k}$ of $0,\dots,k'$ witness $K \precsim_{\WCKn}^{<\o} K_0$. Let $g$ be the homeomorphism witnessing that $K\equiv_{\Kn} C[\bigoplus_{i \leq k} \arc{B_{j_i}}{f_{j_i}}]$. 
Then there exists $i \leq k$ such that $g(x) \in \bar  B_{j_i}$ and $\arc{B_{j_i}}{f_{j_i}}$ is wild and contains an element of $I\Sigma_{K_0}$, which belongs to $\bigoplus_{\N}^\partial\, \arc{B_n^p}{f_n^p}$ for some $n \geq 0$.
This implies that there exists $m \in \N$ such that $C[{\bigoplus^\partial_{t \geq m}\, \arc{B_n^p}{f_n^p}}] \precsim_{\WCKn}^{<\o} K$.
But by the argument of the previous paragraph $C[{\bigoplus^\partial_{t \geq m}\, \arc{B_n^p}{f_n^p}}] \not\precsim_{\WCKn}^{<\o} K_{n+1}$, 
hence $K \not\precsim_{\WCKn}^{<\o} K_{n+1}$, which is a contradiction.
\end{proof}

When $K \in \Kn$, we say that $x \in K$ is \textbf{isolable in $K$} if there exists a subarc $\arc{B}{f}$ of $K$ such that $x \in I\Sigma_{\arc{B}{f}}$. 
Notice that every $x \in I\Sigma_K$ is isolable in $K$, but some point which is isolable in $K$ can fail to belong to $I\Sigma_K$ because it is an accumulation point of other singular points only from one side. 
It is immediate that the set of points isolable in $K$ is countable.

\begin{thm} \label{thm:antichains_for_knots} 
Every \( \precsim_{\WCKn}^{<\o}\)-antichain is contained in a \( \precsim_{\WCKn}^{<\o}\)-antichain of size \( 2^{\aleph_0} \). In particular, there are no maximal \( \precsim_{\WCKn}^{<\o}\)-antichains of size smaller than \( 2^{\aleph_0}\), and every \( K \in \WCKn \) belongs to a \( \precsim_{\WCKn}^{<\o}\)-antichain of size \( 2^{\aleph_0} \).
\end{thm}

\begin{proof}
Let $\{\arc{B_i^p}{f_i^p} \mid i \in \N\}$, $\mathcal{P}$ and $K_S \in \WCKn$, with $S \in\mathcal{P}$, be as in the proof of Theorem \ref{thm:basis_for_knots}\ref{thm:basis_for_knots-1}. 
Following the proof of Proposition \ref{prop:everylobelongstoantichain}, it is enough to prove that the set \( \{ S \in \mathcal{P} \mid K_S \precsim_{\WCKn}^{<\o} K \} \) is countable for each $K \in \WCKn$. 
Suppose that $S \subseteq \N$ is such that $K_S\precsim_{\WCKn}^{<\o} K$, as witnessed by $\{(\bar{B}_j,f_j) \mid j \leq k'\}$ and the subsequence $(j_i)_{i \leq k}$ of $0,\dots,k'$. 
There exist $i \leq k$ and $m$ such that $\bigoplus_{i \in S, i \geq m}^\partial\, \arc{B_i^p}{f_i^p} \precsim_{\Ar} 
\arc{B_{j_i}}{f_{j_i}} $; thus $\bigoplus_{i \in S, i \geq m}^\partial\, \arc{B_i^p}{f_i^p}$ is a subarc of $K$.
Therefore there exists $x_S \in \Sigma_K$ which is determined by a tail of $\bigoplus_{i \in S}^\partial\, \arc{B_i^p}{f_i^p}$. Notice that $x_S$ is isolable in $K$.
If $S,S' \in \mathcal{P}$ are distinct 
and $x_S=x_{S'}$ 
then  by Lemma \ref{lem:fewpredecessors}\ref{lem:fewpredecessors-2} (and recalling that $S$ and $S'$ have finite intersection) the images of $\bigoplus_{i \in S, i \geq m}\, \arc{B_i^p}{f_i^p}$ and $\bigoplus_{i \in S', i \geq m'}\, \arc{B_i^p}{f_i^p}$ approach $x_S$ from opposite sides. 
Hence, $|\{S \in \mathcal{P} \mid x_S=x\}| \leq 2$ for every $x$ isolable in $K$.
Since the set of isolable points in $K$ is countable, \( \{ S \in \mathcal{P} \mid K_S \precsim_{\WCKn}^{<\o} K \} \) is countable as well.
\end{proof}

\section{Open Problems}

We conclude by discussing some open problems that need more investigation. The main one concerns the classification of $\cvxeq$ with respect to Borel reducibility. In Section \ref{sec:complexity of cvx} we showed that ${\iso} \leq_B {\cvxeq}$ and that ${\cvxeq} \leq_{\text{\scriptsize \textit{Baire}}} {\iso}$. 

\begin{question}
    Is $\cvxeq$ Borel reducible to $\iso$ or at least to an equivalence relation induced
    by a continuous action of a Polish group?
\end{question}

In Section \ref{sec:cvx_co} we introduced the quasi-order $\pccvx$ of piecewise convex embeddability for countable circular orders and proved that ${\pccvxeq} \nleq_B {\cvxeq}$, so that also ${\pccvx} \nleq_B {\cvx}$. 

\begin{question}
    Is it true that ${\cvx} <_B{\pccvx}$? What about ${\cvxeq} <_B {\pccvxeq}$?
\end{question}

Considering the combinatorial properties of the quasi-orders defined in the previous sections, many other natural questions can be explored. For example:

\begin{question}
    Is there a basis for $\cvx$ (respectively: $\pccvx$, $\precsim_{\WAr}$, or $\cskr$) which is an antichain?
\end{question}

It is open if the classifications of proper arcs and knots with respect to \textit{equivalence} are related.

\begin{question}
    Is ${\equiv_{\Ar}} \leq_B {\equiv_{\Kn}}$ and/or ${\equiv_{\Kn}} \leq_B {\equiv_{\Ar}}$?
\end{question}

A problem in attacking the second question is the existence of knots which are not equivalent to the circularization of any proper arc.

Recalling that for linear and circular orders the isomorphism relation is Borel reducible to (piecewise) convex embeddability, it is natural to ask the following questions:

\begin{question}
    Does ${\equiv_{\Ar}} \leq_B {\approx_{\Ar}}$? Does ${\equiv_{\Kn}} \leq_B {\approx_{\Kn}^{<\o}}$?
\end{question}

\newcommand{\etalchar}[1]{$^{#1}$}

\end{document}